\documentclass[paper=a4, fontsize=12pt]{article} % A4 paper and 11pt font size   scrartcl
\DeclareOldFontCommand{\rm}{\normalfont\rmfamily}{\mathrm}
\usepackage{graphicx}

\usepackage[margin=1in]{geometry}

\usepackage[english]{babel}
\usepackage{amsfonts}
\usepackage{amsmath}
\usepackage{amsbsy}
\usepackage{amscd}
\usepackage{amsthm}
\usepackage{latexsym} 
\usepackage{float}
\usepackage{setspace}
\usepackage{fancyhdr}
\usepackage{color}
\usepackage{stmaryrd} 
\usepackage{accents}

\newcommand*{\DOT}[1]{%
	\accentset{\mbox{\bfseries .}}{#1}}

\usepackage{bm,algorithm}
\usepackage{lipsum}
\DeclareMathOperator\atanh{atanh}
\usepackage{epstopdf}

\usepackage{subfigure}

\usepackage{tikz}
\usepackage{breakcites}

\newtheorem{remark}{Remark}[section] % add * for  unnumbered remarks 

\newcommand{\bs}[1]{\boldsymbol{#1}}
\newcommand{\uu}[1]  {{\bs{#1}} }

\newcommand{\n}      { \uu{n}         }

\newcommand{\dep}{ {\bs d} }
\newcommand{\vel}{\DOT{\dep}}

\newcommand{\strain}{{{\bs{\epsilon}}}}

\newcommand{\jump}[1]{\llbracket #1 \rrbracket }
\newcommand{\average}[1]{\lbrace\hspace{-1.3mm}\lbrace #1 \rbrace\hspace{-1.3mm}\rbrace }
\newcommand{\diverg}{ \text{div} }

\newcommand{\grad}   { \uu{\nabla}    }

\newcommand{\eps}   {\varepsilon}

\newcommand{\stress}{{\bs\sigma}_{\rm f}}
\newcommand{\sstress}{{\bs\sigma}_{\rm s}}

\newcommand{\bbR}{\mathbb{R}}

\def\eqd{{:=}}

\newcommand{\bphi}{\boldsymbol{\phi}}
\newcommand{\bv}{{\bs v}}

\newcommand{\bx}{{\bs x}}

\newcommand{\bV}{{\bs V}}
\newcommand{\bW}{{\bs W}}

\newcommand{\bd}{{\bs d}}

\newcommand{\bz}{{\bs z}}

\newcommand{\bu}{\uu{u}}

\newcommand{\bw}{{\bs w}}
\newcommand{\bn}{{\bs n}}

\ifpdf
\DeclareGraphicsExtensions{.eps,.pdf,.png,.jpg}
\else
\DeclareGraphicsExtensions{.eps}
\fi

\usepackage{xcolor}

\allowdisplaybreaks

\usepackage{enumitem}

\usepackage{relsize}

\title{	
\LARGE 
A mechanically consistent model for fluid-structure interactions with contact including seepage
}

\author{Erik Burman, Miguel A. Fern\'andez,  Stefan Frei and Fannie M. Gerosa} % Your name

\date{\normalsize\today} % Today's date tor a custom date

\begin{document}

\maketitle % Print the title

\abstract{
We present a new approach for the mechanically consistent modelling and simulation of fluid-structure interactions with contact. 
The fundamental idea consists of combining a relaxed contact formulation with the modelling of seepage through a porous layer 
of co-dimension 1 during contact. For the latter, a Darcy model is considered in a thin porous layer attached to a solid boundary in the limit of infinitesimal thickness. In combination with a relaxation of the contact conditions the computational model is both mechanically consistent and simple to implement. We analyse the approach in detailed numerical studies with both thick- and thin-walled solids, within a fully Eulerian and an immersed approach for the fluid-structure interaction and using fitted and unfitted finite element discretisations.
}
%\tableofcontents
\section{Introduction}
%{Numerically the treatment of the changing in time of the contact area, i.e. the real change in topology, is not trivial. In the following work we will avoid the explicit treatment of the change in topology, by considering an extension of the relaxed contact formulation presented in \cite{BurmanFernandezFrei}. The idea is to place a virtual obstacle at distance $\varepsilon \eqd \varepsilon(h)$ from the porous layer $\Sigma_p$. 
%	
%	The advantages of this method are clearly its simplicity, in fact we avoid the numerical difficulties deriving form any change in topology, and its richness from the modelling point of view. During contact an infinitesimal layer always remains between the structure and the porous layer, furthermore the fluid stresses and velocity are very close to the porous pressure and velocity respectively. Differently from the relaxed contact formulation proposed in \cite{BurmanFernandezFrei}, the porous-contact approach gives a physical meaning to the stresses generated in the infinitesimal fluid layer, allowing the problem to be physically consistent also during contact. }
%%%%%%%%%%%%%%%%%%%%%%%%%%%%%%%%%%

The design and analysis of computational methods for systems where several solids are immersed in a fluid and that can come into contact is an outstanding problem. Already fluid-structure interaction (FSI) without contact is challenging due to the moving geometries and the stiff coupling between the solid and the fluid systems. If contact between solids is to be modelled as well, the complexity increases drastically. Indeed, the addition of contact introduces several important aspects, such as: 
\begin{itemize}
\item Topological changes in the fluid domain;
\item Non-linearly changing interface conditions: The interface condition changes from a fluid-solid interaction to a solid-solid contact problem which is described by variational inequalities;
\item Important differences in the characteristic scales of the different physical phenomena: The contact represents a singular phenomenon in time. In three space dimensions, the contact zone is a two dimensional subset of the solid-solid interface and there is also a one dimensional subset of the contact zone forming the solid-solid-fluid line.
\end{itemize}
Ideally, a computational method should be consistent with the physics, be amenable to mathematical analysis and convenient to implement in a computational software. 

In the case of fluid-structure interaction with contact, an additional complication is that it is unclear what mathematical modelling will produce the best results. Indeed, it is known that a naive imposition of no-slip conditions on one of the boundaries of a solid will prevent contact of smooth bodies in the solution of the PDE system~\cite{HeslaPhD, Hillairet2d, GerardVaretetal2015}, contrary to what is observed in experiments~\cite{HagemeierTheveninRichter2020}.

Therefore, the design of computational methods for FSI-contact problems can not be completely dissociated from the  problem of modelling, but it is important to keep a certain flexibility concerning the contact modelling to be able to include a wide range of physics, depending on the characteristics of the considered system. In this work, we will build on previous work for FSI with contact \cite{BurmanFernandezFrei2020}. There, recent techniques merging the ideas of weak imposition of fluid-structure interface conditions \cite{MR2087009, BurmanFernandez2014} with a multiplier free formulation for contact \cite{MR3045657} were developed, leading to an automatic handling of both fluid-solid and 
solid-solid coupling conditions. The main idea was to merge different versions of Nitsche's method using an augmented Lagrangian formulation for variational inequalities dating back to Rockafellar \cite{MR371416}. The resulting method is consistent and shown in numerical examples to be both accurate and robust. It can also easily be combined with tools developed to facilitate the handling of the moving interfaces such as cutFEM \cite{MR3416285,BurmanFernandez2014}, XFEM~\cite{XFEM}, GFEM~\cite{GFEM} or fitted finite element approaches~\cite{FreiRichter2014}.

Already in \cite{BurmanFernandezFrei2020} some modelling aspects were developed. In order to avoid the singularity of vanishing pressure in the contact zone the fluid was extended into the solid in a porous medium model. Alternatively, the distance between the contacting bodies can be lower bounded by some small value (proportional to the mesh size) representing the idea that a fluid layer always remain between the contacting bodies. A similar approach was taken in \cite{MR4007823}, but here the physical modelling went further, introducing porous layers on the solids and modelling the full 3D poro-elastic fluid-structure interaction with contact. In the computational model for contact, a Lagrange multiplier technique was applied in contrast to the Nitsche-approach in \cite{BurmanFernandezFrei2020}.
 In \cite{ZoncaAntoniettiVergara2020} a polygonal DG discretisation is used within a penalty method in combination with a switch between Navier-slip to slip conditions to enable the transition to contact.
The recent work \cite{Ageretal2020} showcases the potential of combining the Nitsche FSI-contact conditions of \cite{BurmanFernandezFrei2020} with the FSI-cutFEM approach from \cite{BurmanFernandez2014} using realistic physical models in some impressive computational examples.

In the present work, we wish to build on the ideas of \cite{BurmanFernandezFrei2020} by considering a model for contact with seepage due to microscopic roughness
of the contacting bodies. Let $\Omega\subset \mathbb{R}^d$ for $d=2,3$ be the overall domain consisting of fluid and solid subdomains. The seepage is modelled by the introduction of a $d-1$-dimensional porous layer that adheres to a solid boundary, where contact might take place. This can be considered as a generalised boundary condition, or a bulk surface coupling in the spirit of \cite{MR3047936}. In our case, however, the free-flow Navier-Stokes' system is coupled to a surface Darcy equation. This model goes back to~\cite{MartinJaffreRobert2005} and is of interest in its own right, as discussed in the note \cite{BurmanFernandezFreiGerosa2019}. By combining this porous layer approach for seepage with the contact approach of \cite{BurmanFernandezFrei2020}, herein extended to the case of thin-walled solids, we obtain an approach that inherits the simplicity, accuracy and robustness of \cite{BurmanFernandezFrei2020}, but provides a mechanically consistent model for fluid-structure interaction with contact. 

We implement the approach using different coordinate systems, discretisations and solid models. Concerning coordinate systems, we consider both an Immersed approach going back to Peskin~\cite{Peskin1972} as well as a Fully Eulerian approach~\cite{Dunne2006, Cottetetal2008, Richter2012b, FreiPhD}. For discretisation, we use the unfitted finite element method of \cite{MR2087009,BurmanFernandez2014,alauzet-et-al-15} 
and the two-scale interface fitting approach of \cite{FreiRichter2014}. 
We illustrate the modelling capacity in a series of computational examples in two dimensions, including a beam solid model and a thick-walled solid model.

Indeed, depending on if no-slip conditions are imposed or if the porous medium approach proposed here is used the approximations will converge to different solutions for mesh size $h\to 0$. If we consider the case of a bouncing ball the use of no-slip conditions will lead to a sequence of solutions that converge to a ball that does not bounce, whereas the solutions obtained with the porous medium approach converge to a certain bouncing height that depends on the parameters of the Darcy model. Recent comparisons of computational methods with experimental studies (%Mendeley Data, doi: 10.17632/mf27c92nc3.1],
\cite{HagemeierTheveninRichter2020}, \cite{2020arXiv201108691V}) confirm that the second behavior is the physical one. Of course the parameters of the model need to be fixed through experimental studies, or otherwise.

An outline of the paper is as follows. In Section~\ref{sec.eq}, we introduce the Navier-Stokes-Darcy coupling as well as the FSI-contact model. The variational formulation and the discretisation based on Nitsche's method is described in Section~\ref{sec:disc}. In Section~\ref{sec:num} we give detailed numerical studies both in the case of a beam model and a thick solid. We conclude in Section~\ref{sec:concl}.

%{\sfrei Some further notes by Erik:
%Relaxation avoids the need of eliminating isolated pockets of liquids.
% Motivation for seepage by creation of vacuum
%}

\section{Equations}
\label{sec.eq}

In this section, we derive the Navier-Stokes-Darcy coupling, and subsequently the
equations for fluid-structure-porous-contact interaction. 
For simplicity, we will consider that contact takes place at a given fixed plane surface. This can be either an exterior rigid wall or a symmetry boundary within the fluid domain, which is relevant for example in the case of contact between two symmetric valves. The case of to two-body contact is not considered here. %, but would be possible, for example by building on the the techniques from~\cite{Ageretal2020}.
%\mf{I do not think this case is just a mater of using \cite{Ageretal2020}, I thin-walled poro-elastic model has to be derived and it deserves a separate paper.}

%We consider both thin and thick structures and their interaction with an incompressible viscous fluid. 
The fluid equations in $\Omega^{\rm f}(t) \subset \mathbb{R}^d$ will be coupled to a fixed $(d-1)$-dimensional porous layer $\Sigma_{\rm p}$ on the exterior boundary, where contact might take place. The fluid is described by the Navier-Stokes equations in Eulerian formalism and the structure by a possibly non-linear solid model. We consider both $(d-1)$-dimensional thin-walled solids and $d$-dimensional thick-walled solids. % governed by a hyperelastic material law. 
%In the first case, the solid models will be formulated using a Lagrangian formalism, while in the latter case we will also consider a fully Eulerian approach.

\subsection{Problem setting}

\begin{figure}[t]
	\centering
	\subfigure[Coupling with a thick-walled solid]{\includegraphics[width=0.45\linewidth]{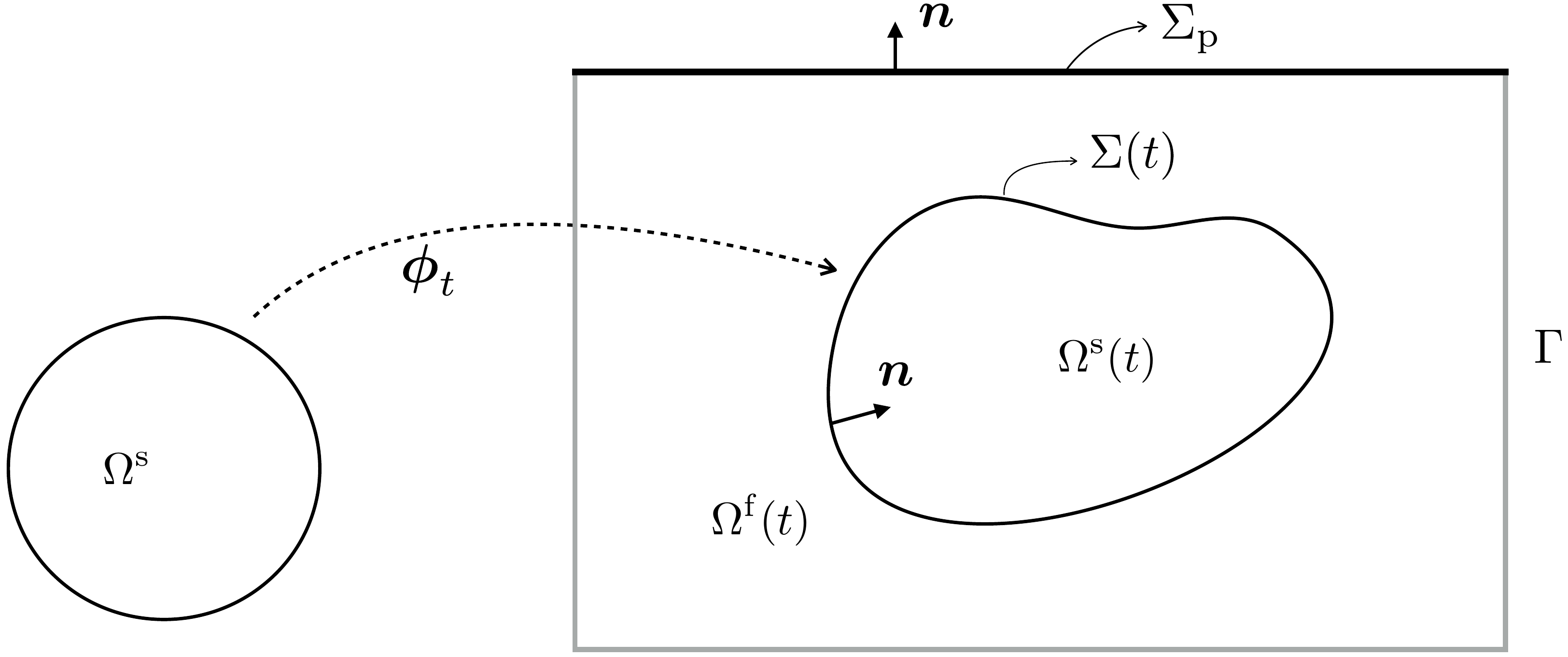}} \qquad\qquad
	\subfigure[Coupling with a closed thin-walled solid]{\includegraphics[width=0.45\linewidth]{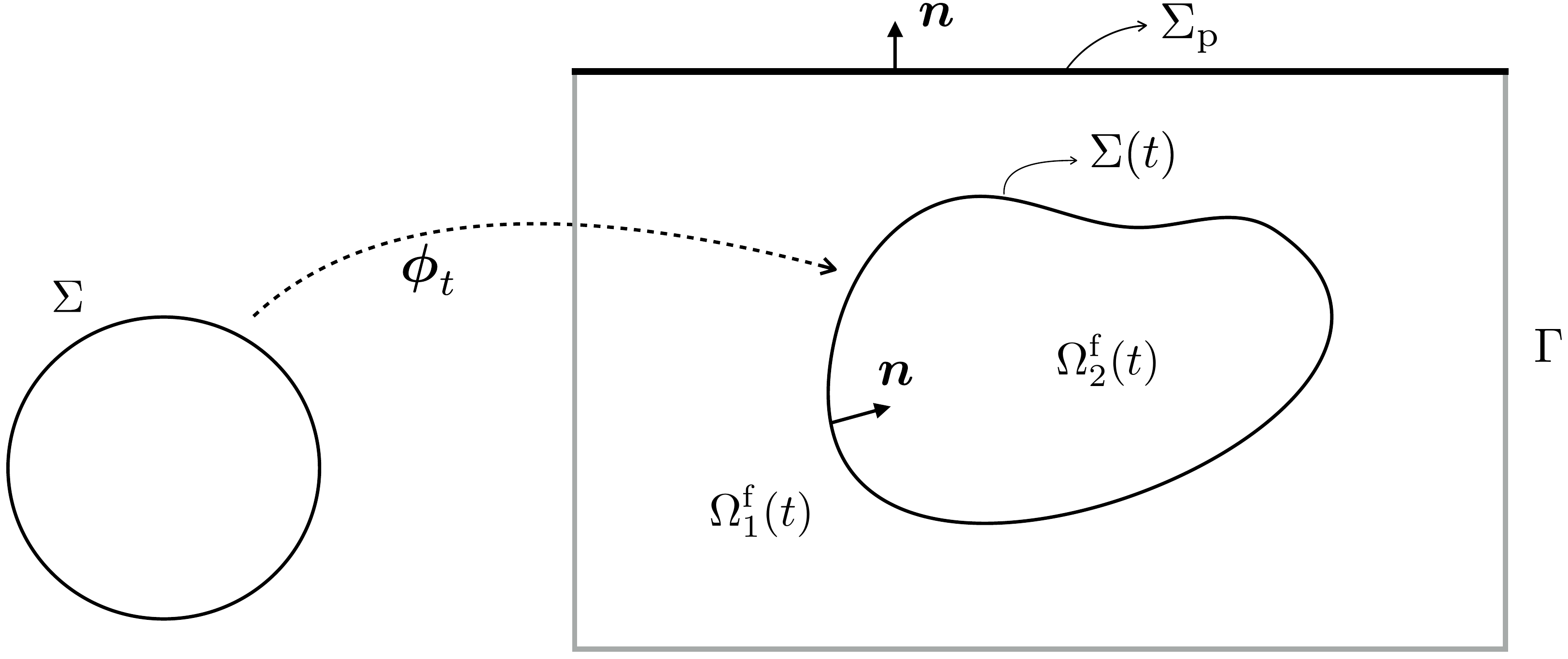}}
	\caption{Geometric configurations of the fluid and solid domains.}
	\label{fig:domains}
\end{figure}

Let $\Omega= \overline{\Omega^{\rm s}(t) }\cup \Omega^{\rm f}(t) \subset \mathbb{R}^d$  be a current configuration of the complete domain of interest, with boundary $\partial \Omega \eqd \Gamma \cup \Sigma_{\rm p} $, 
where $\Sigma_{\rm p}$ denotes the part of the boundary where contact might take place (see Figure~\ref{fig:domains}). 
There, a thin porous fluid layer is considered.  
The solid domain $\Omega^{\rm s}(t)$ can be either a surface (actually the solid mid-surface) or a domain with positive volume in $\mathbb{R}^d$ in the case of the coupling with a thick-walled solid. The current fluid-structure interface is denoted by $\Sigma(t)$ and coincides with $\Omega^{\rm s}(t)$ in the case of a thin-walled solid. The corresponding reference configurations are denoted by $\Sigma$ and $\Omega^{\rm s}$.

The structure is allowed to move freely within the domain $\Omega$. The current position of the interface $\Sigma(t)$ and the solid domain $\Omega^{\rm s}(t)$ are described in terms of a deformation map $\bphi: \Omega^{\rm s} \times \mathbb{R}^+ \longrightarrow \mathbb{R}^d$ such that $\Omega^{\rm s}(t) = \bphi(\Omega^{\rm s},t)$ and $\Sigma(t)=\bphi(\Sigma,t)$, with $\bphi \eqd \bs I_{\Omega^{\rm s}} + \bd $ and where $\dep$ denotes the solid  displacement. To simplify the notation we will refer to $\bphi_t \eqd \bphi(\cdot,t)$, so that we can also write $\Omega^{\rm s}(t) = \bphi_t(\Omega^{\rm s}), \,\Sigma(t) = \bphi_t(\Sigma)$. The fluid domain is time-dependent, namely $\Omega^{\rm f}(t)\eqd \Omega \backslash (\Omega^{\rm s}(t)\cup \Sigma(t)) \subset  \mathbb{R}^d$ with boundary $\partial \Omega^{\rm f}(t) = \Sigma(t) \cup \Gamma \cup \Sigma_p$. 
In the case of a closed thin-walled structure, the solid domain $\Omega^{\rm s}(t)$ divides 
$\Omega^{\rm f}(t)$ into two subdomains $\Omega^{\rm f}(t) = \Omega_1^{\rm f}(t)\cup \Omega_2^{\rm f}(t)$, with respective 
unit normals  $\bn_1 \eqd \bn$ and $\bn_2 \eqd -\bn $, as shown in Figure~\ref{fig:domains}(b). %The normal unit vector $\bn $ is given by the orientation of the surface $\Omega^{\rm s}(t)$. 
Similarly, in the case of a thick-walled solid, the interface $\Sigma(t)$ divides $\Omega$ into a solid part $\Omega^{\rm s}(t)$ and a fluid part $\Omega^f(t)$. We write $H^1_\Gamma(\Omega)$ for the first-order Sobolev space with vanishing trace on $\Gamma\subset \partial\Omega$.

For a given field $f$ defined in $\Omega$ (possibly discontinuous across the interface), we can define its one-sided restrictions, denoted by $f_1$ and $f_2$, as 
$$f_1(\bx) \eqd \lim_{\xi \to 0^{-}} f(\bx + \xi \bn_1), \quad f_2(\bx) \eqd \lim_{\xi \to 0^{-}} f(\bx + \xi \bn_2 ), $$
%$$f_i(\bx) \eqd \lim_{\xi \to 0^{+}} f(\bx + \xi \bn_i),\quad i=1,2,
%f_2(\bx) \eqd \lim_{\xi \to 0^{+}} f(\bx + \xi \bn_2 )
%$$
for all $ \bx \in \Sigma(t)$, and the following jump and average operators across $\Sigma(t)$:
\begin{align} \label{eq:discont}
\jump{f} \eqd f_1 - f_2 \quad 
\jump{f\bn  } \eqd f_1 \bn_1 + f_2\bn_2,\quad \average{f} \eqd \frac12\big(f_1 + f_2\big).
\end{align}
In the case of a thin structure that has a boundary inside the fluid domain (for example with a tip), these quantities can be defined similarly. For the details, we refer to~\cite{alauzet-et-al-15} and Remark~\ref{rem.open} below.

%Let first consider the coupled problem without contact, in which the fluid is coupled with the porous medium 

While the fluid and solid equations are standard and will be introduced in Section~\ref{sec.coupled}, we give some details on the porous medium model in the following section.

\subsection{Porous medium model and Navier-Stokes-Darcy coupling}

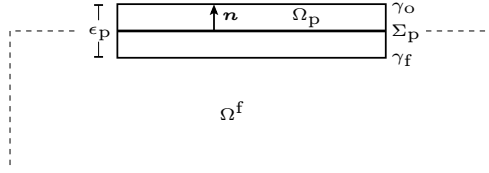
\begin{figure}[h!]
\centering
\resizebox*{0.4\textwidth}{!}{
\begin{picture}(0,0)%
\includegraphics{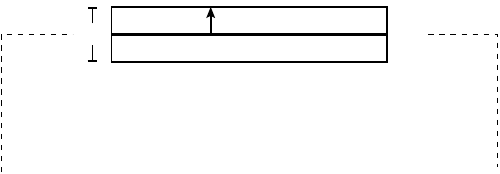}%
\end{picture}%
\setlength{\unitlength}{1160sp}%
\begingroup\makeatletter\ifx\SetFigFont\undefined%
\gdef\SetFigFont#1#2{%
  \fontsize{#1}{#2pt}%
  \selectfont}%
\fi\endgroup%
\begin{picture}(8144,2779)(1779,-3683)
\put(8191,-1906){\makebox(0,0)[lb]{\smash{{\SetFigFont{5}{6.0}{\color[rgb]{0,0,0}$\gamma_{\rm f}$}%
}}}}
\put(6526,-1231){\makebox(0,0)[lb]{\smash{{\SetFigFont{5}{6.0}{\color[rgb]{0,0,0}$\Omega_{\rm p}$}%
}}}}
\put(5356,-1231){\makebox(0,0)[lb]{\smash{{\SetFigFont{5}{6.0}{\color[rgb]{0,0,0}$\bn$}%
}}}}
\put(8191,-1051){\makebox(0,0)[lb]{\smash{{\SetFigFont{5}{6.0}{\color[rgb]{0,0,0}$\gamma_{\rm o}$}%
}}}}
\put(5311,-2896){\makebox(0,0)[lb]{\smash{{\SetFigFont{5}{6.0}{\color[rgb]{0,0,0}$\Omega^{\rm f}$}%
}}}}
\put(8191,-1488){\makebox(0,0)[lb]{\smash{{\SetFigFont{5}{6.0}{\color[rgb]{0,0,0}$\Sigma_{\rm p}$}%
}}}}
\put(3106,-1456){\makebox(0,0)[lb]{\smash{{\SetFigFont{5}{6.0}{\color[rgb]{0,0,0}$\epsilon_{\rm p}$}%
}}}}
\end{picture}%
}
\caption{\label{fig.porous} Porous medium domain $\Omega_{\rm p}$ with interface $\gamma_{\rm f}$ to $\Omega^{\rm f}$ and exterior boundary $\gamma_o$.}
\end{figure}

We consider the configuration sketched in Figure~\ref{fig.porous}, where a thin porous layer $\Omega_{\rm p} = \Sigma_{\rm p} 
\times (-\frac{\epsilon_{\rm p}}{2},\frac{\epsilon_{\rm p}}{2}) \in \mathbb{R}^d$ ($d=2,3$) with midsurface $\Sigma_{\rm p}$
is coupled to a surrounding fluid in a fixed domain $\Omega^{\rm f}\subset \mathbb{R}^d$. The surrounding fluid is governed by the Navier-Stokes equations
\begin{align}
\left\{
\begin{aligned}
\rho^{\rm{f}}\big(\partial_t\bu +  \bu \cdot \nabla  \bu \big) - \diverg\, \stress(\bu,p) = \bs 0 %\boldsymbol{f} 
&\quad \mbox{in}\quad \Omega^{\rm f} ,  \\
\diverg\, \bu = 0 &\quad \mbox{in}\quad \Omega^{\rm f}.
\end{aligned}
\right.
\end{align}
Here, $\bu$ denotes the fluid velocity and $p$ the fluid pressure in $\Omega^{\rm f}$.  The Cauchy stress tensor is given by 
$$\stress :=  2\mu_{\rm f} \boldsymbol  \varepsilon(\bu) - p \boldsymbol I, \quad \boldsymbol \varepsilon(\bu)=\frac12\left( \nabla \bu +\nabla \bu^{\rm T}\right),$$
where  $\bs I$ denotes the identity matrix. 

In the porous domain $\Omega_p$, we assume a Darcy law
%The model for the porous layer is derived from a Darcy law in
%a thin porous domain $\Omega_p = \Sigma_p \times (-\frac{\epsilon_{\rm p}}{2},\frac{\epsilon_{\rm p}}{2}) \in \mathbb{R}^d$ for $d=2,3$ with midsurface $\Sigma_p$.
% The Darcy equations in $\Omega_p$ read  
\begin{equation}
\left\{
\begin{aligned}
 \bu_{\rm l}+ \bs K \nabla p_{\rm l} = 0 \quad\mbox{in}\quad \Omega_{\rm p},\\
\nabla \cdot \bu_{\rm l} = 0 \quad\mbox{in}\quad \Omega_{\rm p},
\end{aligned}\right.
\end{equation}
where $\bu_{\rm l}$ denotes the Darcy velocity, $p_{\rm l}$ the Darcy pressure and $\bs K$ is a $d\times d$ matrix such that the following decomposition holds 
\begin{align*}
\bs K\nabla p_{\rm l} = K_\tau \nabla_{\tau} p_{\rm l} + K_n \bs \partial_n p_{\rm l},
\end{align*}
with $K_\tau,K_n\in \mathbb{R}^+$.
Here, $\bn$ is the unit normal vector of the mid-surface $\Sigma_{\rm p}$ that points towards the exterior boundary $\gamma_{\rm o}$, $\bs \partial_n = \bn \partial_n $ and $\nabla_{\tau} := P_\tau \nabla$ stands for the corresponding tangential part of the gradient 
$$P_\tau := (\boldsymbol I - \bn \otimes \bn).$$ 
Furthermore, we assume that the porous layer is very thin and consider the limit case $\epsilon_{\rm p}\to 0$. Let the outer boundary of $\Omega_{\rm p}$ be denoted by $\gamma_{\rm o}$ 
and the interior boundary connecting to the fluid domain $\Omega^{\rm f}$ by $\gamma_{\rm f}$, see Figure~\ref{fig.porous}.
We assume zero normal velocity ($\bu_{\rm l}\cdot \bn=0$) on the outer boundary $\gamma_{\rm o}$ and continuity of normal velocities and normal stresses on $\gamma_{\rm f}$. For the tangential fluid stresses, we consider the Beavers-Joseph-Saffman coupling conditions~\cite{Saffman71}.
%two possibilities: either the %Beavers-Joseph-Saffman coupling conditions~\cite{Saffman71} or a symmetry condition (vanishing tangential stresses). 
By $\bu_\tau := P_\tau \bu$, we denote the tangential part of the velocity vector and by $\sigma_{\mathrm f,n}=\bn^{\rm T}\stress \bn$ and $\sigma_{\mathrm f,\tau} = P_\tau \stress \bn$ the normal and tangential part of the Cauchy stress tensor $\stress$ introduced above. The coupling conditions between porous medium $\Omega_{\rm p}$ and fluid $\Omega^{\rm f}$ read
 \begin{equation}\label{eq:intf}
 \left\{
 \begin{aligned}
 \sigma_{\mathrm f,n} = - p_{\rm l}  & \qquad\mbox{on}\quad \gamma_{\rm f}, \\
  \bu\cdot \bn = \bu_{\rm l}\cdot \bn& \qquad\mbox{on}\quad \gamma_{\rm f}, \\
 \sigma_{\mathrm f,\tau}  = -\frac{\alpha}{\sqrt{K_\tau\epsilon_{\rm p}}} {\bu}_\tau  & \qquad\mbox{on}\quad \gamma_{\rm f}.
 \end{aligned}\right.  \end{equation}
 We note that the condition for the tangential stresses, in the last line, corresponds to a Navier-slip boundary condition for the fluid. In contrast to this boundary condition for the fluid, the normal velocity $\bu\cdot \bn$ is not zero here, as the fluid can enter the porous layer. 
 
 The appropriate choice of the parameter $\alpha$ in the last line of \eqref{eq:intf} depends on the application. In the case that $\gamma_{\rm f}$ correponds to a symmetry boundary within a larger fluid domain, where contact can take place, for example between two contacting valves, it is appropriate to set $\alpha=0$ (pure slip). If the porous layer is, however, placed at a rigid wall, the Beavers-Joseph-Saffman condition with $\alpha> 0$ is more appropriate~\cite{Saffman71, MikelicJaeger2000}. The parameter $\alpha$ depends on the structure of the porous layer. Values $0.01<\alpha<5$ have been suggested in~\cite{NieldBejan2006}. We will consider both kind of conditions in the numerical examples of Section~\ref{sec:num}.
 
 Introducing the averaged porous pressure $P_{\rm l}$ as \begin{equation}\label{eq:meanp}
P_{\rm l}= \frac{1}{2}\left(p_{\rm l}\vert{\gamma_{\rm f}} + p_{\rm l}\vert{\gamma_{\rm o}}  \right)  \quad\mbox{in}\quad \Sigma_{\rm p},
\end{equation}
 the following equations can be derived in the limit case $\epsilon_{\rm p}\to 0$ (see \cite{MartinJaffreRobert2005,BurmanFernandezFreiGerosa2019})
\begin{align}
\label{eq:darcyfluid}
\left\{
\begin{aligned}
- \nabla_\tau \cdot \left( \epsilon_{\rm p}K_\tau \nabla_\tau P_{\rm l}  \right)  =   \bu \cdot \bn & \quad\mbox{on}\quad \Sigma_{\rm p},\\
\sigma_{\mathrm f,n}  = - P_{\rm l} - \frac{\epsilon_{\rm p} K_n^{-1} }{4} \bu \cdot \bn  & \quad\mbox{on}\quad \Sigma_{\rm p},\\
\sigma_{\mathrm f,\tau}  = -\frac{\alpha}{\sqrt{K_\tau\epsilon_{\rm p}}} \bu_\tau  &\quad\mbox{on}\quad \Sigma_{\rm p}.
\end{aligned}
\right. % \qquad\qquad\mbox{in}\quad \Sigma_{\rm p}.
\end{align}
Note that the only remaining porous medium variable is the averaged pressure $P_{\rm l}$. In the limit $K_n, K_\tau\to 0$, the coupling conditions turn into a Navier-slip boundary condition for the fluid on $\Sigma_{\rm p}$. 

\subsection{Fluid-structure-porous-contact interaction model}
\label{sec.coupled}
We assume that  $\Omega^{\rm f}(t)$ is filled by an incompressible fluid governed by the Navier-Stokes equations. The domain $\Omega^{\rm s}(t)$ is occupied by a solid media described by 
a beam or shell solid model (given in terms of an abstract surface differential operator $\bs L$) on a $(d-1)$-dimensional domain $\Sigma$ or by the elastodynamics equations 
in the case of a $d$-dimensional domain $\Omega^{\rm s}$. The fluid and solid equations are coupled with no-slip interface conditions on the fluid-structure interface $\Sigma(t)$.
The solid is constrained to not penetrate into the porous medium $\Sigma_{\rm p}$ via the (relaxed) unilateral frictionless contact conditions 
\begin{equation}\label{eq:contact_ineq}
\bd\cdot \bn - g_\epsilon \leq 0, \quad \lambda \leq 0, \quad \lambda (\bd\cdot \bn - g_\epsilon)  = 0 \quad \text{ on } \Sigma,
\end{equation}
Here, $g_\epsilon :=g-\epsilon_g$, where $g$ denotes the gap function to $\Sigma_p$  and $\epsilon_g > 0 $ is a small parameter. 
The symbol $\lambda$ stands for the normal component of the contact traction, which corresponds to the Lagrange multiplier associated to
 the no-penetration condition.

The proposed fluid-structure-porous-contact interaction model is hence formulated as follows: Find the fluid velocity and pressure $\bu:\Omega^{\rm f}  \times \mathbb{R}^+ \rightarrow \mathbb{R}^d$, $p:\Omega^{\rm f} \times \mathbb{R}^+\rightarrow \mathbb{R}$, 
the solid displacement and velocity $\dep:\Omega^s \times \mathbb{R}^+ \rightarrow \mathbb{R}^d $, $\vel:\Omega^s \times \mathbb{R}^+ \rightarrow \mathbb{R}^d $, the Darcy porous pressure $P_{\rm l}:\Sigma_p \times \mathbb{R}^+\rightarrow \mathbb{R}$  and the Lagrange multiplier $\lambda : \Sigma \times \mathbb{R}^+ \rightarrow \mathbb{R}$ such that, for all $t\in \bbR^+$, the following relations are satisfied 
\begin{itemize}
\item Fluid problem:
\begin{align}
\label{eq:fluid-nl}
\left\{
\begin{aligned}
\rho^{\rm{f}}\big(\partial_t\bu +  \bu \cdot \nabla  \bu \big) - \diverg\, \stress(\bu,p) = \bs 0 %\boldsymbol{f} 
&\quad \mbox{in}\quad \Omega^{\rm f}(t) ,  \\
\diverg\, \bu = 0 &\quad \mbox{in}\quad \Omega^{\rm f}(t),   \\
\bu =  \bs 0&  \quad \mbox{on} \quad \Gamma,
\end{aligned}
\right.
\end{align}

\item Porous layer:
\begin{align}
\label{eq:darcy}
\left\{
\begin{aligned}
- \nabla_\tau \cdot \left( \epsilon_{\rm p} K_\tau \nabla_\tau P_{\rm l}  \right) =   \bu_{\rm l} \cdot \bn & \quad\mbox{on}\quad \Sigma_{\rm p},\\
\epsilon_{\rm p} K_\tau  \tau \cdot \nabla_\tau P_{\rm l} =  0  &\quad\mbox{on}\quad \partial \Sigma_{\rm p},
%\sigma_{f,n} = \underbrace{- P_l  - \frac{\epsilon_{\rm p} K_n^{-1} }{4} \bu \cdot \bn}_{=:\sigma_p} & \quad\mbox{on}\quad \Sigma_p,\\
% \sigma_{f,\tau}  = -\frac{\alpha}{\sqrt{K_\tau\epsilon_{\rm p}}} u_\tau \quad & \quad\mbox{on}\quad \Sigma_p.
%\\
%\tau^T \sigma_f n =0 & \quad\mbox{on}\quad \Sigma_p.
\end{aligned}
\right.  %\quad\mbox{in}\quad \Sigma_l.
\end{align}
%where we have considered $u_{l,n} =0$ on $\gamma^-$.

\item Solid problem: 
\begin{align}
\label{eq:solid-thin}
\left\{
\begin{aligned}
\rho^{\rm{s}}\epsilon^{\rm{s}}\partial_t  \vel 
+ \bs  L( \dep )    
=  \bs T   & \quad  \mbox{on}\quad \Omega^{\rm s}=\Sigma,\\
\vel = \partial_t \dep  & \quad \mbox{on}\quad  \Omega^{\rm s}=\Sigma,\\
\dep = \bs 0 & \quad \mbox{on}\quad \partial\Omega^{\rm s}\cap \Gamma,
\end{aligned}
\right. 
\end{align}
in case of a thin-walled solid, or  
\begin{align}
\label{eq:solid-nl}
\left\{
\begin{aligned}
\rho^{\rm{s}}\partial_t  \vel 
- \bs  \diverg\, \sstress\left( \dep \right)   
=  0   & \quad  \mbox{on}\quad \Omega^{\rm s},\\
\vel = \partial_t \dep  & \quad \mbox{on}\quad  \Omega^{\rm s},\\
\dep = \bs 0 & \quad \mbox{on}\quad \partial\Omega^{\rm s}\cap \Gamma,
\end{aligned}
\right.
\end{align}
in case of a thick-walled solid.
\item Contact conditions:
\begin{equation}\label{eq:contact_ineq2}
\bd\cdot \bn - g_\epsilon \leq 0, \quad \lambda \leq 0, \quad \lambda (\bd\cdot \bn - g_\epsilon)  = 0 \quad \text{ on } \Sigma.
\end{equation}
\item Fluid-structure coupling conditions:
\begin{align}
\label{eq:coupling-nl}
\left\{
\begin{aligned}
&  \bphi =  \bs I_{ \Omega^{\rm s}  } + \bd ,\quad   \Omega_s(t) = \bphi_t (\Omega^{\rm s}), \quad  \Omega^{\rm f}(t) = \Omega \backslash \Omega^{\rm s}(t),\\
&\bu = {  \vel \circ \bphi_t^{-1}}  \quad \mbox{on} \quad \Sigma(t), 
\end{aligned}
\right.
\end{align}
and 
\begin{align}\label{eq:coupling-thin}
& \int_\Sigma  (\bs T- \lambda \bn)  \cdot \bw = - \int_{\Sigma(t)}\jump{ \stress(\bu, p )\bn} \cdot {\bw\circ \bphi_t^{-1} },
\end{align}
or
\begin{align}\label{eq:coupling-thick}
& \int_{\Sigma}  (\bs \sigma_{\rm s} - \lambda \bs I)\bn  \cdot \bw = - \int_{\Sigma(t)}\stress(\bu, p )\bn \cdot {\bw\circ \bphi_t^{-1} }
\end{align}
for all test functions $\bw$, respectively in the case of the coupling with a thin- or thick-walled solid. 

\item Fluid-porous coupling conditions:
\begin{align}
\label{eq:darcy-coupling}
\left\{
\begin{aligned}
 \bu_{\rm l} \cdot \bn  =   \bu \cdot \bn & \quad\mbox{on}\quad \Sigma_{\rm p},\\
%\epsilon_{\rm p} K_\tau  \tau \cdot \nabla_\tau P_l =  0  &\quad\mbox{on}\quad \partial \Sigma_p,\\
\sigma_{\mathrm f,n} = \underbrace{- P_{\rm l}  - \frac{\epsilon_{\rm p} K_n^{-1} }{4} \bu \cdot \bn}_{\displaystyle =:\sigma_{\rm p}} & \quad\mbox{on}\quad \Sigma_{\rm p},\\
 \sigma_{\mathrm f,\tau}  = -\frac{\alpha}{\sqrt{K_\tau\epsilon_{\rm p}}} u_\tau \quad & \quad\mbox{on}\quad \Sigma_{\rm p}.
%\\
%\tau^T \sigma_f n =0 & \quad\mbox{on}\quad \Sigma_p.
\end{aligned}
\right.  %\quad\mbox{in}\quad \Sigma_l.
\end{align}
%where we have considered $u_{l,n} =0$ on $\gamma^-$.
\end{itemize}

The relations in \eqref{eq:coupling-nl}-\eqref{eq:coupling-thick} enforce the geometrical compatibility and the kinematic and the dynamic coupling at the interface between the fluid and the solid, respectively. 
It should be noted that  the no-penetration condition in \eqref{eq:contact_ineq2} is already imposed at an $\epsilon_g$-distance to the porous layer $\Sigma_{p\rm }$.
This modeling simplification circumvents most of the numerical difficulties  associated with the topological change in the fluid domain induced by  
the exact contact condition (i.e., with $\epsilon =0$), such as switching between the contact and fluid-solid interfaces and presence of isolated small fluid regions 
(see~\cite{Ageretal2020}). Moreover, it also facilitates the explicit treatment of the geometric condition in the fluid-structure coupling (see Section~\ref{sec:disc}).

\subsubsection{Mechanical consistency}

In the fluid-structure-porous-contact interaction model \eqref{eq:fluid-nl}-\eqref{eq:darcy-coupling}, a very thin fluid layer always remains between the solid and porous medium during contact. Owing to the relations \eqref{eq:darcy-coupling}, the behavior of the fluid confined in the contact layer is expected to be very close to the one of the porous fluid. Indeed, this is a consequence of the kinematic-dynamics relations \eqref{eq:darcy-coupling}$_{1,2}$, which are enforced both during and in absence of contact.
If a part of $\Sigma(t)$ is in contact with  $\Sigma_{\rm p}$ accoding to \eqref{eq:contact_ineq}, the value of $\sigma_{\mathrm f,n}$ (resp. $\bu\cdot \bn$) on this part $\Sigma(t)$ will be close to $\sigma_{\rm p}$  (resp. $\bu_{\rm l}\cdot \bn$) on the corresponding part of 
$\Sigma_{\rm p}$. As a result, all the kinematic and dynamic relations acting on the solid during contact have a physical meaning, which guarantees the mechanical consistency of the 
proposed fluid-structure-contact interaction model.

%Hence, the fluid stresses are very close to the "porous stresses" $\sigma_p$ for small $\epsilon_{\rm p}$.
More precisely, owing to \eqref{eq:coupling-thick}, in the case of a thick-walled solid the Lagrange multiplier for the no-penetration condition will formally 
assume the form 
\begin{align*}
\lambda = \underbrace{\sigma_{\mathrm s,n} - \sigma_{\mathrm f,n}}_{\displaystyle =:\jump{\sigma_n}} \approx \sigma_{\mathrm s,n} - \sigma_{\rm p}  \circ \pi \quad \text{on}\quad  \Sigma,
\end{align*}
where we write $\sigma_{\mathrm s, n}:=\bn^{\rm T} \bs \sigma_{\rm s} \bn$ for the solid normal traction and $\pi$ denotes a (closest-point) projection from $\Sigma(t)$ to $\Sigma_{\rm p}$. 
Both, the solid stresses $\sigma_{\mathrm s,n}$ and the "porous stresses" $\sigma_{\rm p}$ have a physical meaning during solid-porous contact. 
Hence, this porous-contact approach gives a physical meaning to the stresses generated in the infinitesimal fluid layer, 
in contrast to the relaxed contact formulation in \cite{BurmanFernandezFrei2020}, where the fluid stresses $\sigma_{\mathrm f,n}$ 
did not allow for a direct physical interpretation. A similar argumentation holds true for the thin-walled solid case, where 
\begin{align*}
\lambda = \big(\rho^{\rm{s}}\epsilon^{\rm{s}}\partial_t  \vel 
+ \bs  L( \dep )\big)\cdot \bn + \jump{\sigma_{\mathrm f,n}} \,\approx\, \big( \rho^{\rm{s}}\epsilon^{\rm{s}}\partial_t  \vel 
+ \bs  L( \dep )\big) \cdot \bn + \sigma_{\rm p} \circ \pi - \sigma_{\mathrm f,n}\vert_2 \quad \text{on} \quad  \Sigma.
\end{align*}
 
%allowing the problem to be physically consistent also during contact. 

In the spirit of \cite{AlartCurnier91}, the Lagrange multiplier can further be eliminated, which in the case of a thick-walled solid results in the non-linear contact condition 
\begin{align}\label{sigmaPgamma}
\jump{\sigma_n} = -\gamma_{\rm C} \big[  \underbrace{\bd\cdot \bn - g_\epsilon -\gamma_{\rm C}^{-1} \jump{\sigma_n}}_{\displaystyle =:P_{\gamma_{\rm C}}(\bd\cdot \bn,\jump{\sigma_n})} \big]_+ \quad \text{on}\quad  \Sigma, 
\end{align}
for $\gamma_{\rm C}>0$. 
This can be embedded in an elegant way in the variational formulation using a Nitsche-based approach, see~\cite{BurmanFernandezFrei2020} and Section~\ref{sec:disc}.
For a thin solid, a similar approach is possible, with the additional difficulty that the normal solid traction on the mid-surface is given in terms of the normal PDE residual $\big(\rho^{\rm{s}}\epsilon^{\rm{s}}\partial_t  \vel 
+ \bs  L( \dep )\big)\cdot \bn$, which is rarely available at the discrete level. On the other hand, it has been shown (see \cite{BurmanetalObstacle,Scholz84, ChoulyHild12}) for the case of a thin-walled 
solid that a pure penalty approach (i.e.$\,$neglecting the normal traction $\lambda$ in the term $P_{\gamma_c}$) leads to a first-order approximation. The detailed variational formulations for both the case of a thick- and a thin-walled solid will be given in the next section.

The main advantage of this method is its simplicity. Moreover, it is expected from a mechanical point of view that the behavior of the fluid confined in the contact layer is close to the porous fluid, as explained above. 
%\begin{equation*}
%[v]_+ = \left\{
%\begin{aligned}
%v & \quad \mbox{if} \,\,v>0,\\
%0 & \quad \mbox{if} \,\,v\leq0.
%\end{aligned}
%\right.
%\end{equation*}
%
%{The final ingredient towards the discretisation of problem~\eqref{eq:fluid-nl}-\eqref{eq:darcy} is the contact treatment. 	%In order to avoid the different difficulties related to the exact treatment of contact. We build on the relaxed contact formulation presented in \cite{BurmanFernandezFrei}. The idea is to 
%virtual place the contacting wall at distance $\varepsilon \eqd \varepsilon(h)$ from the porous layer $\Sigma_p$. 
%\begin{equation}\label{eq:contact-penalty}
%\frac{\gamma_{\rm{c}}E\epsilon^{\rm{s}}}{h^2} \big( \big[\dep_h^n \cdot \bn_{\Sigma_p} - g+\eps_h \big]_+,\bw_h \big)_{\Sigma}, 
%\end{equation}
%where $[x]_+ \eqd \max\{0,x\}$, $\gamma_{\rm{c}}>0$ is a (dimensionless) user-defined parameter, $\eps_h > 0$ is a contact tolerance,  $\bn_{\Sigma_p}$ denotes the exterior unit normal to $\Sigma_p$  and $g:\Sigma\rightarrow\mathbb{R}^+$ refer to the gap function between $\Sigma$ and $\Sigma_p$. The gap function is defined as the initial distance of a point on $\Sigma$ to the porous layer $\Sigma_p$ in the direction of $\bn_{\Sigma_p}$.
The porous medium and the structure are always coupled with the fluid only and never directly to each other. This avoids switches in the variational formulation, which would be necessary in the transition between fluid-solid and solid-porous interaction (\cite{BurmanFernandezFreiGerosa2019}). On the other hand, the solid perceives indirectly the presence of the porous layer through the fluid stresses and velocity during contact.
The resulting numerical approach is highly competitive in terms of computational costs compared to approaches using Lagrange multipliers and/or active-sets. 
%also because no switch between contact area and no contact area is needed.}

\subsubsection{Seepage}

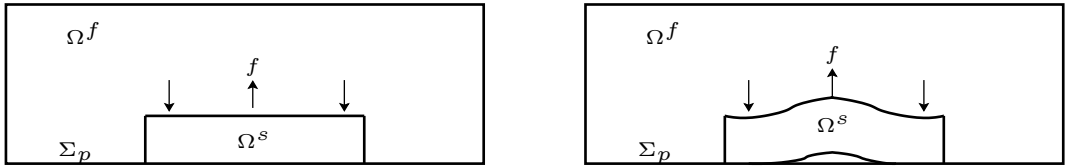
\begin{figure}[h!]
\centering
\resizebox*{0.4\textwidth}{!}{
\begin{picture}(0,0)%
\includegraphics{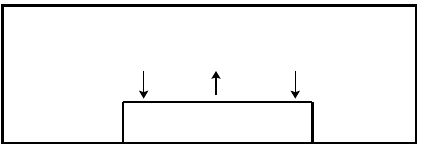}%
\end{picture}%
\setlength{\unitlength}{1450sp}%
\begingroup\makeatletter\ifx\SetFigFont\undefined%
\gdef\SetFigFont#1#2{%
  \fontsize{#1}{#2pt}%
  \selectfont}%
\fi\endgroup%
\begin{picture}(5466,1866)(1768,-2794)
\put(2476,-1411){\makebox(0,0)[lb]{\smash{{\SetFigFont{5}{6.0}{\color[rgb]{0,0,0}$\Omega^f$}%
}}}}
\put(2386,-2671){\makebox(0,0)[lb]{\smash{{\SetFigFont{5}{6.0}{\color[rgb]{0,0,0}$\Sigma_p$}%
}}}}
\put(4411,-2581){\makebox(0,0)[lb]{\smash{{\SetFigFont{5}{6.0}{\color[rgb]{0,0,0}$\Omega^s$}%
}}}}
\put(4501,-1726){\makebox(0,0)[lb]{\smash{{\SetFigFont{5}{6.0}{\color[rgb]{0,0,0}$f$}%
}}}}
\end{picture}%
}\hfil
\resizebox*{0.4\textwidth}{!}{
\begin{picture}(0,0)%
\includegraphics{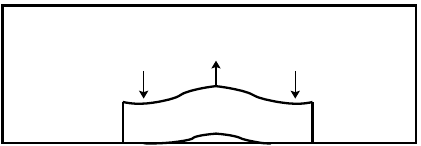}%
\end{picture}%
\setlength{\unitlength}{1450sp}%
\begingroup\makeatletter\ifx\SetFigFont\undefined%
\gdef\SetFigFont#1#2{%
  \fontsize{#1}{#2pt}%
  \selectfont}%
\fi\endgroup%
\begin{picture}(5466,1867)(1768,-2795)
\put(4411,-2401){\makebox(0,0)[lb]{\smash{{\SetFigFont{5}{6.0}{\color[rgb]{0,0,0}$\Omega^s$}%
}}}}
\put(2476,-1411){\makebox(0,0)[lb]{\smash{{\SetFigFont{5}{6.0}{\color[rgb]{0,0,0}$\Omega^f$}%
}}}}
\put(2386,-2671){\makebox(0,0)[lb]{\smash{{\SetFigFont{5}{6.0}{\color[rgb]{0,0,0}$\Sigma_p$}%
}}}}
\put(4501,-1591){\makebox(0,0)[lb]{\smash{{\SetFigFont{5}{6.0}{\color[rgb]{0,0,0}$f$}%
}}}}
\end{picture}%
}
\caption{\label{fig.vacuum} Illustrative example to motivate the role of seepage in fluid-structure interaction with contact: Contact of a solid body with the lower wall is released in the central part of the contact surface, when a specific force $f$ is applied. Without seepage on $\Sigma_p$ a vacuum would be created.
}
\end{figure}

The proposed fluid-structure-porous-contact interaction model \eqref{eq:fluid-nl}-\eqref{eq:darcy-coupling}
allows for seepage in the sense that fluid can flow through the porous layer $\Sigma_{\rm p}$, for example to connect a cavity in the central part of the contact surface with the exterior fluid. These could emerge when the impact of the structure happens in the lateral parts of the structure first or when contact of the solid is released in a central part of the contact surface only. 
This is an important aspect in the modelling of fluid-structure-contact interaction, as otherwise unphysical configurations might result. As an example consider the situation sketched in Figure~\ref{fig.vacuum}, where a solid body is in contact with the lower wall $\Sigma_{\rm p}$ at initial time (left sketch). When a (sufficiently strong) force $\bs f$ is applied in the central part of $\Omega^{\rm s}$, while the body is fixed at the lateral parts, contact will be released in the central part only. If no seepage along $\Sigma_{\rm p}$ is allowed, a vacuum would emerge between $\Omega^{\rm s}$ and $\Sigma_{\rm p}$.
While one could argue that this paradox is already circumvented by using the relaxed contact conditions \eqref{eq:contact_ineq}, we note that only the porous layer gives a physical meaning to the fluid filling 
the contact layer.

\section{Numerical methods}
\label{sec:disc}

This section is devoted to the numerical discretisation of the fluid-structure-porous-contact interaction model \eqref{eq:fluid-nl}-\eqref{eq:darcy-coupling}.
Two numerical approaches will be considered which basically depend on thin- or thick-walled nature of the solid model and 
on the formalism used for the fluid-structure coupling (mixed Lagrangian-Eulerian or fully Eulerian formalisms). 
For an accurate discretisation of the fluid-structure-porous-contact interaction model \eqref{eq:fluid-nl}-\eqref{eq:darcy-coupling}, 
two different strategies will be considered in the numerical examples reported in Section~\ref{sec:num}. The first strategy is based on an unfitted Nitsche-XFEM method, 
drawing on~\cite{BurmanFernandez2014,alauzet-et-al-15}. The second strategy is a fitted finite element method, following~\cite{FreiRichter2014}.
In both cases, we will use equal-order finite element methods.

%, that are both suitable for problems involving substantial changes in the fluid domain, including topology changes. 
%This is on one hand an immersed mixed-coordinate framework, where both fluid and solid models are formulated in their traditional Eulerian and Lagrangian coordinate frameworks, respectively, as introduced in Section~\ref{sec.eq}. As a second possibility, we consider a Fully Eulerian approach for FSI, where both sub-problems are formulated in the current, Eulerian framework~\cite{Dunne2006, FreiPhD}.

%The respective discrete spaces will be introduced in Sections~\ref{subsec.xfem}-\ref{subsec.fitted}. In the case of a thin-walled solid, the discrete fluid  pressure gradient of the velocity are allowed to be discontinuous across the interface $\Sigma(t)$.
%We will formulate the variational formulations in a mixed-coordinate framework first, i.e.\,using a Lagrangian description for the solid equations and remark afterwards in Remark~\ref{rem.Eulerian}, how the algorithm changes in a Fully Eulerian framework.

%Further, we will consider both thin and thick structures. 

To fix ideas, we will present the numerical approaches for two specific combinations of these components (mixed  Lagrangian-Eulerian vs. Fully Eulerian, fitted vs. unfitted finite elements, thick-walled vs. thin-walled structures), namely a mixed Lagrangian-Eulerian approach with a thin-walled solid using unfitted finite elements in Section~\ref{sec.immersed} and a fully Eulerian approach with a thick-walled 
solid using fitted finite element discretisation in Section~\ref{sec.Eulerian}. Different combinations are possible as well, but will not be considered in the remainder of this article.

\subsection{Lagrange-Eulerian formalism with immersed thin-walled solids}
\label{sec.immersed}

%%%%%%%%%%%%%%% discretisation nxfem 

In what follows, the parameter $\delta t>0$ stands for the time-step length and $t_n := n\delta t$ denotes the time instant at time level $n\in \mathbb{N}$. 
The symbol $x^n$ generally denotes an approximation of $x(t_n)$, for a given time valued function $x(t)$. 
We also introduced the notation  
$$\partial_{\delta t} x^{n}  := \frac{1}{\delta t }\left(x^{n}-x^{n-1}\right),$$ for the first-order backward difference.

We  consider the fluid-structure-porous-contact interaction problem \eqref{eq:fluid-nl}-\eqref{eq:darcy-coupling} in the case of the coupling with immersed thin-walled solids. 
The time discretisation is performed with a  backward-Euler scheme, including a semi-implicit treatment of the the convective term in \eqref{eq:fluid-nl} and an explicit treatment of the 
geometric coupling \eqref{eq:coupling-nl}$_1$. As regards the spatial discretisation, an unfitted finite element approximation with overlapping meshes is considered for 
 the fluid-solid coupling, by drawing on the Nitsche-XFEM method reported in \cite{alauzet-et-al-15,BurmanFernandez2014}. The fluid-porous system is discretised 
 by a standard fitted finite element approximation. 

%discretisation \eqref{eq:darcy-coupling} in a consistent unfitted finite element discretisation of the problem~\eqref{eq:fluid-nl}-\eqref{eq:darcy} which is independent of $\Sigma$, . 
%The fluid variables  $(\bu^n,p^n)$ will hence be approximated in triangulations of $\Omega^{\rm f}$. 
%To this purpose, it is important to note that the velocity gradient $\grad \bu^n$ and   the pressure $p^n$ are possibly discontinuous across $\Sigma(t)$. 
%In the following, the closed spaces $H^1_\Gamma(\omega)$, of functions in $H^1(\omega)$ with zero trace on $\Gamma$, and $L^2_0(\omega)$, of functions in $L^2(\omega)$ with zero mean in $\omega$, will   be considered. The scalar product in $L^2(\omega)$ is denoted by $(\cdot ,\cdot )_\omega$. 
%In order to simplify the presentation, we assume that both $\Omega^{\rm f}$ and $\Sigma$ are polyhedral. 

For the solid, we start by assuming that there exists a positive form
$a^{\rm s}:\bW \times \bW \longrightarrow \bbR$, linear with respect to the second argument, such that 
$$
a^{\rm s}\big(\dep,\bw) = \big(\bs L(\dep), \bw  \big)_{\Sigma} 
$$
for all $\bw \in  \bW := [H^1_{\Gamma\cap \partial\Sigma}(\Sigma)]^d$, where $\bW$ stands for the space of admissible displacements.
Let  $\{ \mathcal{T}_h^{\rm{s}} \}_{0<h < 1}$ be a family of triangulations of $\Sigma$. 
We consider the standard space of continuous 
piecewise affine functions
$$X_h^{\rm{s}}   \eqd  \left\{v_h\in C^0(\overline{\Sigma})\, \big| \, v_{h\vert K} \in \mathbb{P}_1(K), \quad \forall K\in \mathcal{T}_h^{\rm{s}} \right\}, \qquad \bW_h :=[X_h^{\rm{s}}]^d \cap \bW.$$
The contact condition \eqref{eq:contact_ineq2} is approximated via a penalty method (see, e.g., \cite{Scholz84}), by adding 
the following penalty term into the solid discrete problem:
$$ 
\frac{\gamma_{\rm{c}}E\epsilon^{\rm s} }{h^2} \big( \big[\dep_h^n \cdot \bn - g_\eps \big]_+,\bw_h \big)_{\Sigma}, 
$$
where $E$ is the solid Young modulus and $\gamma_{\rm{c}}>0$ is the (dimensionless) penalty parameter.

For a given discrete displacement approximation $\dep_h^{n} \in \bW_h$ at time $t_n$, we define its associated deformation map by $\bphi_{h}^{n} \eqd  \bs I_{\Sigma} + \dep_h^{n}$. This map characterises the current solid configuration (i.e., at time level $n$), as $\Sigma^{n} \eqd \bphi_{h}^{n}(\Sigma)$. As indicated above, we consider an explicit update for the physical fluid domain in  \eqref{eq:coupling-nl}$_1$, namely,  
	\begin{equation}\label{eq:exp-geo}
	\Omega^{\mathrm f,n} \eqd \Omega \backslash \Sigma^{n-1},
	\end{equation}
which has the effect of removing the geometrical non-linearities in the fluid problem \eqref{eq:fluid-nl}.

\begin{figure}[t]
	\centering
	\includegraphics[width=0.4\linewidth]{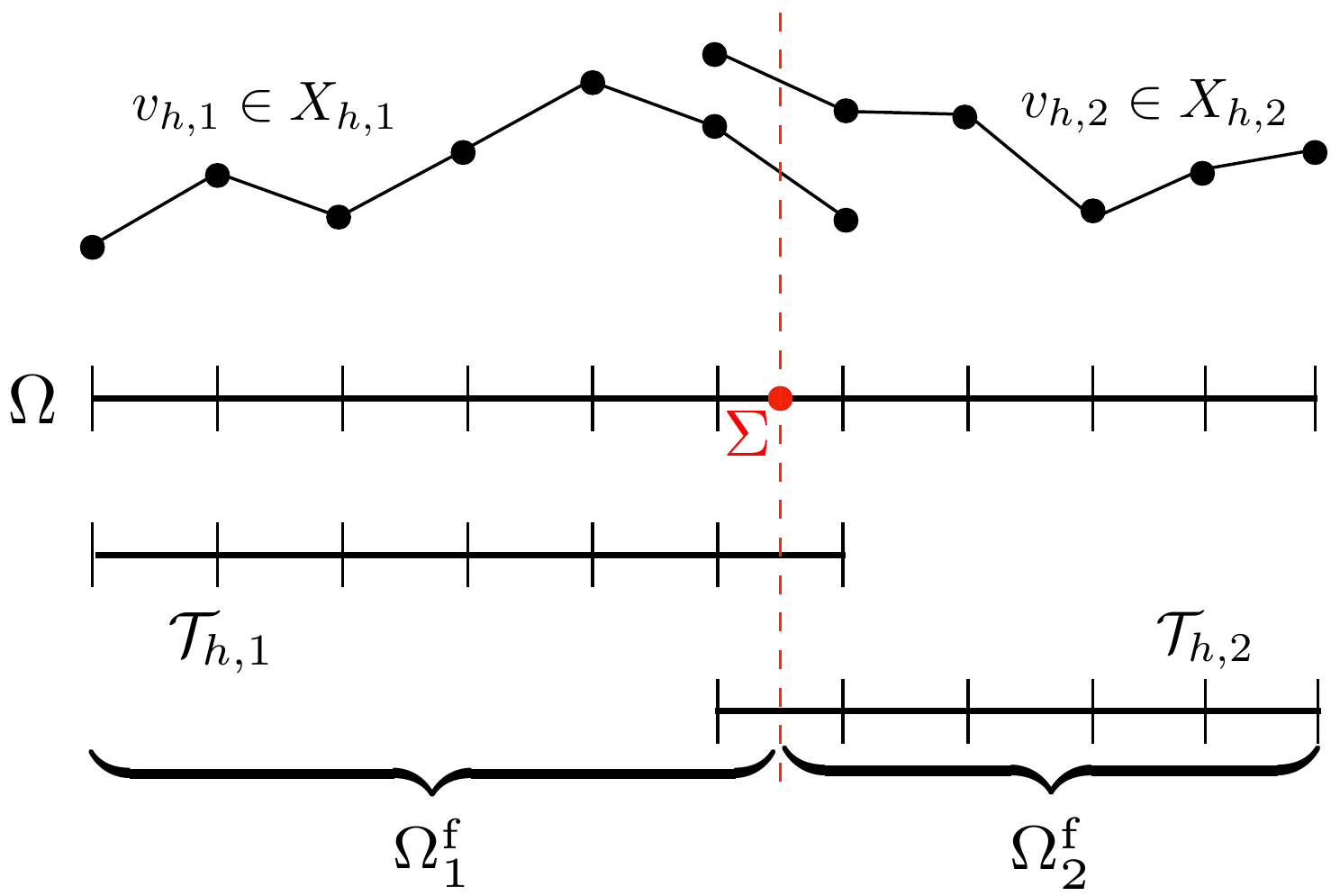}
	\caption{One dimensional illustration of the overlapping meshes $\mathcal{T}_{h,i}^n$ and of the construction of the discrete spaces $X_{h,i}^n$.}
	\label{fig:xfem}
\end{figure}

Let  $\{ \mathcal{T}_{h} \}_{0<h< 1}$ a family of triangulations of $\Omega$. 
Owing to \eqref{eq:exp-geo}, for each $\mathcal{T}_{h}$ we defined two overlapping meshes $\mathcal{T}_{h,i}^n$, $i=1,2$, such that 
$\mathcal{T}_{h,i}^n$ covers the $i$-th fluid region $\Omega_i^{\mathrm f,n}$ defined by $\Sigma^{n-1}$ through \eqref{eq:exp-geo}. Note that each triangulation 
$ \mathcal{T}_{h,i}^n$ is fitted to the exterior boundary $\Gamma \cup \Sigma_{\rm p}$, but in general not to  $\Sigma^{n-1}$ (nor $\mathcal{T}_h^{\rm{s}}$), see Figure~\ref{fig:xfem}.  
There will be  duplicated elements, i.e., such that $K \in  \mathcal{T}_{h,1}^n  \cap  \mathcal{T}_{h,2}^n $, 
but this is only allowed when $K\cap \Sigma^{n-1}\neq \emptyset$. %Note that $ \mathcal{T}_{h,1} \cup  \mathcal{T}_{h,2} $ is a triangulation of the whole fluid domain $\Omega$. 
We denote by
$\Omega_{h,i}^{\mathrm f,n} $ the domain  covered by $\mathcal{T}_{h,i}^n$, viz., 
\begin{align*}
\Omega_{h,i}^{\mathrm f,n}  \eqd \mbox{int}\left( \underset{K \in \mathcal{T}_{h,i}^n}{\bigcup}K\right).
\end{align*}
We can hence introduce the following spaces of continuous 
piecewise affine functions
$$
\begin{aligned}
X_{h,i}^n  & \eqd  \left\{v_h\in C^0(\overline{\Omega_{h,i}^{\mathrm f,n} })\,\big|  \, v_{h\vert K} \in \mathbb{P}_1(K), \quad \forall K\in \mathcal{T}_{h,i}^n \right\}.
%\bV_{h,i}^n & \eqd \left \{ \bv_{h,i} \in  [X_{h,i}^{\mathrm f,n}]^2  \big/ \, \bv_{h,i}
%\vert_{\Gamma} = 0 \right\},\quad Q_{h,i}^n\eqd X_{h,i}^{n}.
\end{aligned}
$$
%Associated with $X_{h,i}^n$, we define the spaces  
and set
$$ \bV_{h,i}^n  \eqd  [X_{h,i}^n ]^d \cap [H^1_\Gamma(\Omega^{\mathrm f,n})]^d,\quad  Q_{h,i}^n  \eqd X_{h,i}^n \cap L^2_0(\Omega^{\mathrm f,n}).
$$
For the approximation of the fluid velocity and pressure we will consider the following time-dependent discrete product spaces
\begin{equation}\label{eq:espaces}
\bV_h^n  \eqd \bV_{h,1}^n \times \bV_{h,2}^n, \quad  Q_h^n  \eqd  Q_{h,1}^n \times  Q_{h,2}^n.
\end{equation}
%where $H^1_\Gamma(\Omega)$ denotes the standard subspace of $H^1(\Omega)$ with functions  vanishing on $\Gamma$ and $L^2_0(\Omega)$ the subspace 
% of $L^2(\Omega)$ made of zero mean functions.  
%which guarantee that interfacial (strong and weak) discontinuities are included in the discrete approximation of both the fluid velocity and pressure. 
%while we use in the thick case
%\begin{equation}\label{eq:espaces_thick}
%\bV_h  \eqd [X_{h,1} ]^d \cap [H^1_\Gamma(\Omega^{\rm f})]^d,\quad  Q_h  \eqd X_{h,1}\cap L^2_0(\Omega^{\rm f}).
%\end{equation}
The functions in \eqref{eq:espaces}
are continuous in the physical fluid domain $\Omega_i^{\mathrm f,n}$,
but discontinuous across the interface location $\Sigma^{n-1}$ (see Figure~\ref{fig:xfem}).

\begin{algorithm}[h!]
	For $n\geq 1$:
	\begin{enumerate}
		\item Interface update:  
		\begin{equation*}
		\bphi_h^{n-1} =  \bs I_{\Sigma} + \dep_h^{n-1} , \quad
		\Sigma^{n-1}  = \bphi_h^{n-1}(\Sigma),  \quad  \Omega^{\mathrm f,n} = \Omega \backslash \Sigma^{n-1}. 
		\end{equation*}
		
		\item Find $\big( \bu_h^{n}, p_h^{n}, \vel_h^{n}, \dep_h^{n}, P_{\mathrm l,h}^n\big) \in \bV_h^n  \times Q_h^n \times \bW_h \times \bW_h \times \bs S_h$ with $\vel_h^{n} = \partial_{\delta t} \dep_h^{n}$ such that 
		\begin{multline}\label{eq:algofsiDarcy}
		\rho^{\rm f} \big( \partial_{\delta t} \bu_h^{n} ,  \bv_h \big)_{\Omega^n } + a_h^{\mathrm f,n}\big(\bu_h^{n-1}; ( \bu_h^{n} , p_h^{n}), (\bv_h, q_h)\big)  
		+ \rho^{\rm s}\epsilon^{\rm{s}}  \big(\partial_{\delta t} \vel_h^{n} , \bw_h  \big)_{\Sigma}  
		+ a^{\rm s}(\dep_h^{n}, \bw_h ) \\
		%%%% contact term 
		+\frac{\gamma_{\rm{c}} E \epsilon^{\rm s}}{h^2} \big( \big[\dep_h^n \cdot \n - g_\epsilon \big]_+,\bw_h  \cdot \n \big)_{\Sigma} 
		%%%% coupling FS
		-  \sum_{i=1}^2 \big(\stress(\bu_{h,i}^{n},  p_{h,i}^n  )\n_i , \bv_{h,i} - \bw_h \big)_{\Sigma^{n-1}} \\
		-\sum_{i=1}^2\big( \bu_{h,i}^{n}  - \vel_h^{n} , \stress(\bv_{h,i},-q_{h,i}))\n_i   \big)_{\Sigma^{n-1}} 
		+ \frac{\gamma \mu}{h} \sum_{i=1}^2  \big(\bu_{h,i}^{n } - \vel_h^{n} ,  \bv_{h,i} - \bw_h\big)_{\Sigma^{n-1}}  
		\\
		%%%%%%%%%% fsi side for contact 
		%		-  \big(\stress(\bu_{h,s_2}^{n},  p_{h,s_2}^n  )\n_{s_2} , \bv_{h,s_2} - \bw_h \big)_{\Sigma^n_c} 
		%		-\big( \bu_{h,s_2}^{n}  - \vel_h^{n} , \stress(\bv_{h,s_2},-q_{h,s_2}))\n_{s_2} \big)_{\Sigma^n_c} 
		%		\\
		%		+ \frac{\gamma \mu}{h} \big(\bu_{h,s_2}^{n } - \vel_h^{n} ,  \bv_{h,s_2} - \bw_h\big)_{\Sigma^n_c}  \\
		%%%% pourous stresses
		 +\frac{\alpha}{\sqrt{K_\tau\epsilon_{\rm p}}} (\bu_{h,\tau}^n, \bv_{h,\tau})_{\Sigma_{\rm p}}
		- (\sigma_{\rm p},   \bv_h\cdot \bn)_{\Sigma_{\rm p} } 
		%%%% darcy
		+ (\epsilon_{\rm p} K_\tau \nabla_\tau P_{\mathrm l,h}^n, \nabla_\tau q_{\mathrm l,h})_{\Sigma_{\rm p}}
		- \big( \bu^{n}_{h}\cdot \n, q_{\mathrm l,h}\big)_{\Sigma_{\rm p} } 
		= 0
		\end{multline}
		for all $(\bv_h , q_h, \bw_h,q_{\mathrm l,h} ) \in \bV_h^n  \times Q_h^n \times \bW_h \times \bs S_h$, where the porous stress $\sigma_{\rm p}$ is given by
		\begin{align*}
		\sigma_{\rm p}%=- P^n_{\rm l}  {\sfrei -} \frac{\epsilon_{\rm p} K_n^{-1} }{4} u_{l,n}|_{\gamma^+} 
		:= - P^n_{\mathrm l,h} { -} \frac{\epsilon_{\rm p} K_n^{-1} }{4}  \bu^n_{h} \cdot \bn  \quad \text{on}\quad \Sigma_{\rm p} .
		\end{align*}
		
	\end{enumerate}
	\caption{Strongly coupled scheme for fluid-structure-porous-contact  interaction (thin-walled solid).}
	\label{alg:FSIDarcy}
\end{algorithm}

We can now introduce the corresponding fluid discrete tri-linear form (see \cite{alauzet-et-al-15}):
$$
\begin{aligned}\label{eq:afn}
a_{h}^{{\mathrm f,n}} \big( \bz_h; (\bu_h, p_h),(\bv_h, q_h)\big) \eqd & 
%c_{h }^n (\bz_h, \bu_h, \bv_h) +
2 \mu \big(\strain(\bu_h),\strain(\bv_h)  \big)_{\Omega^{\mathrm f,n}} 
+
\rho^{\rm f} \big( \bz_h \cdot  \grad \bu_h, \bv_h \big)_{\Omega^{\mathrm f,n}}+ \frac{\rho^{\rm f}}{2}\big((\diverg \bz_h) \bu_h,\bv_h   \big)_{\Omega^{\mathrm f,n}} \\
&- \rho^{\rm f}\big(\average{\bz_h}\cdot \bn \jump{\bu_h}, \average{\bv_h}\big)_{\Sigma^{n-1}   } - \frac{\rho^{\rm f}}{2} 
\big(\jump{\bz_h\cdot\bn} ,\average{\bu_h \cdot \bv_h}\big)_{\Sigma^{n-1}  }\\
& -(p_h, \diverg\, \bv_h)_{\Omega^{\mathrm f,n}}+(\diverg\, \bu_h, q_h)_{\Omega^{\mathrm f,n}} \\
& + s_{\mathrm v,h}^{n}(\bz_h;\bu_h, \bv_h) + s_{{\rm p},h}^n(\bz_h;p_h, q_h) + g_{h}^{n}(\bu_h,\bv_h).
\end{aligned}
$$
For consistency the bulk terms  are integrated in the physical domain $\Omega^{\mathrm f,n}$, which requires a specific track of the interface intersections within the fluid domain (see e.g. \cite{MLL13,fernandez-alauzet-landajuela-16,zonca-et-al-18}). 
The terms $s_{\mathrm v,h}^{n}$ and $ s_{{\rm p},h}^n$  respectively correspond to the 
continuous interior penalty velocity and pressure stabilisation operators,
  given by (see, e.g., \cite{burman-fernandez-hansbo-06}):
$$
\begin{aligned}
s_{\mathrm v,h}^{n}(\bz_h;\bu_h, \bv_h)  \eqd  & \gamma_{\rm v} h^2 \sum_{i=1}^2   \sum_{F \in  \mathcal{F}_{h,i}^n }  \xi\big({\rm Re}_F(\bz_h)\big)  \| \bz_h \cdot \bn \|_{L^\infty(F)} \big(\jump{  \grad \bu_h }_F,\jump{  \grad \bv_h }_F\big)_F  
%\\&+  \gamma_{\rm v, 2} h^{2} \sum_{i=1}^{2} \sum_{F \in \mathcal{F}_{h, i}^{n}}\left\|z_{h}\right\|_{L^{\infty}(F)}\left( \| \operatorname{div} \boldsymbol{u}_{h}\right]_{F},\left\|\operatorname{div} \boldsymbol{v}_{h}\right\|_{F} )_{F}
, \\
s_{{\rm p},h}^n(\bz_h;p_h, q_h)   \eqd & \gamma_{\rm p}h^2 \sum_{i=1}^2   \sum_{F \in  \mathcal{F}_{h,i}^n }  \frac{ \xi\big(  {\rm Re}_F(\bz_h)\big)}{\|\bz_h\|_{L^\infty(F)}}  \big( \jump{\grad p_h}_F , \jump{\grad q_h}_F\big)_F ,
\end{aligned}
$$
where  $\mathcal{F}_{h,i}^n$  denotes the set of interior edges or faces of $\mathcal{T}_{h,i}^{n}$, 
${\rm Re}_F(\bz_h) \eqd   \rho^{\rm f} \| \bz_h \|_{L^\infty(F)} h \mu^{-1} $ denotes the local Reynolds number, $\xi(x) ~\eqd \min\{1, x\}$ is 
a cut-off function and $ \gamma_{\rm p}, \gamma_{\rm v} >0$ are user-defined parameters.
Finally, $g_{h}^{n}$ is the so-called ghost-penalty operator, given by 
$$
g_{h}^{n}(\bu_h,\bv_h) \eqd  \gamma_{\rm g} \mu h \sum_{i=1}^2   \sum_{F \in  \mathcal{F}^{\Sigma^{n-1}}_{i,h}  }   \big(\jump{ \grad \bu_{i,h} }_F , \jump{  \grad \bv_{i,h} }_F\big)_F,
$$
where $\mathcal{F}^{ \Sigma^{n-1}}_{i,h}$  denotes the set of interior edges or faces of the elements intersected by $\Sigma^{n-1}$.  
This term is added 
to guarantee robustness independent of the way the interface cuts the fluid mesh. The underlying idea is to extend the coercivity of the bi-linear form to the whole computational domain, see \cite{Bu10} or \cite{LehrenfeldOlshanskii2019} for different possibilities.

For the approximation of the porous system, we consider a family of triangulation $\{ \mathcal{T}_h^{\rm p}\}_{h>0}$  of $\Sigma_{\rm p}$, so that each $\mathcal{T}_h^{\rm p}$
is fitted $\mathcal{T}_h$. We then consider the standard space of continuous 
piecewise affine functions for the approximation of the porous pressure $P_{\rm l}$
$$\bs S_h   \eqd  \left\{v_h\in C^0(\overline{\Sigma_p})\, \big| \, v_{h\vert K} \in \mathbb{P}_1(K), \quad \forall K\in \mathcal{T}_h^{\rm p} \right\}.$$
%Finally, the following discrete space is considered for the approximation of the porous pressure  
%$$\bs S_h \eqd X_h^{p} \cap L^2(\Sigma_p).$$
%%%%%%%%%%%%%%% algorithm 

% For consistency the bulk terms in \eqref{eq:algofsiDarcy} are integrated in the physical domain, this implies integration over cut elements (see e.g. \cite{MLL13,fernandez-alauzet-landajuela-16,zonca-et-al-18} ) and therefore a specific track of the interface intersection within the fluid domain. At each time level a new interface intersection and a new sub-division of the cut-elements have to be computed. Consequently, Algorithm~\ref{alg:FSIDarcy} involves integrals of functions associated with different time levels, such as $$ \frac{\rho^{\rm f}}{\tau} \big( \bu_{h}^{n-1},\bv_h\big)$$. 
% This requires the integration of products of functions that might be discontinuous at different locations in the same element.
% In order avoid the simultaneous intersection of different interface locations with the same fluid element, we consider the approach introduced in  
% \cite{fernandez-alauzet-landajuela-16} (see also \cite{FriesZilian09}), which basically consists in locally shifting the discontinuity at time $t^{\star}$ to the structure location at time $t^n$, where $t^{\star}$ refers to $t^{n-1}$ and $t^{n-2}$  respectively.
% In the case where the discontinuities are located in different elements, the quadrature is performed in a standard fashion since we keep track of the (previous) intersections at different times and we can treat each discontinuity separately. 

 In summary, the resulting fully discrete method is reported in Algorithm~\ref{alg:FSIDarcy}. 
 We use the notation $\bu_{h,i}^n, \bv_{h,i}^n, p_{h,i}^n, q_{h,i}^n\, (i=1,2)$ introduced in \eqref{eq:discont} for the two parts of the discontinuous functions across $\Sigma^{n-1}$.
Note that the kinematic-dynamic interface 
coupling \eqref{eq:coupling-nl}$_2$-\eqref{eq:coupling-thin} is enforced in a consistent and strongly coupled fashion through Nitsche's method (see \cite{BurmanFernandez2014,alauzet-et-al-15}).
 
%\begin{remark}
%	It is worth noting that in the thin-walled solid case, when contact occurs between $\Sigma$ and $\Sigma_p$, the fluid element is duplicated. Hence, for consistency reasons, the porous layer is connected to the physically relevant component of the fluid domain.
%\end{remark}

{
\begin{remark}\label{rem.open}
	If the interface has a boundary inside the fluid domain (the so-called), we consider the construction of the fluid and solid discrete spaces proposed in \cite{alauzet-et-al-15} (see \cite[Chapter 6]{gerosa-21} for an extension to the 3D case). 
	A virtual interface $\widetilde \Sigma^{n-1} $ is introduced by connecting the interface tip with the fluid vertex opposite to the edge intersected by the interface and therefore the fluid domain is closed. 
	Afterwards, we enforce the kinematic/dynamic continuity of the fluid on  $\widetilde \Sigma^{n-1}$ in a discontinuous Galerkin fashion (see, e.g., \cite{di-pietro-ern-12}). 
	More precisely, the following terms are added to~\eqref{eq:algofsiDarcy}
	\begin{equation*}%\label{}
	-  \big(\average{\stress(\bu_{h}^{n}, p_h^n )}\n  , \jump{\bv_{h}}  \big)_{{\widetilde \Sigma}^{n-1} } 
	-   \big(\average{\stress(\bv_{h}, - q_h)}\n, \jump{\bu_{h}^{n}}  \big)_{{\widetilde \Sigma}^{n-1 }}
	+ \frac{\gamma \mu}{h} \big( \jump{\bu_{h}^{n}}  ,  \jump{\bv_{h}} \big)_{{\widetilde \Sigma}^{n-1}  }.
	\end{equation*}
\end{remark}
}

%Notice that the schemes are strongly coupled, no weakly coupled strategies are considered.
%In the case of a thin solid, the contact term needs to be scaled with $h^{-2}$~\cite{ChoulyHild12}. To solve the coupled problem we will use a semi-smooth Newton method.

%%%%%%%%%%%%%%% contact treatment 

\subsection{Fully Eulerian formalism with immersed thick-walled solids}
\label{sec.Eulerian}

In a fully Eulerian approach the current displacement $\bd(\bx,t)$ is defined by the relation
\begin{align}\label{EulerianDispl}
\bx - \bd(\bx,t) = \hat{\bx},
\end{align}
where $\hat{\bx}$ is the corresponding point in the reference configuration $\Omega^{\rm s}$. This means that the displacement can be used to trace back points $\bx\in \Omega^{\rm s}(t)$ to their reference position $\hat{\bx}$ in $\Omega^{\rm s}$ and hence to determine the domain affiliation of a point $\bx\in\Omega$ at time $t$.
 As in the previous section, we use again an explicit approach to avoid the issues related to geometrical non-linearities 
\begin{align*}
\Omega^{\mathrm s,n} :=\left\{ \bx\in \Omega, \, \big| \, \bx-\bd_h^{n-1} \in \Omega^{\rm s}\right\}, \quad \Sigma^{n} :=\left\{ \bx\in \Omega, \, \big| \, \bx-\bd_h^{n-1} \in \Sigma\right\}, \quad \Omega^{\mathrm f,n} := \Omega \setminus \left(\Omega^{\mathrm s,n}\cup \Sigma^{n}\right).  
\end{align*}
Numerically, the domain affiliations can be determined by the Initial Point Set (backward characteristics) method (see, e.g.,~\cite{Dunne2006, Cottetetal2008, FreiPhD}). To evaluate the displacement $\bd_h^{n-1}$ in points $\bx\in\Omega^{\mathrm f,n-1}$ near the interface $\Sigma^{n-1}$, that could belong to $\Omega^{\mathrm s,n}$ in the next step, an extension of the solid displacement to a small layer around the interface is required (see, e.g.,~\cite{Richter2012b}). The domain affiliation can be computed in a separate step before setting up the variational formulation as shown in Algorithm~\ref{alg:thick}, or "on-the-fly" while setting up the finite element formulation.

\begin{figure}[h!]
\centering
\begin{picture}(0,0)%
\includegraphics{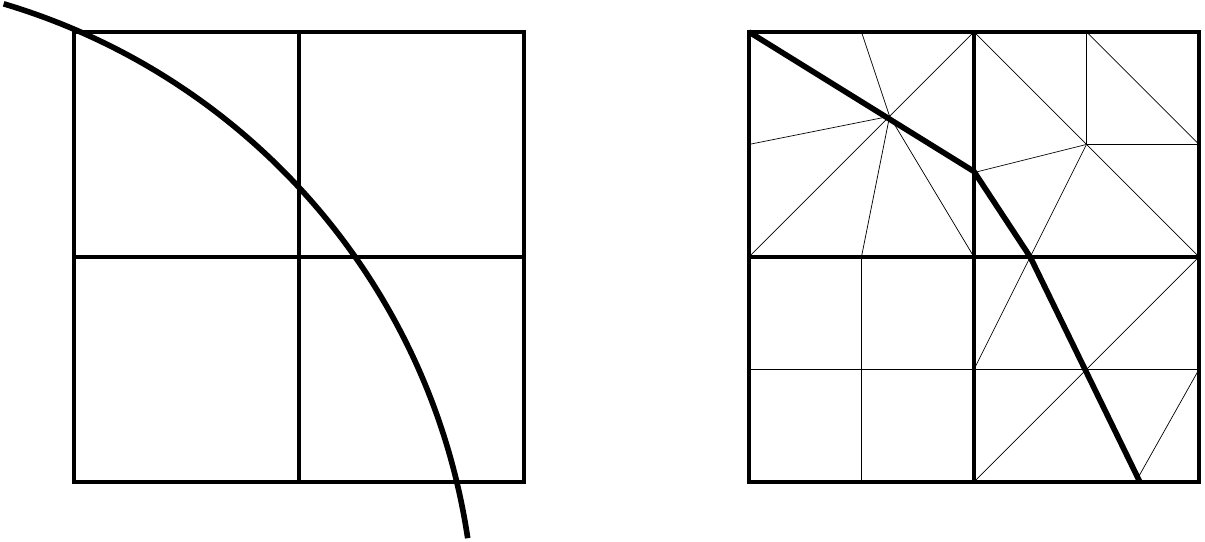}%
\end{picture}%
\setlength{\unitlength}{1776sp}%
\begingroup\makeatletter\ifx\SetFigFont\undefined%
\gdef\SetFigFont#1#2{%
  \fontsize{#1}{#2pt}%
  \selectfont}%
\fi\endgroup%
\begin{picture}(12836,5776)(413,-5799)
\put(12710,-5613){\makebox(0,0)[lb]{\smash{{\SetFigFont{8}{9.6}{\color[rgb]{0,0,0}$\Sigma_h$}%
}}}}
\put(5506,-5611){\makebox(0,0)[lb]{\smash{{\SetFigFont{8}{9.6}{\color[rgb]{0,0,0}$\Sigma$}%
}}}}
\put(2401,-5611){\makebox(0,0)[lb]{\smash{{\SetFigFont{8}{9.6}{\color[rgb]{0,0,0}$\Omega^{\rm f}$}%
}}}}
\put(6301,-2461){\makebox(0,0)[lb]{\smash{{\SetFigFont{8}{9.6}{\color[rgb]{0,0,0}$\Omega^{\rm s}$}%
}}}}
\end{picture}%
\caption{Example of a locally fitted finite element discretisation with 4 coarse cells.\label{fig.lmfem}}
\end{figure}
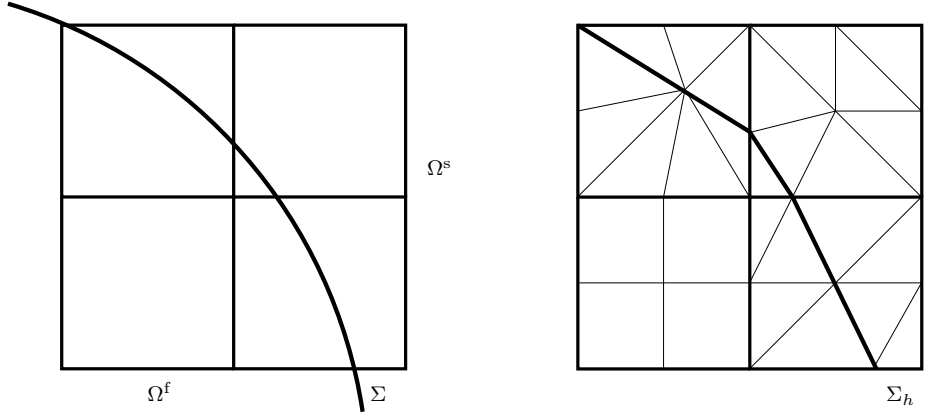

As an alternative to the unfitted finite element method presented in the previous section, we consider here a fitted finite element method for spatial discretisation. We briefly describe the locally modified finite element method as an example in two space dimensions  (see ~\cite{FreiRichter2014}). The method is based on a fixed coarse triangulation ${\cal T}_{2h}$, which is independent of the interface position, and a further subtriangulation of the coarse elements $P\in {\cal T}_{2h}$, which resolves the interface, see Figure~\ref{fig.lmfem}. We restrict ourselves to linear finite elements and a linear interface approximation. A second-order approximation has been presented recently in~\cite{FreiJudakovaRichter20}.

In order to resolve the interface locally, we split each coarse cell cut by the interface into 8 subtriangles $K_1,...,K_8$ and move some of the interior vertices to the interface, such that a linear interface approximation is obtained. The position of the 9 degrees of freedom $x_1,...,x_9$ in each coarse cell are described by a piecewise linear reference map
from the reference patch $\hat{P}=[0,1]^2$ 
\begin{align*}
\xi_P: \hat{P} \mapsto P, \qquad 
\xi_P \in Q_P :=\big\{ \phi_h \in C^0(\hat{P})\, \big|\,\, \xi|_{K_i} \in P_1(K_i),\, i =1,...,8\big\} 
\end{align*}
that fulfills the 9 conditions $\xi_P(\hat{x}_i) = x_i, \, i=1,...,9$, where $\hat{x}_i$ denotes the (fixed) Lagrangian points on the reference patch. In elements that are not affected by the interface piecewise bilinear shape functions can be used on four quadrilaterals alternatively, see the lower left coarse cell in Figure~\ref{fig.lmfem}.
The locally modified finite element space is then given by 
\begin{align*}
X_h^{{\rm lmfem},n}(\Omega) := \big\{ \phi_h\in C^0(\overline{\Omega}), \,\big|\, (\phi_h\circ \xi^{-1})|_P \in Q_P \, \forall P\in {\cal T}_{2h}\big\}. 
\end{align*}
By $X_h^{\rm lmfem,n}(\Omega^i)$ we denote the space that results by eliminating all degrees of freedom that do not lie in a subdomain $\Omega^i\subset\Omega$ or on its boundary.
Then, we define the spaces 
\begin{align*}
\bV_h^n := \big[ X_h^{{\rm lmfem},n}(\Omega^{\mathrm f,n})\cap H^1_\Gamma(\Omega^{\mathrm f,n})\big]^2, \quad Q_h^n := X_h^{{\rm lmfem},n}(\Omega^{\mathrm f,n}), \quad \bW_h^n := \big[ X_h^{{\rm lmfem},n}(\Omega^{{\rm s},n})\cap H^1_\Gamma(\Omega^{\rm s,n})\big]^2.
\end{align*}
The space $S_h$ is defined as in Section~\ref{sec.immersed}. We use an implicit form of the fluid semi-linear form
$$
\begin{aligned}
a_{h}^{{\mathrm f,n}} \big(\bu_h, p_h;\bv_h, q_h\big) \eqd  
%c_{h }^n (\bz_h, \bu_h, \bv_h) +
2 \mu \big(&\strain(\bu_h),\strain(\bv_h)  \big)_{\Omega^{\mathrm f,n}} 
+\rho^{\rm f} \big( \bu_h \cdot  \grad \bu_h, \bv_h \big)_{\Omega^{\mathrm f,n}} -(p_h, \diverg\, \bv_h)_{\Omega^{\mathrm f,n}}\\
&+(\diverg\, \bu_h, q_h)_{\Omega^{\mathrm f,n}}
+ s_{{\rm p},h}^n(\bz_h;p_h, q_h).
\end{aligned}
$$
To cope with the lack of inf-sup-stability of the discrete spaces, we use the (anisotropic) CIP pressure stabilisation developed in~\cite{Frei2019} for $s_{{\rm p},h}^n$. In contrast to \eqref{eq:afn}, we have omitted the convection stabilisation, which will not be needed in the examples with a thick solid below. Moreover, a ghost-penalty term is not required, as we use a fitted finite element discretisation.

For the solid, we assume a hyperelastic material law with a corresponding variational formulation of the form
\begin{align*}
a^{\mathrm s,n}(\bd_h, \bw_h) := (\sstress(\bd_h), \nabla \bw_h)_{\Omega^{\mathrm s,n}},
\end{align*}
where $\sstress$ denotes the Cauchy stress tensor.

Concerning time discretisation let us first note that the variables $\bu_h^{n-1}$ and $\bd_h^{n-1}$ are undefined on parts of $\Omega^{\rm s,n}$ and $\Omega^{\mathrm f,n}$, respectively, as both $\Omega^{\rm s}$ and $\Omega^{\rm f}$ are time-dependent. Thus, the method of lines can not be applied in a straight-forward way. To deal with this issue, we use the dG(0) variant of a family of Galerkin schemes that incorporates the characteristics of the domain movement in the Galerkin spaces~\cite{FreiRichter2017}. For the dG(0) variant the difference to a standard backward Euler scheme lies solely in the discretisation of the time derivative. We introduce the notation 
\begin{align*}
\tilde{\partial}_{\delta t} \bu_h^n = \frac{\bu_h^n - \big(\bu_h^{n-1}\circ \Psi\big)}{\delta t} - \partial_t \Psi\cdot \nabla \bu_h^n,
\end{align*}
where $\Psi$ is an (arbitrary) map defined in $\Omega$ that maps $\Omega^{{\rm f}, n}$ to $\Omega^{{\rm f}, n-1}$ and $\Omega^{{\rm s}, n}$ to $\Omega^{{\rm s}, n-1}$, respectively. Alternatively, one could use Eulerian time-stepping schemes with suitable extension operators. An implicit extension my means of ghost-penalties has been studied in~\cite{LehrenfeldOlshanskii2019, BurmanFreiMassing2020}. 

\begin{algorithm}[t]
	For $n\geq 1$:
	\begin{enumerate}
		\item Update the domain affiliations  
		\begin{align*}
\Omega^{{\rm s},n} :=\left\{ \bx\in \Omega, \, \big| \, \bx-\bd_h^{n-1} \in \Omega^{\rm s}\right\}, \quad \Sigma^{ n} :=\left\{ x\in \Omega, \, \big| \, \bx-\bd_h^{n-1} \in \Sigma\right\}, \quad \Omega^{\mathrm f,n} := \Omega \setminus \left(\Omega^{{\rm s},n}\cup \Sigma^{n}\right).  
\end{align*}
		\item Find $\big( \bu_h^{n}, p_h^{n}, \vel_h^{n} , \dep_h^{n}, P_{\mathrm l,h}^n\big)  \in \bV_h^n  \times Q_h^n \times \bW_h^n \times \bW_h^n \times \bs S_h$ with $\vel_h^{n} = \tilde{\partial}_{\delta t} \dep_h^{n}$ and
		\begin{multline}\label{varf:thick}
		\rho^{\rm f} \big( {\tilde{\partial}_{\delta t}} \bu_h^{n} ,  \bv_h \big)_{\Omega^{f,n}} + a_h^{\mathrm f,n}\big( \bu_h^{n} , p_h^{n}; \bv_h, q_h\big)  
		+ \rho^{\rm s} \big(\tilde{\partial}_{\delta t} \vel_h^{n} , \bw_h  \big)_{\Omega^{\rm s,n}}  
		+ a^{\rm s,n}(\dep_h^{n}, \bw_h ) \\
		%%%% contact term 
		+\frac{\gamma_{\rm{c}} E}{h} \big( \big[\tilde{P}_{\gamma_c}(\dep_h^n)]_+,\bw_h  \cdot \n \big)_{\Sigma^{\rm n}}
		%%%% coupling FS
		- \big(\stress(\bu_{h}^{n},  p_{h}^n  )\n , \bv_{h} - \bw_h \big)_{\Sigma^n} \\
		-\big( \bu_{h}^{n}  - \vel_h^{n} , \stress(\bv_{h},-q_{h}))\n   \big)_{\Sigma^n} 
		+ \frac{\gamma \mu}{h} \big(\bu_{h}^{n } - \vel_h^{n} ,  \bv_{h} - \bw_h\big)_{\Sigma^n}  
		\\
		%%%%%%%%%% fsi side for contact 
		%		-  \big(\stress(\bu_{h,s_2}^{n},  p_{h,s_2}^n  )\n_{s_2} , \bv_{h,s_2} - \bw_h \big)_{\Sigma^n_c} 
		%		-\big( \bu_{h,s_2}^{n}  - \vel_h^{n} , \stress(\bv_{h,s_2},-q_{h,s_2}))\n_{s_2} \big)_{\Sigma^n_c} 
		%		\\
		%		+ \frac{\gamma \mu}{h} \big(\bu_{h,s_2}^{n } - \vel_h^{n} ,  \bv_{h,s_2} - \bw_h\big)_{\Sigma^n_c}  \\
		%%%% pourous stresses
		+\frac{\alpha}{\sqrt{K_\tau\epsilon_{\rm p}}} (\bu_{\tau,h}^n, \bv_{h,\tau})_{\Sigma_{\rm p}}
		- (\sigma_p,   \bv_h\cdot n)_{\Sigma_{\rm p} } 
		%%%% darcy
		+ (\epsilon_{\rm p} K_\tau \nabla_\tau P_{\mathrm l,h}^n, \nabla_\tau q_{\mathrm l,h})_{\Sigma_p}
		- \big( \bu^{n}_{h,n}, q_{\mathrm l,h}\big)_{\Sigma_{\rm p} } 
		= 0
		\end{multline}
		for all $(\bv_h , q_h, \bw_h,q_{\mathrm l,h} ) \in \bV_h^n  \times Q_h^n \times \bW_h^n \times \bs S_h  $, with the porous stress $\sigma_{\rm p}$ as defined in Algorithm~\ref{alg:FSIDarcy} and the contact term $\tilde{P}_{\gamma_c}$ defined in~\eqref{Pgammatilde}.
%		\begin{align*}
%		\sigma_p
%		%=- P^n_{\rm l}  {\sfrei -} \frac{\epsilon_{\rm p} K_n^{-1} }{4} u_{l,n}|_{\gamma^+} 
%		= - P^n_{\rm l}  {\sfrei -} \frac{\epsilon_{\rm p} K_n^{-1} }{4}  \bu^n_{h,i} \cdot \bn  \quad \text{ on } \Sigma_p .
%		\end{align*}
	\end{enumerate}
	\caption{Strongly coupled Eulerian approach for an FSI-contact problem with a thick-walled solid.\label{alg:thick}}
\end{algorithm}

Finally, let us note that due to the different meaning of the displacement $\bd_h^n$ in the current frame $\Omega^{{\rm s},n}$ (see \eqref{EulerianDispl}), the no-penetration condition in \eqref{eq:contact_ineq2} becomes 
\begin{align*}
\bd_h^n - \big(\bd_h^{n-1}\circ \Psi\big) - \tilde{g}_\epsilon^n \leq 0
\end{align*} 
where $\tilde{g}_\epsilon^n$ denotes the current distance to the lower wall $\Sigma_{\rm p}$ minus $\epsilon$ and $\Psi$ is a map from $\Sigma^{\rm n}$ to $\Sigma^{\rm n-1}$. 
The contact term takes the form 
\begin{align}\label{Pgammatilde}
\tilde{P}_{\gamma_C}(\bd_h^n) = \bd_h^n - \big(\bd_h^{n-1}\circ \Psi\big) - \tilde{g}_\epsilon^n - \gamma_C^{-1}\jump{\sigma_n},
\end{align} 
where $\gamma_C = \frac{\gamma_c E}{h}, \, E$ denotes the elasticity modulus of the solid and $\gamma_c$ is a (dimensionless) contact parameter. We note that in contrast to the ${\cal O}(h^{-2})$-weighting in the thin case, a weighting of ${\cal O}(h^{-1})$ is needed here~\cite{ChoulyHild12}.
The resulting numerical method is reported in Algorithm~\ref{alg:thick}.

\begin{remark}{(Stability)}
In~\cite{BurmanFernandezFrei2020} a stability result has been derived for a very similar variational formulation with slip- or no-slip boundary conditions on $\Sigma_p$ instead of the porous medium. An analogous result can easily be
shown for both the variational formulations in~\eqref{eq:algofsiDarcy} and~\eqref{varf:thick} using the same technique.
\end{remark}

\section{Numerical experiments}\label{sec:num}

In this section, we present different numerical examples to investigate the properties and the capabilities of the numerical approaches. First,  we investigate the fluid-porous coupling by considering two disconnected fluid reservoirs that are connected 
through a thin-walled porous media in Section~\ref{sec:numstda}. Then, we investigate the full fluid-structure-porous-contact problem for a thin-walled solid by means of a deflected thin elastic valve in Section~\ref{sec:num.thin}. As introduced in Section~\ref{sec.immersed}, we use an unfitted discretisation and solve the solid equations on the reference domain $\Omega^{\rm s}=\Sigma$. Finally, we investigate the case of a thick solid in Section~\ref{subsec:num.thick}, namely an elastic ball that falls down and bounces within a viscous fluid. Here, a Fully Eulerian approach is used in combination with the locally fitted finite element method, as described in Section~\ref{sec.Eulerian}.

\subsection{Reservoirs connected via porous layer}\label{sec:numstda}

\begin{figure}[t]
	\centering
	\includegraphics[width=0.4\textwidth]{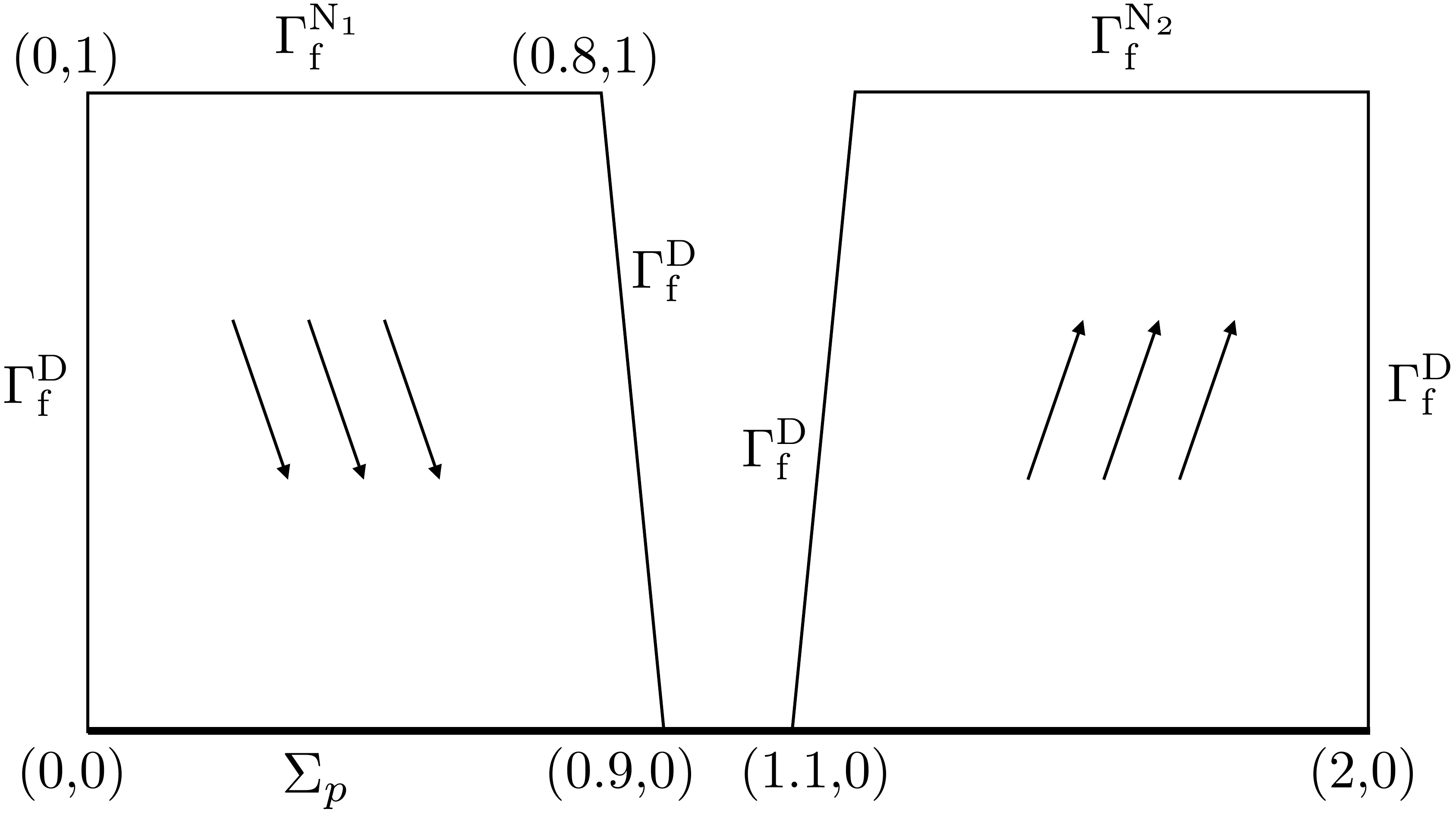}
	\caption{\label{fig.nofsi} Geometrical configuration for the Stokes model with a thin-walled porous medium on the bottom wall.}
\end{figure}

In this example, we consider two disconnected fluid reservoirs, 
connected through a thin-walled porous interface located on the bottom wall $\Sigma_{\rm p}$, as shown in Figure~\ref{fig.nofsi}. 
The fluid domain $\Omega^{\rm f}$ is shown in Figure~\ref{fig.nofsi} and $\Sigma_{\rm p}$ is a segment with extremities $(0,0)$ and $(2,0)$.

\begin{figure}[ht]
	\centering
	\includegraphics[width=0.5\textwidth]{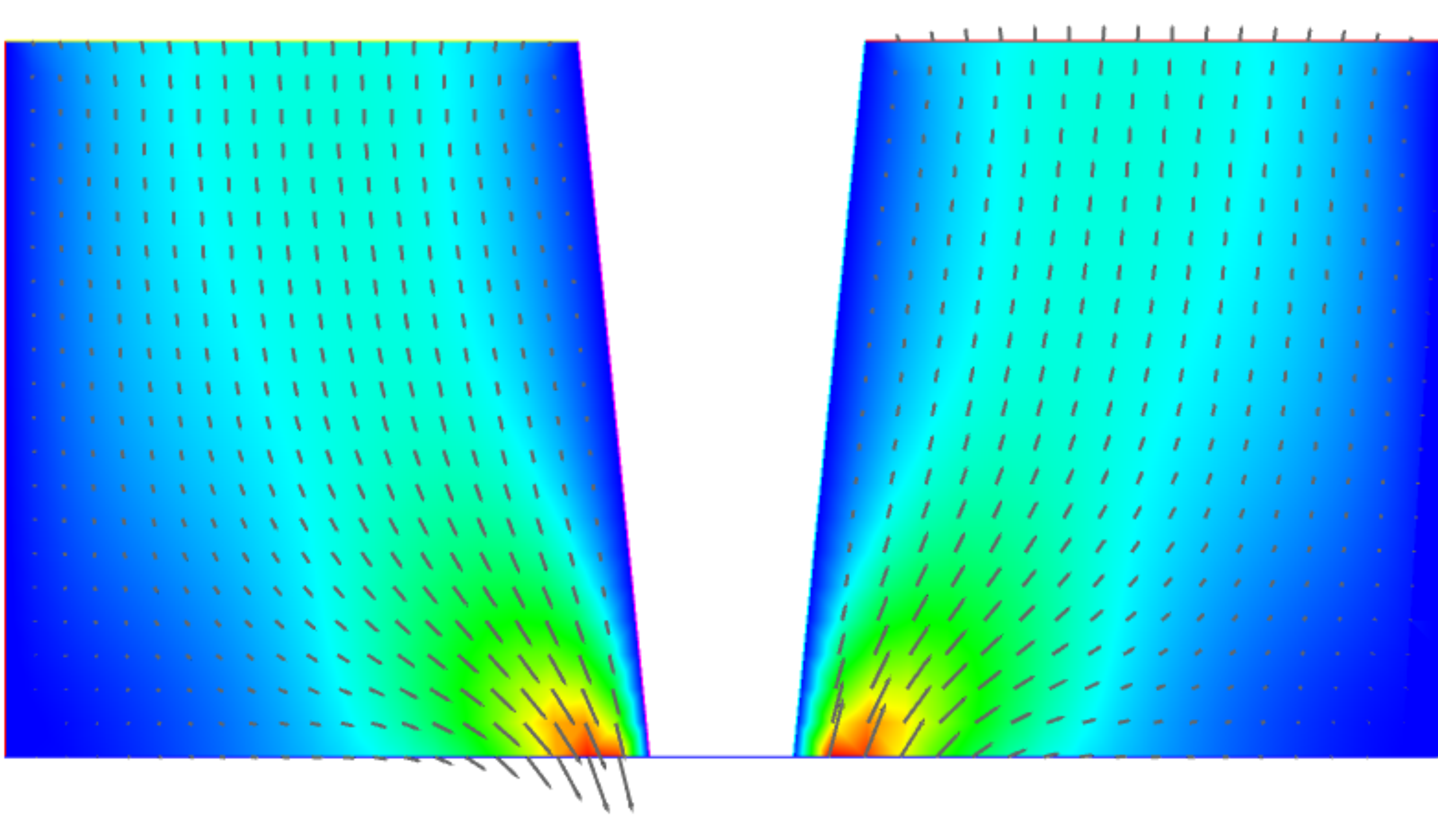}\\
	\includegraphics[width=0.2\textwidth]{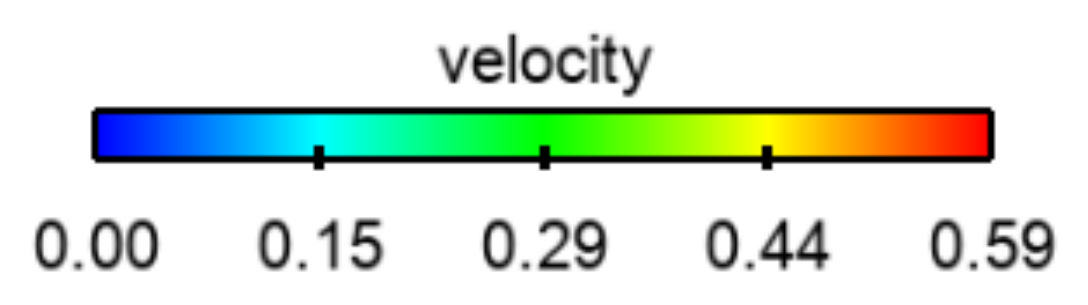}
	\caption{Snapshot of the fluid velocity at time $t=0.25$.}
	\label{fig:nofsi-vel}
\end{figure}

\begin{figure}[ht]
	\centering
	\subfigure[Elevation of the fluid pressure]{\includegraphics[width=0.5\linewidth]{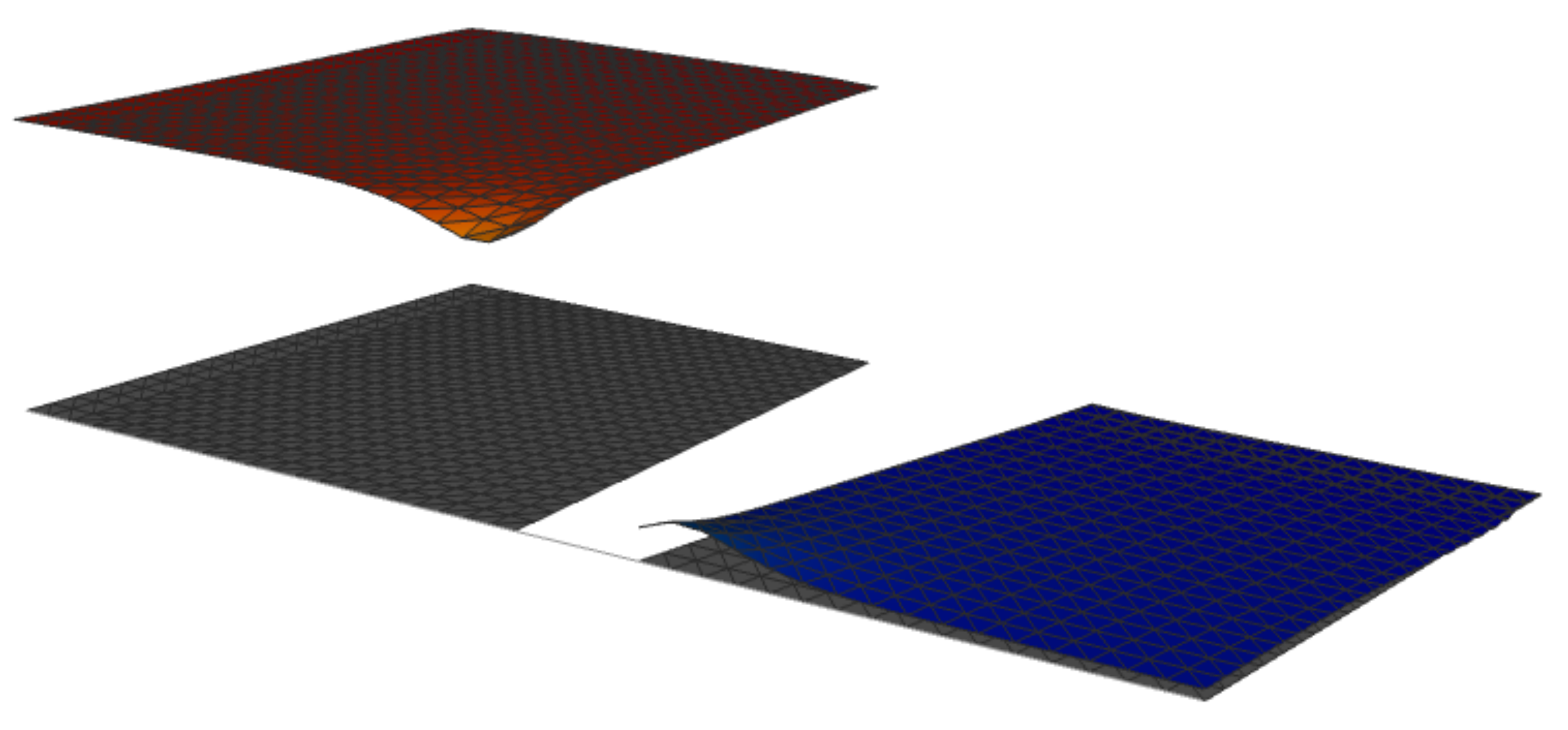}}
	\subfigure[Porous pressure]{\includegraphics[width=0.4\linewidth]{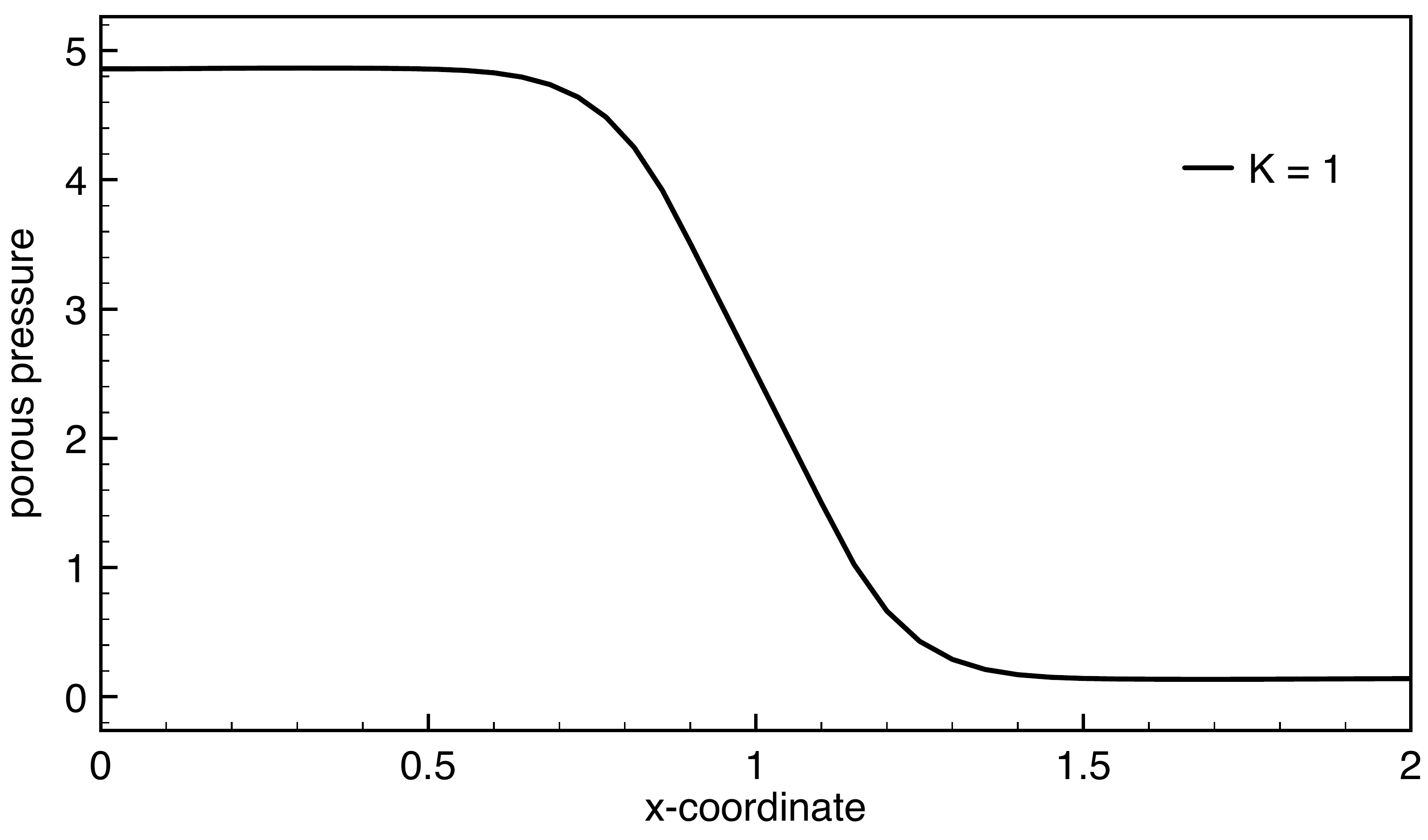}}\\
	\includegraphics[width=0.22\linewidth]{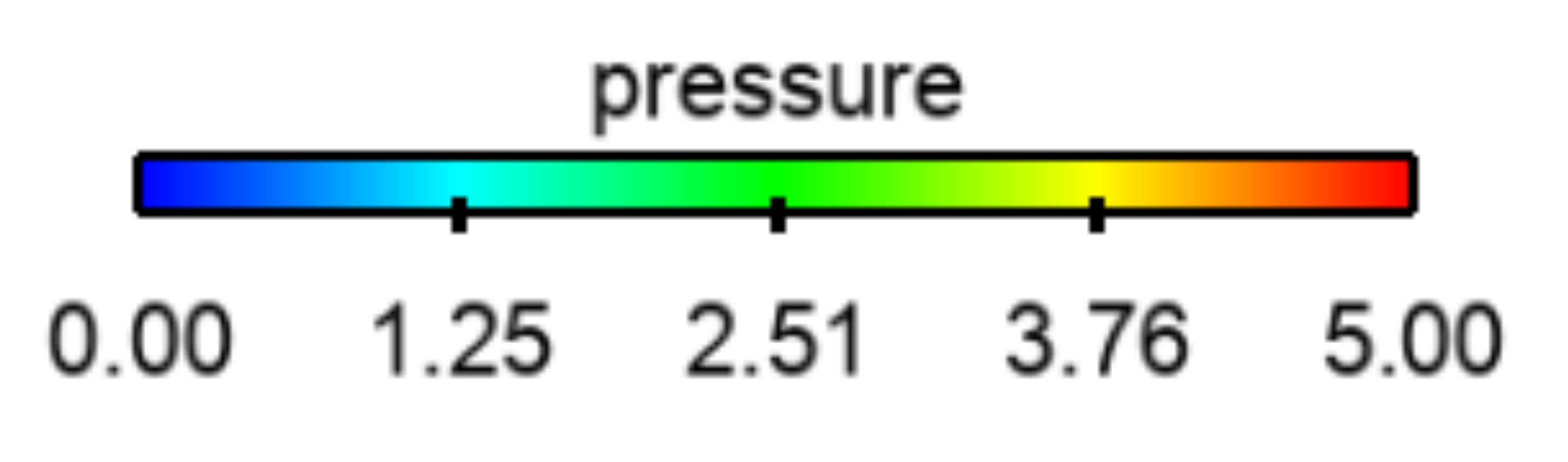}
	\caption{Fluid and porous pressures at $t=0.25$. }
	\label{fig:nofsi-press}
\end{figure}

Regarding the fluid boundary conditions, we impose a pressure drop across the two parts of the top boundary. A traction is imposed on $\Gamma^{\rm N_1}_{\rm f}$ in terms of a sinusoidal time-dependent pressure $p_{\rm {in}}(t) $, namely,  $$p_{\rm {in}}(t) = 5 \cdot  \sin(2\,\pi \,t), \qquad \forall \, t \in \mathbb{R}^+,$$ 
while a zero traction is enforced on $\Gamma^{\rm N_2}_{\rm f}$.
Additionally, a no-slip boundary condition is enforced on $\Gamma^{\rm D}_{\rm f}$. 
The considered physical parameters are $\mu=0.03$, $\rho^{\rm f}=1$, $\varepsilon_p = 0.01$ and $K_\tau=K_n = 1$. We consider an approximation of the Beavers-Joseph-Saffmann condition for the tangential stresses, in which we let $\alpha = 0$.

The purpose of this example is to illustrate how the porous model is able to connect the fluid flow between the two containers. This can be clearly inferred from the results reported in Figure~\ref{fig:nofsi-vel} and Figure~\ref{fig:nofsi-press} at $t=0.25$, which, respectively, show a snapshot of the fluid velocity, the elevation of the fluid pressure and the associated porous pressure. As we can see, the fluid is entering into the porous layer from the left reservoir and leaving the porous interface into the right one. 
%Further studies on the model parameters will be considered in the next numerical examples.
%%%%%%%

%%%%%%%%%%%%%%%%%%%%%%%%%%

\subsection{Idealised valve with contact}
\label{sec:num.thin}
\begin{figure}[t]
	\centering
	\subfigure[]{ \includegraphics[width=0.45\linewidth]{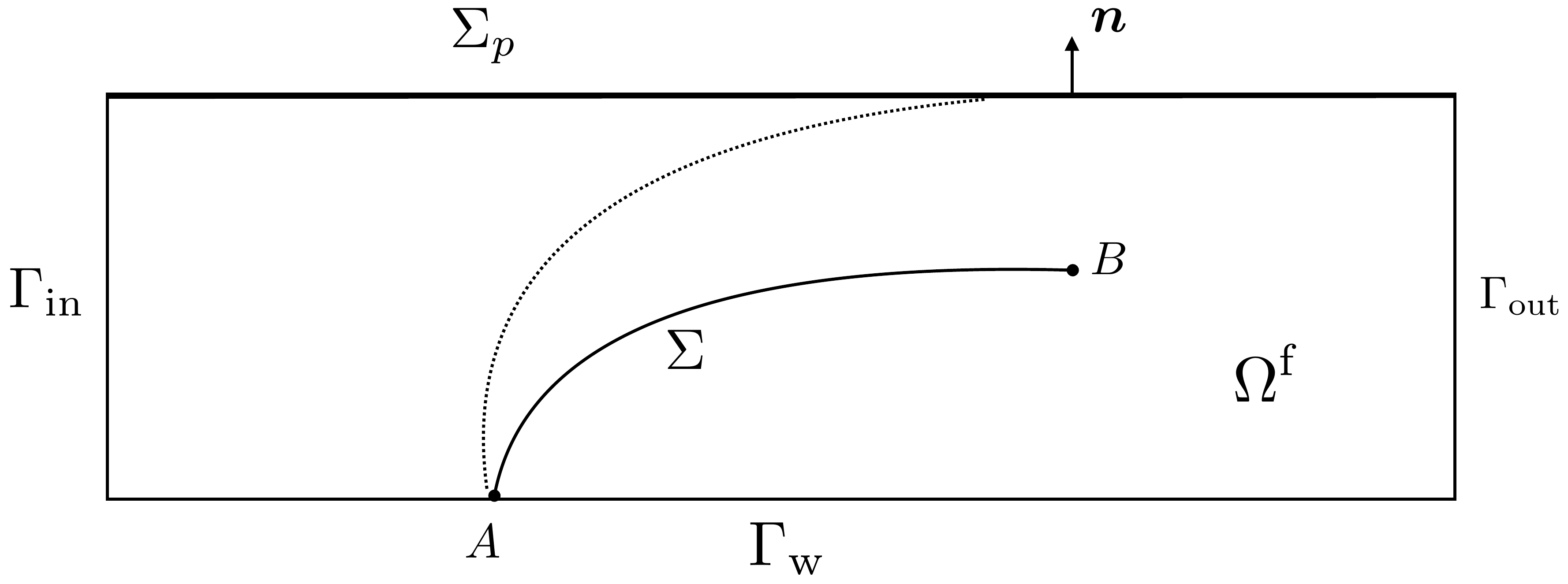} } \qquad
	\subfigure[]{\includegraphics[width=0.45\linewidth]{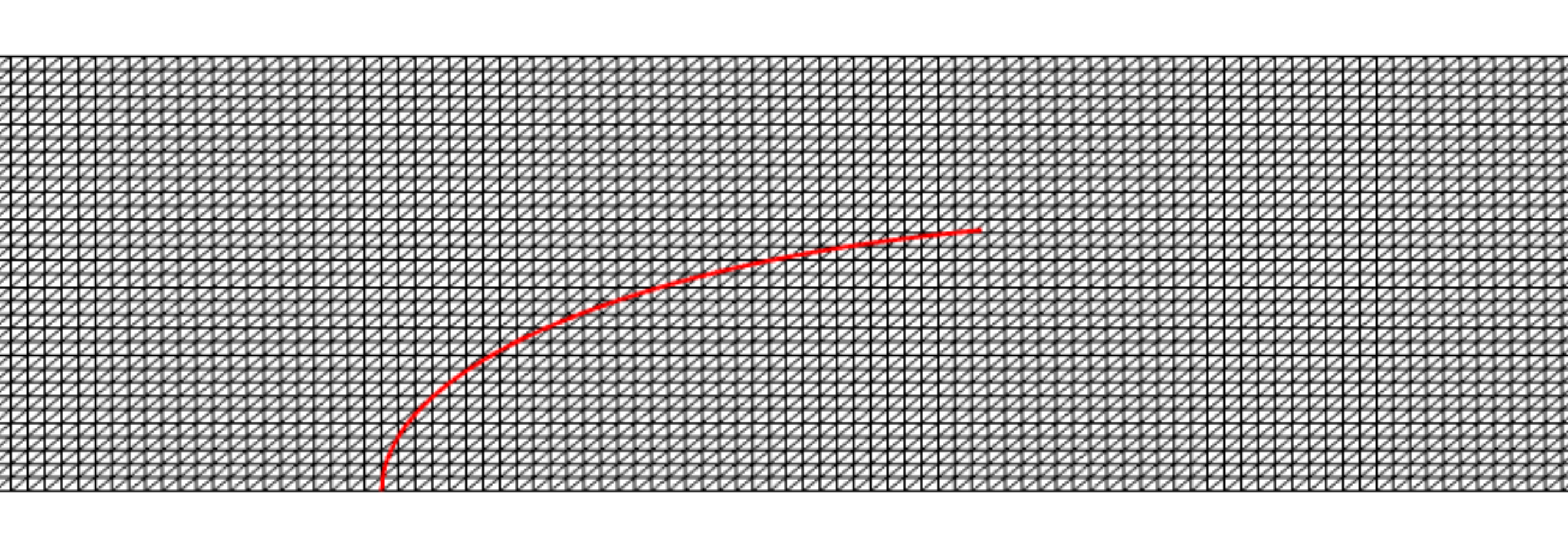} }
	\caption{(a) Geometric configuration of the idealised valve with contact, (b) Zoom of the leaflet mesh and fluid mesh.  }
	\label{fig:contact_geom}
\end{figure}

In this test, we consider a full FSI-contact problem with a thin-walled solid. 
% intro and geometry
This numerical example corresponds to the idealised valve test with possible contact on the porous layer $\Sigma_p$. The geometry is shown in Figure \ref{fig:contact_geom}(a).
The computational domain is a rectangle $\Omega = [0, 8] \times [0, 0.805]$, where the upper boundary is a symmetry axis (we imagine a second symmetrical valve on top), which means that we impose the "slip" condition $\sigma_{f,\tau} = 0$ in \eqref{eq:darcy}, letting $\alpha = 0$.
As reference configuration for the solid, $\Sigma$, we consider a curve segment with end points $A=(4,0)$ and $B=(5.112, 0.483)$, parametrised by the analytical function   
$$ y(x) =  \ \frac12 \sqrt{1 - \dfrac{(x-11/2)^{2}}{(3/2)^{2}}},\quad  x \in [4, 5.112]. $$
% physical parameters
All the following units are given in the CGS units system.
The physical parameters used for the fluid in this test are $\rho^{\rm{f}} = 1$, $\mu = 0.03$. For the solid we have $\rho^{\rm{s}} = 1.2$, $\epsilon^{\rm{s}} = 0.065$, the Young\textquoteright s modulus $E = 10^7$ and the Poisson\textquoteright s ratio $\nu = 0.4$. Regarding the porous medium, we consider $\epsilon_{\rm p} = 0.01$ and we explore the influence of the porous layer on the contact dynamics by changing the hydraulic conductivity parameters $K_\tau=K_n = \{ 10^{-i}\}_{i=1}^3$ and by comparing it with the situation of a simple wall on $\Sigma_{\rm p}$, where we enforce a symmetry condition.

%%%%%%%%%%%%%%%%%%%%%%%%%%%% 
% BC
Regarding the boundary condition, a no-slip condition is enforced on the lower boundary $\Gamma_{w}$, zero traction on the outflow boundary $\Gamma_{\rm{out}}$ and a traction condition on $\Gamma_{\rm{in}}$, in terms of the following  time-dependent pressure: 
\begin{equation*}\label{eq:inpress}
p_{\rm{in}}(t) =  
\left\{
\begin{aligned}
- 200 \atanh (100  t) & \quad \text{if}\quad 0 < t < 0.7, \\
200 & \quad \text{if}\quad t \geq 0.7.
\end{aligned}\right.
\end{equation*}
The final time is $T=1$, which corresponds to one full valve oscillation cycle. 
The fluid and the solid are initially at rest and the beam is pinched at the bottom tip $A$. 
In this test, the solid is described by a non-linear Reissner$-$Mindlin curved beam model with a  MITC spacial discretisation. The ghost penalty parameter has been set to $\gamma_{\rm g} = 1$ and the CIP stabilisation parameters to $ \gamma_{\rm v} = \gamma_{\rm p} = 10^{-2}$.
In this particular test case, the gap function is defined as the initial distance of a point on $\Sigma$ to the wall $\Sigma_p$ in the direction of $\bn_{\rm l}$, namely $g = y_{\Sigma_p} -  y(x) $.
% contact
The contact parameters are given by  $\epsilon_g = 0.01$  and $\gamma_{\rm{c}} = 5 \cdot 10^{-3}$ as in \cite{boilevin-et-al-19}.
The relaxation parameter $\epsilon_g$  is chosen in such a way that the generated artificial gap is below $h$, typically $\varepsilon_h \approx h/2$. 
The penalty parameter $\gamma_{\rm c}$ (independent of $h$) is chosen to avoid penetration (i.e., not very small) and in such a way that 
the term \eqref{sigmaPgamma} does not perturb the convergence of the Newton solver in the solid (the operator $[\cdot]_+$ is not differentiable at $0$).
%The contact parameters are given by  $\epsilon_{\rm p} = 0.01$  and $\gamma_{\rm{c}} = 5 \cdot 10^{-3}$ as in \cite{boilevin-et-al-19}.
%The relaxation parameter $\epsilon_g$ is chosen in such a way that the generated artificial gap is below $h$. 
%The penalty parameter $\gamma_{\rm c}$ (independent of $h$) is chosen to avoid penetration (i.e., not very small) and in such a way that the contact penalty term does not perturb the convergence of the Newton solver in the solid (the operator $[\cdot]_+$ is not differentiable at $0$).
% descretization (mesh, ghost penalty and stab param)

The fluid mesh has $16\,384$ triangles and the solid mesh 50 edges. We have $h\approx 0.04$. The zoom on both meshes is presented in Figure \ref{fig:contact_geom}(b).
The time discretisation parameter is $\delta t = 10^{-3}$ and the Nitsche parameter is set to $\gamma = 100$.
% remaining numerical parameters are set to $\gamma = 100, \gamma_{\rm g} = 1$ and $\gamma_{\rm p} = 10^{-2}$. 

\begin{figure}[ht]
	\centering
	\subfigure[$ t = 0.6 $.]{\includegraphics[width=0.45\linewidth]{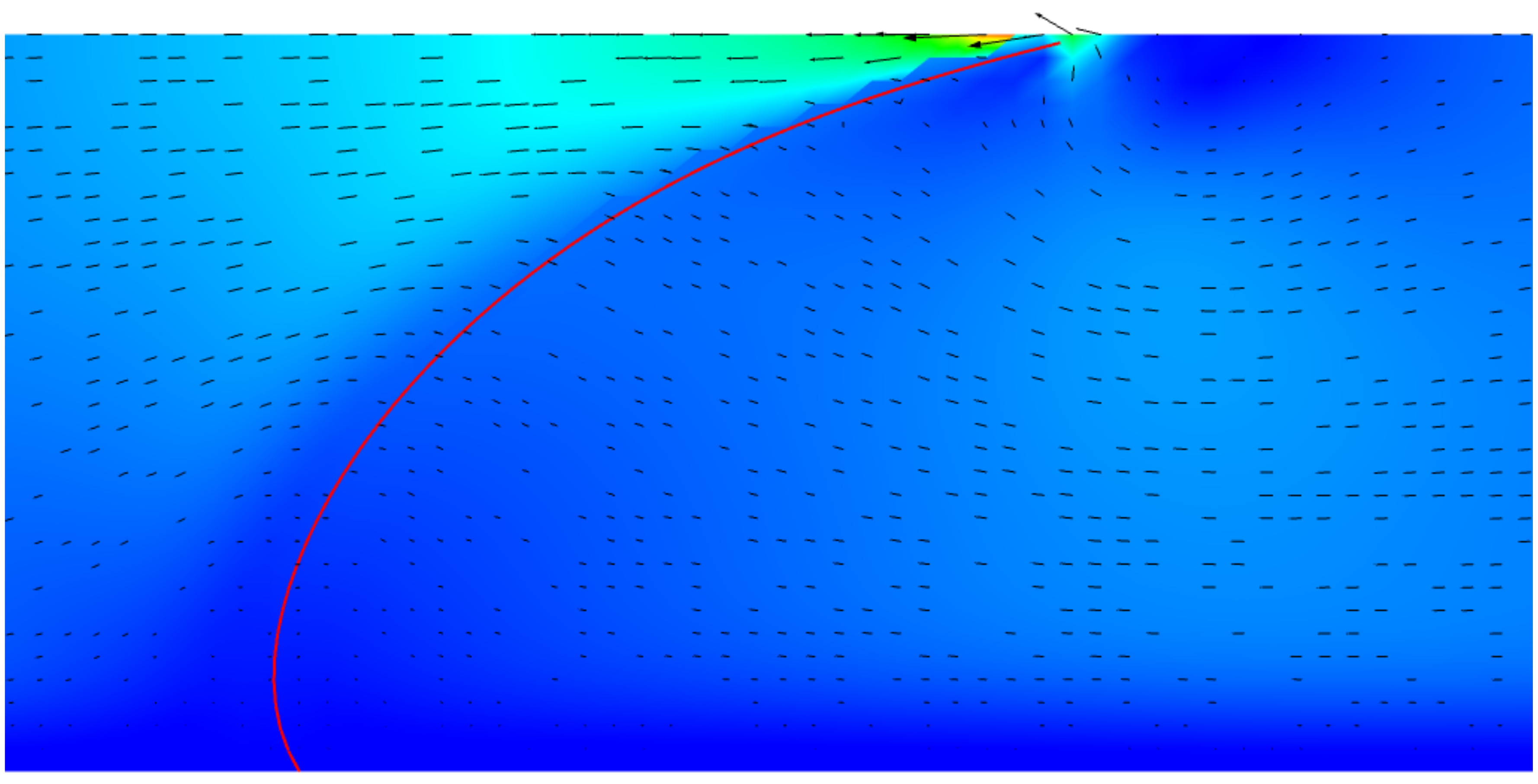}}
	\subfigure[$ t = 1 $.]{\includegraphics[width=0.45\linewidth]{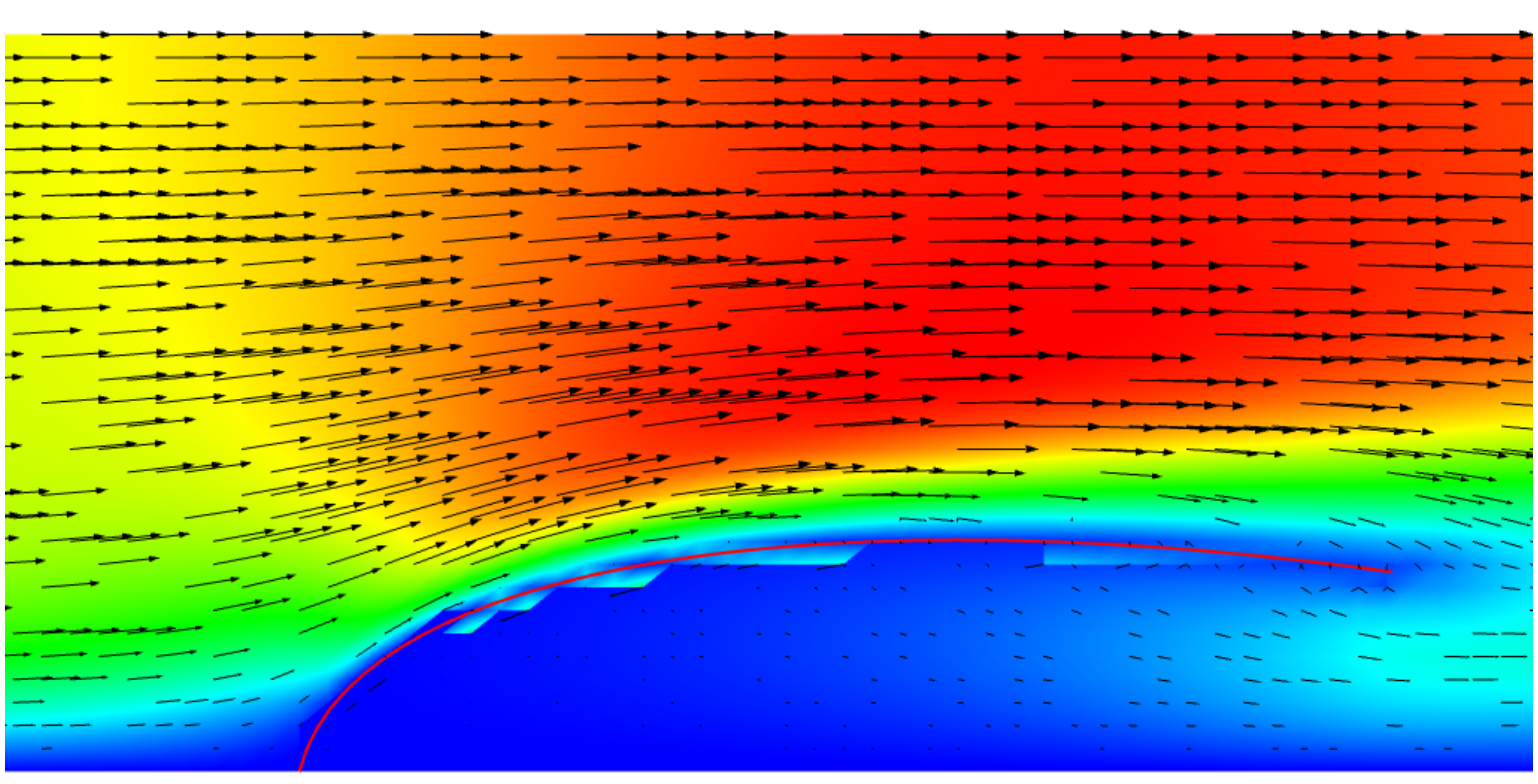}}
	\\
	\includegraphics[width=0.3\linewidth]{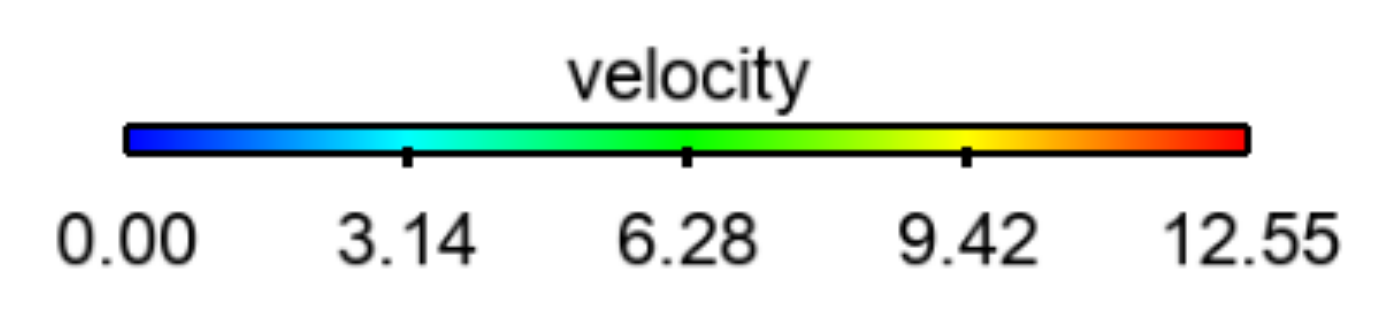}
	\caption{Velocity magnitude snapshots.}
	\label{fig:vel_contact}
\end{figure}
\begin{figure}[h!]
	\centering
	\subfigure[$ t = 0.6 $.]{\includegraphics[width=0.45\linewidth]{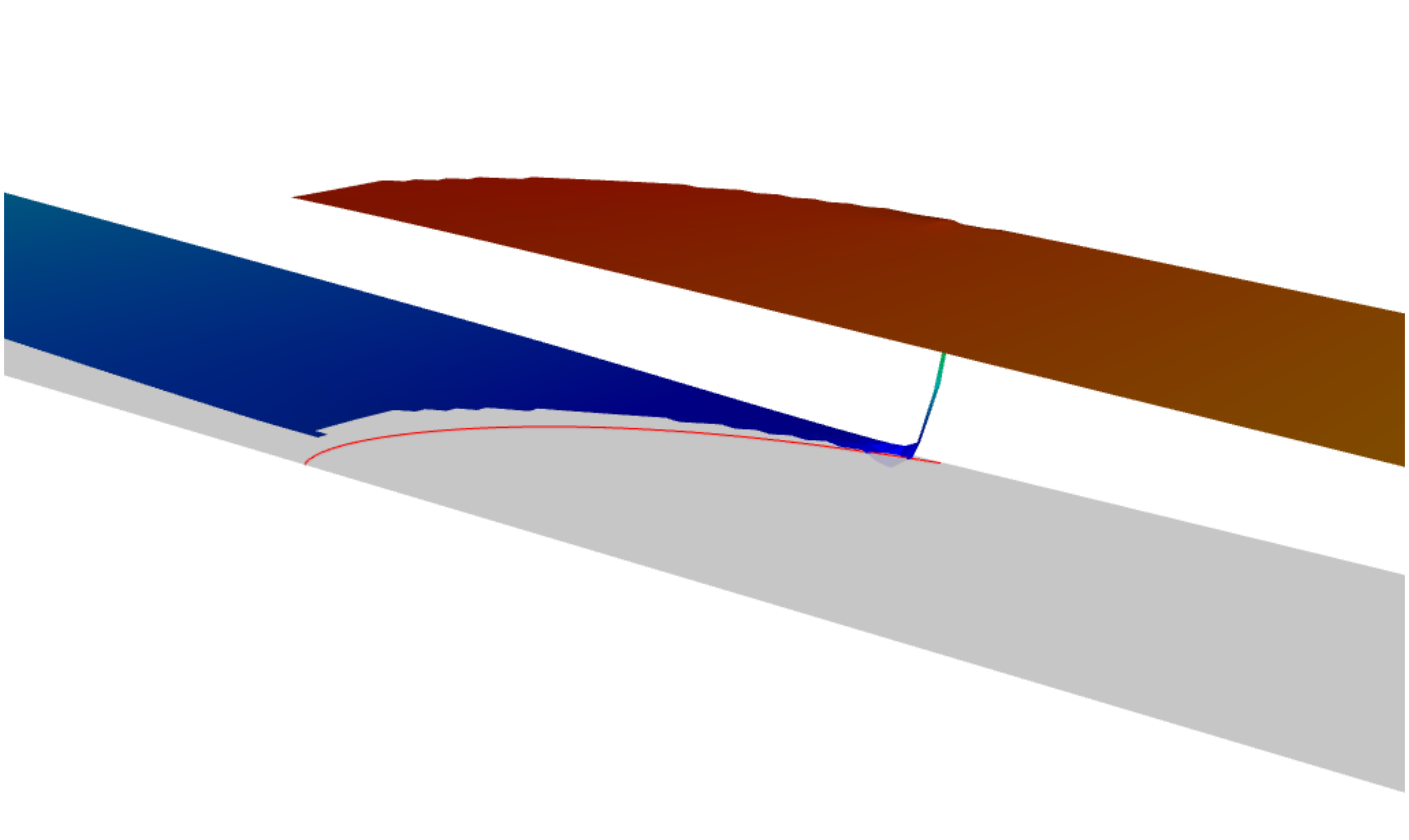}}
	\subfigure[$t = 1 $.]{\includegraphics[width=0.45\linewidth]{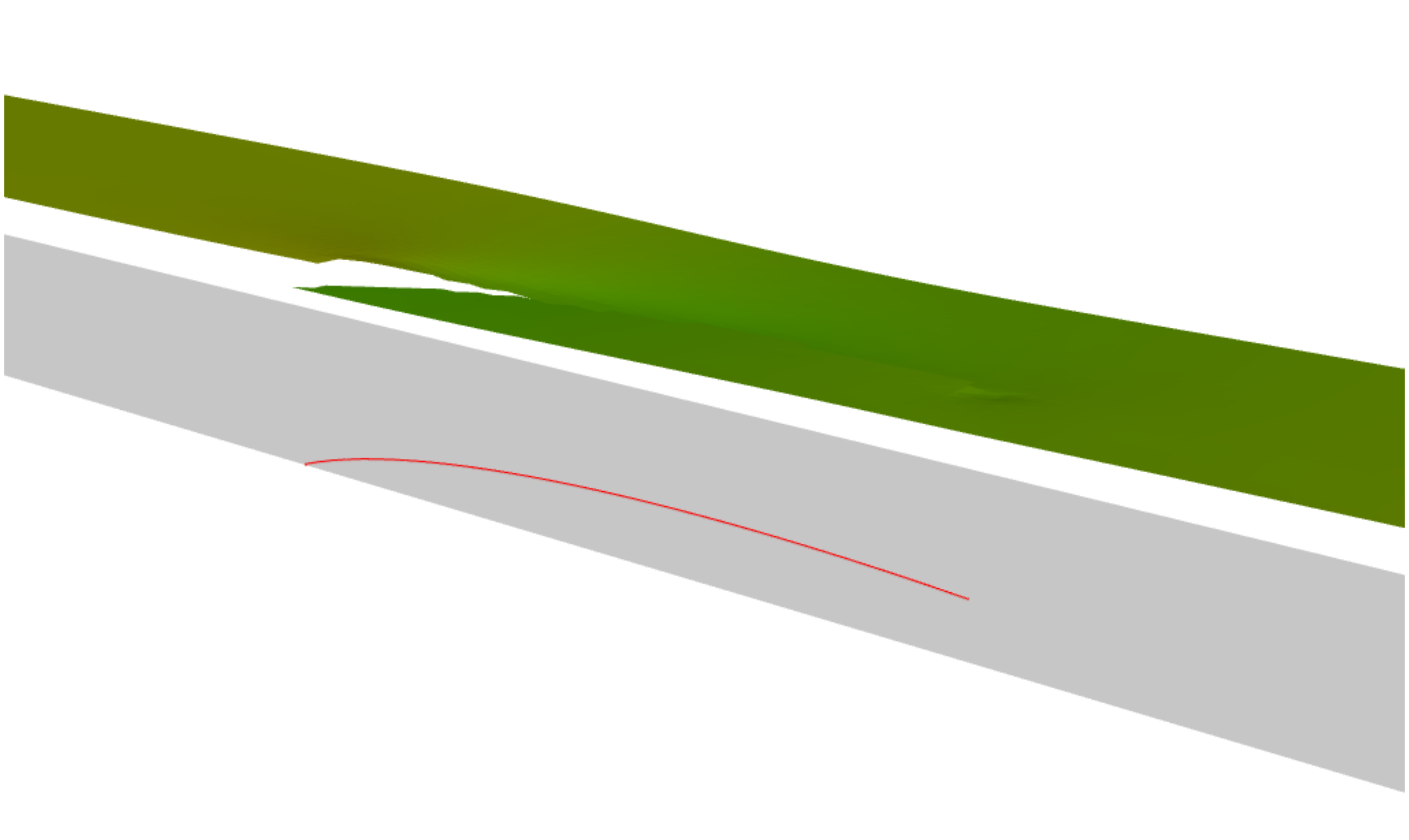}} 
	\\
	\includegraphics[width=0.23\linewidth]{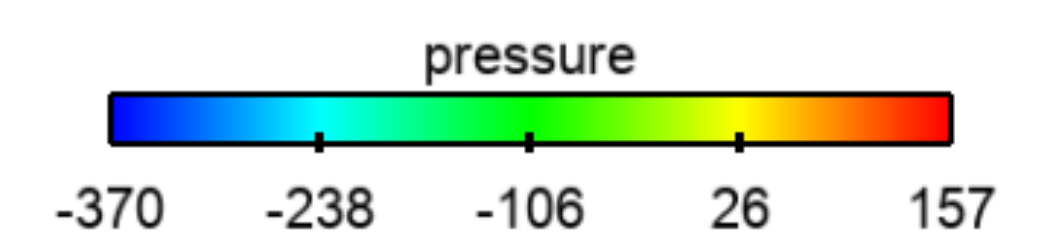}
	\caption{Pressure elevation snapshots.}
	\label{fig:press_contact}
\end{figure}

\begin{figure}[ht]
	\centering
	\subfigure[$x$-displacement.]{\includegraphics[width=0.45\linewidth]{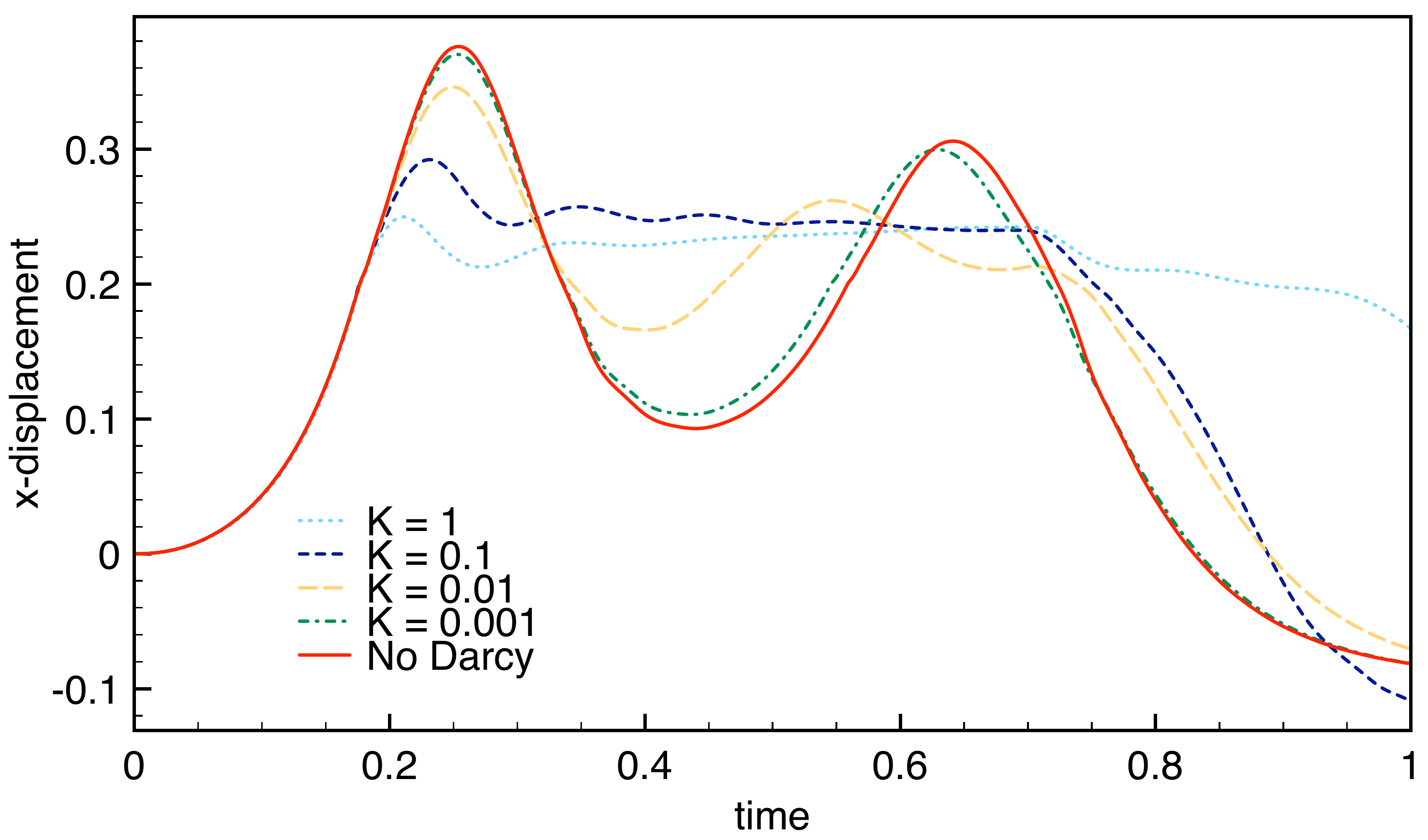}}
	\subfigure[$y$-displacement.]{\includegraphics[width=0.45\linewidth]{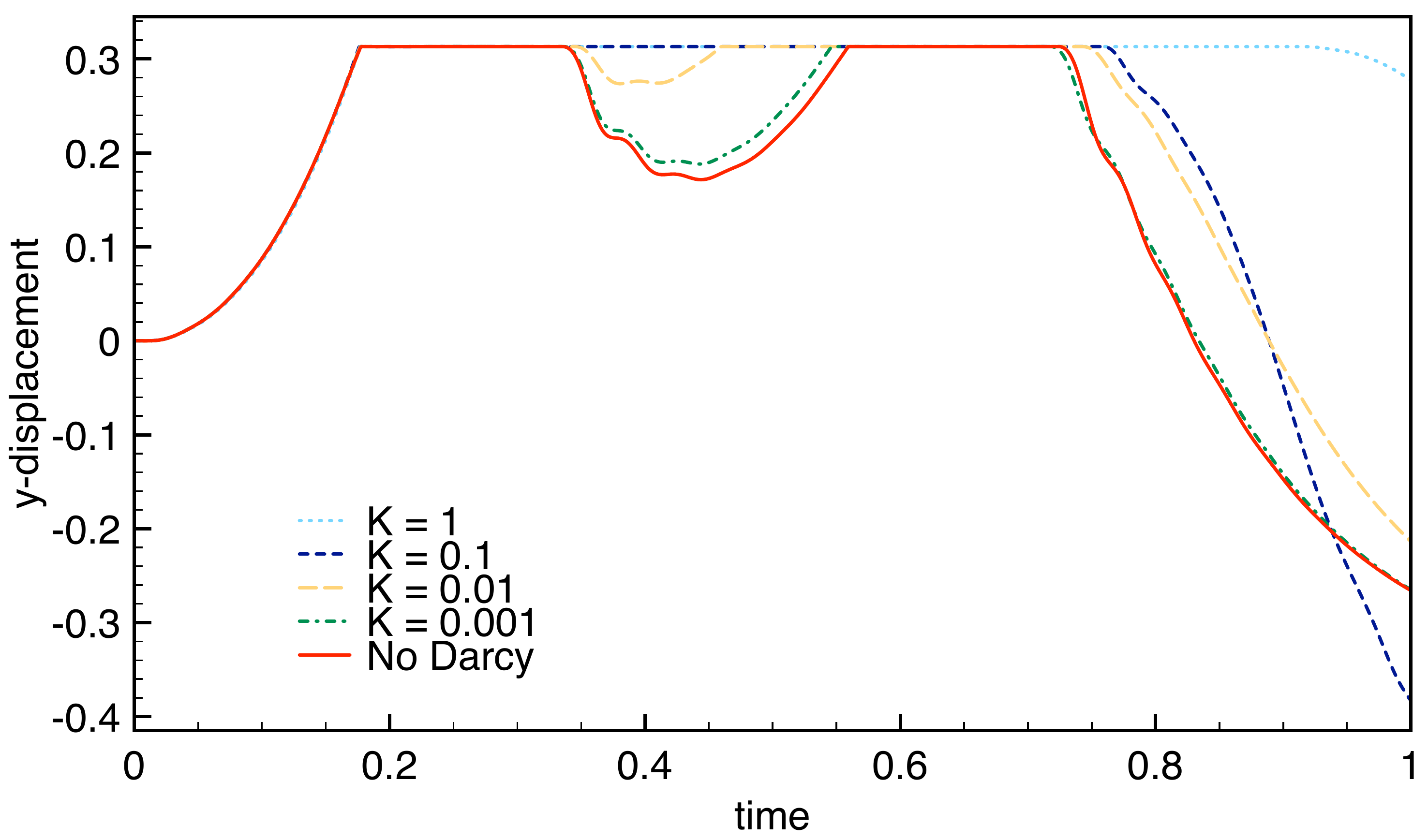}}
	\caption{Time evolution of the $x$ and $y$-displacement for the structure endpoint $B$.}
	\label{fig:dispxy}
\end{figure}
\begin{figure}[ht]
	\centering
	\subfigure[$ t = 0.25$.]{\includegraphics[width=0.45\linewidth]{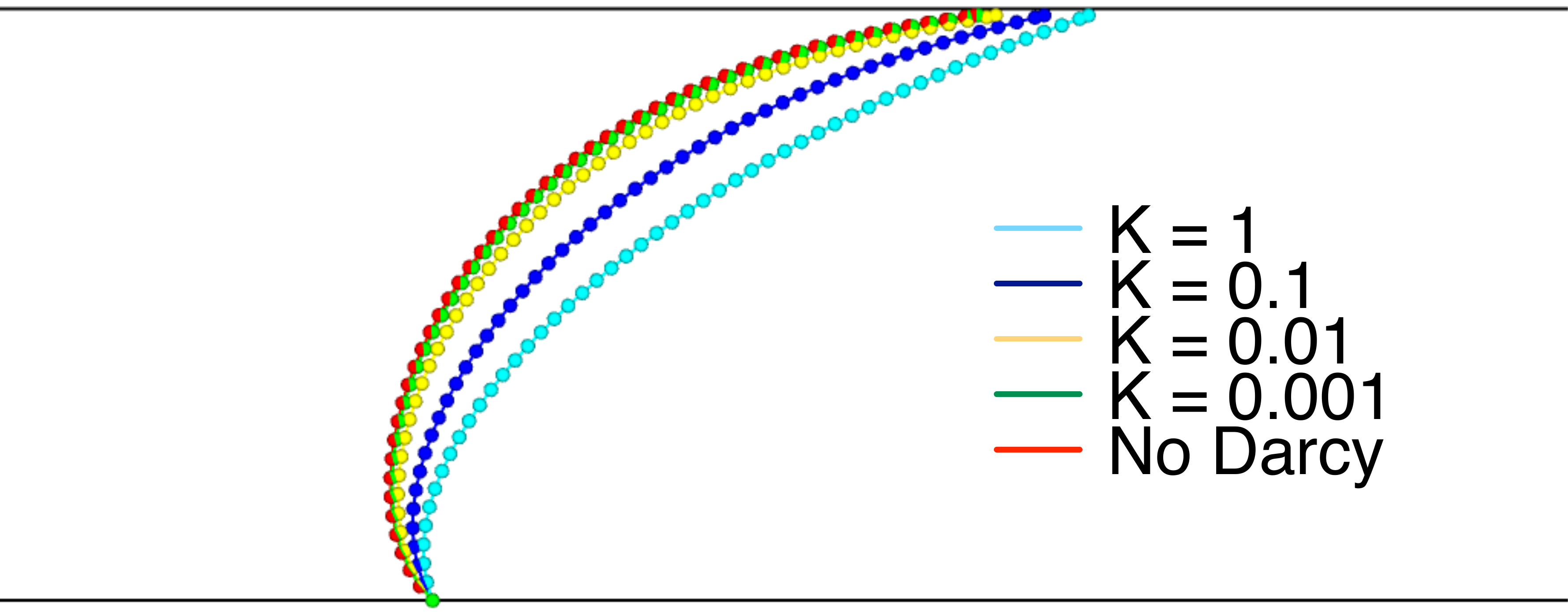}}
	\subfigure[$ t = 0.45 $.]{\includegraphics[width=0.45\linewidth]{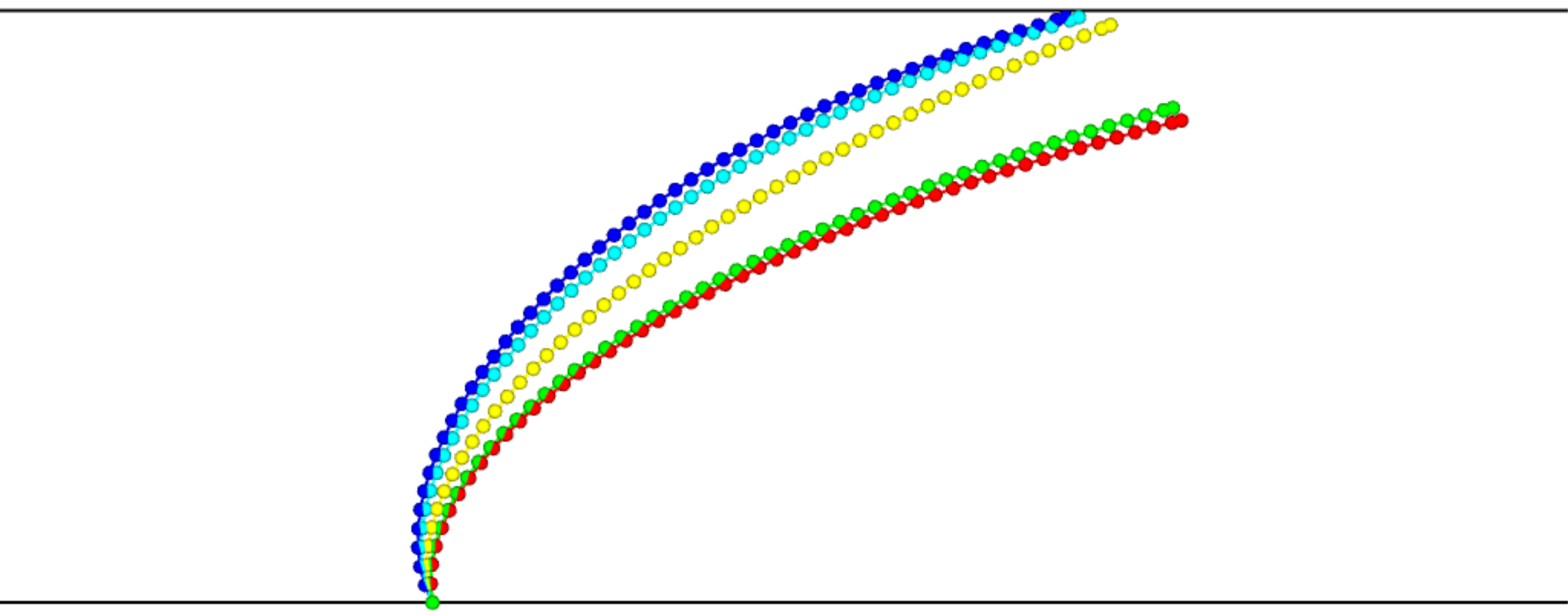}}
	\subfigure[$ t = 1 $.]{\includegraphics[width=0.45\linewidth]{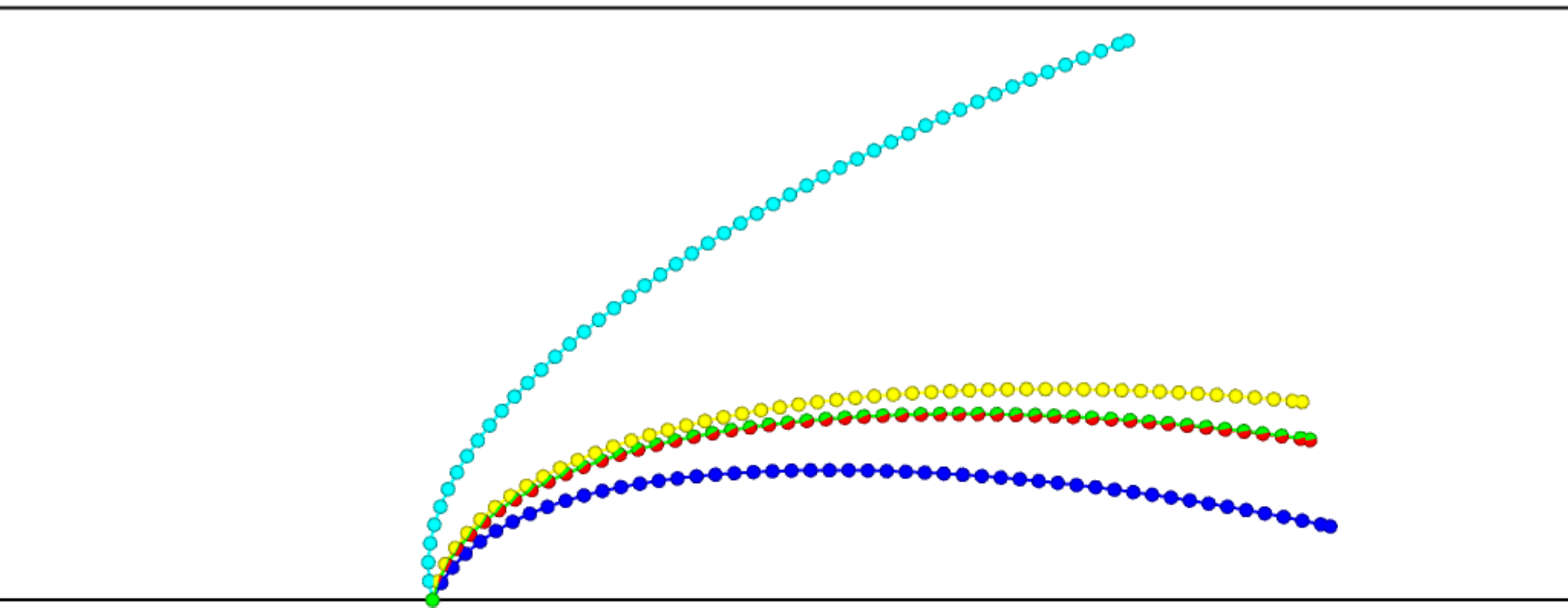}}
	\caption{Interfaces location at time $ t = 0.25$ (a), $ t = 0.45$ (b) and $ t = 1$(c).}
	\label{fig:interfacePos}
\end{figure}

%
%\begin{figure}[h!]
%	\centering
%	\subfigure[$ t = 0.6 $.]{\includegraphics[width=0.45\linewidth]{figures/velocity-close}}
%	\subfigure[$ t = 0.6 $.]{\includegraphics[width=0.4\linewidth]{figures/pressure-close}}
%	\\
%	\subfigure[$ t = 1 $.]{\includegraphics[width=0.45\linewidth]{figures/velocity-open}}
%	\subfigure[$t = 1 $.]{\includegraphics[width=0.4\linewidth]{figures/pressure-open}} 
%	\\
%	\includegraphics[width=0.25\linewidth]{figures/velocity_legend}
%	\includegraphics[width=0.25\linewidth]{figures/pressure_legend}
%	\caption{Velocity magnitude (a),(c) and pressure elevation (b),(d) at time $t=0.6$ and $t=1$, respectively.}
%	\label{fig:velpress_contact}
%\end{figure}
%\begin{figure}[h!]
%	\centering
%	\subfigure[$x$-displacement.]{\includegraphics[width=0.4\linewidth]{figures/dispX}}
%	\subfigure[$y$-displacement.]{\includegraphics[width=0.4\linewidth]{figures/dispY}}
%	\caption{Time evolution of the $x$ and $y$-displacement for the structure endpoint $B$.}
%	\label{fig:dispxy}
%\end{figure}
%\begin{figure}[h!]
%	\centering
%	\subfigure[$ t = 0.25$.]{\includegraphics[width=0.3\linewidth]{figures/solids025}}
%	\subfigure[$ t = 0.45 $.]{\includegraphics[width=0.3\linewidth]{figures/solids045}}
%	\subfigure[$ t = 1 $.]{\includegraphics[width=0.3\linewidth]{figures/solids1}}
%	\caption{Interfaces location at time $ t = 0.25$, $ t = 0.45$ and $ t = 1$.}
%	\label{fig:interfacePos}
%\end{figure}

Let us first consider a test case with $K= K_\tau = K_n = 10^{-3}$. We report in Figure~\ref{fig:vel_contact} the velocity magnitude at two different instants. In Figure~\ref{fig:vel_contact}(a) we report the approximation obtained at time $t = 0.6$. At this instant, the valve is in contact with the upper wall and the fluid velocity decreases globally as a consequence of the closing. Contrarily to the idealised valve test without porous layer at the top wall, here, we allow the flow to enter the porous interface at contact. The fluid is transported through the porous layer, from the right side of the domain to the left side. 
At $t =1$ the valve is open and far from $\Sigma_p$, therefore the fluid flow is reestablished and the velocity increases in the channel. 
The same comparison is performed in Figures~\ref{fig:press_contact}(a) and (b) for the pressure.
We can see the high pressure jump when the valve is in contact with the wall (Figure~\ref{fig:press_contact}(a)), while at $t =1$ the discontinuity between the two sides of the interface is weaker (Figure~\ref{fig:press_contact}(b)). 
\begin{figure}[h!]
	\centering
	\subfigure[$ t = 0.25$.]{\includegraphics[width=0.45\linewidth]{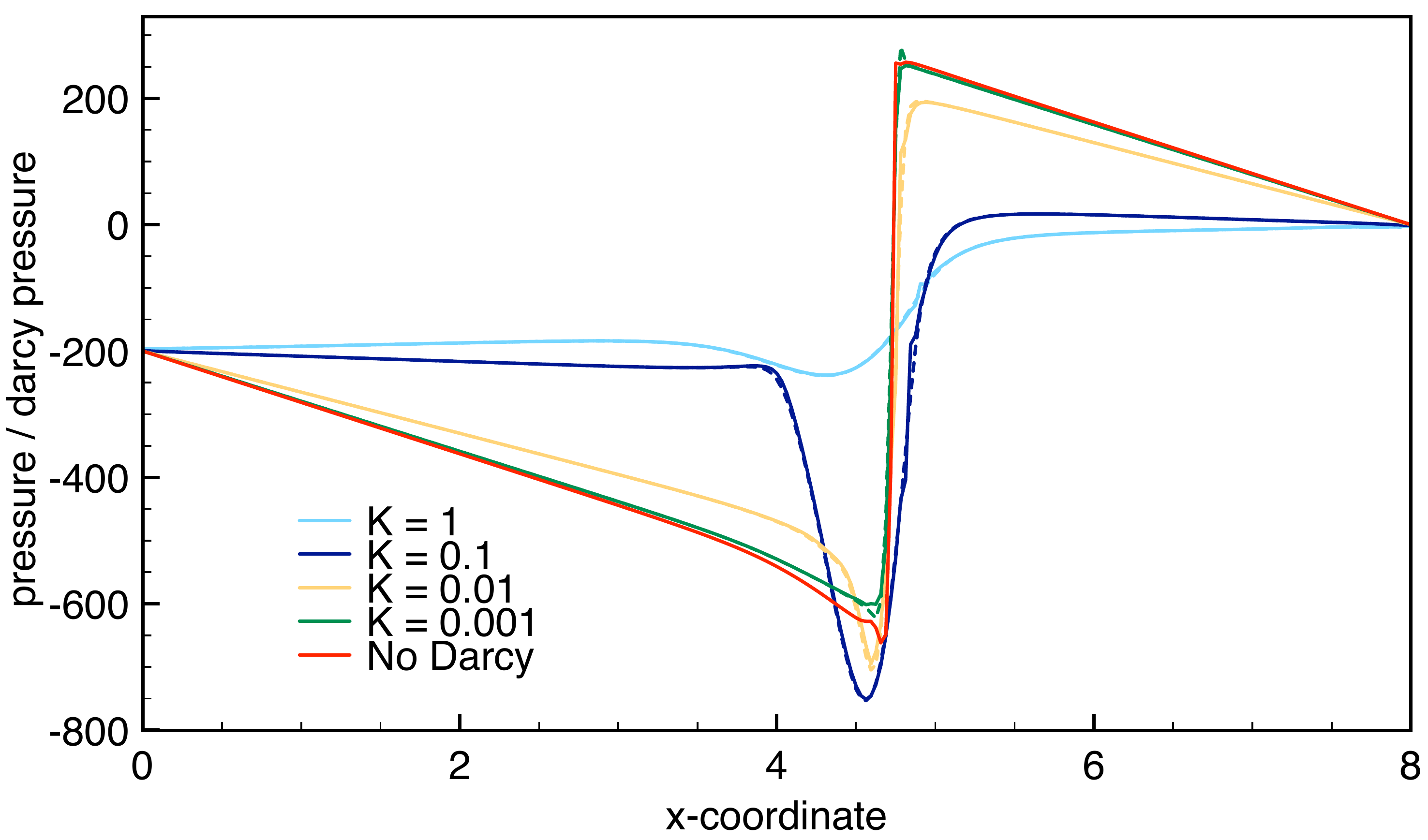}}
	\subfigure[$ t = 0.45 $.]{\includegraphics[width=0.45\linewidth]{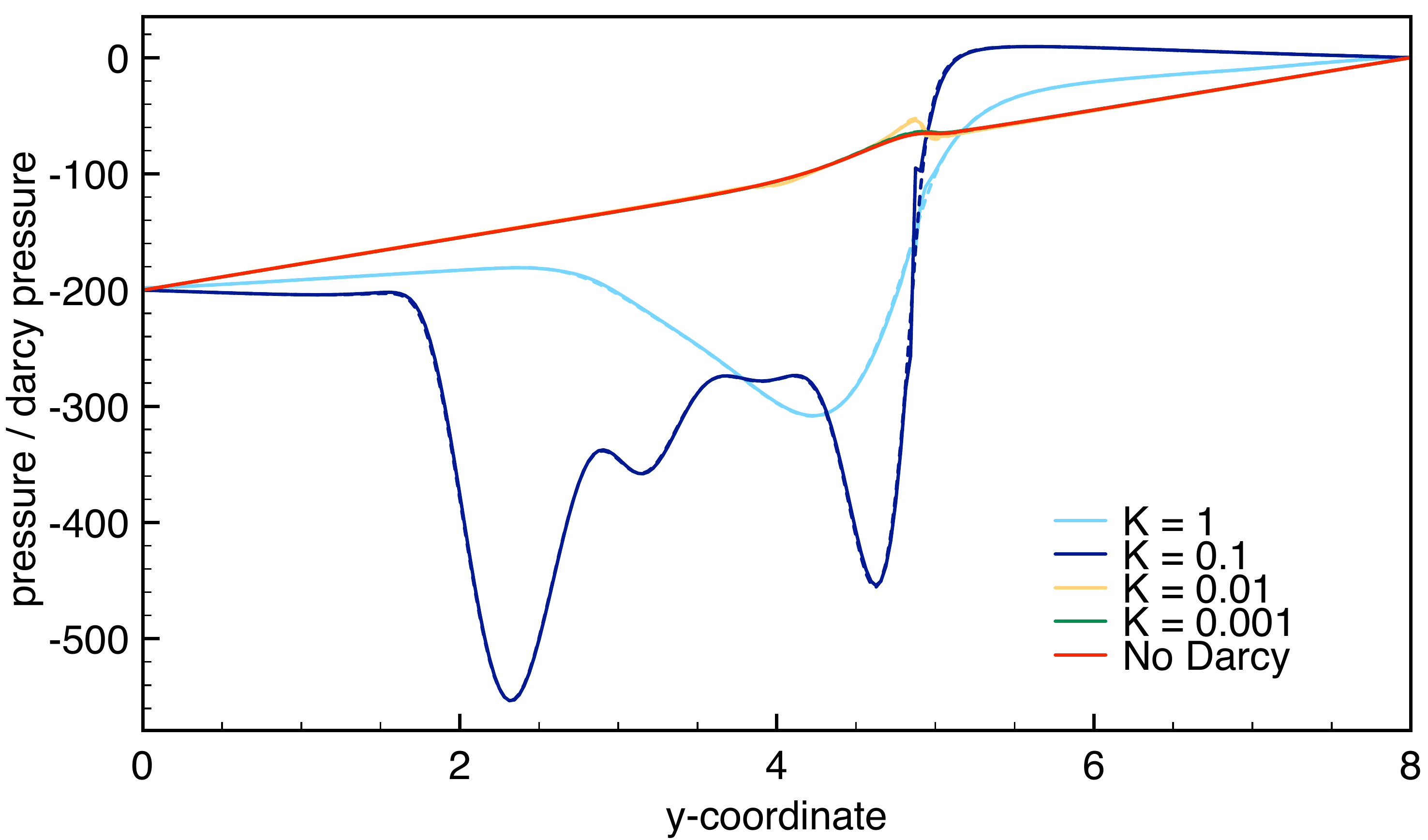}}
	\caption{Fluid pressure (continuous line) and porous pressure (dash line) on $\Sigma_p$ for different value of hydraulic conductivity at time $t=0.25$ and $t=0.45$.}
	\label{fig:fluidpressure}
\end{figure}

	We now consider the case in which we insert a thin porous layer on the top contact wall and investigate the impact of $ K_\tau = K_n$ on the results. 
	\begin{figure}[ht]
		\centering
		\subfigure[x-velocity at $ t = 0.25$.]{\includegraphics[width=0.45\linewidth]{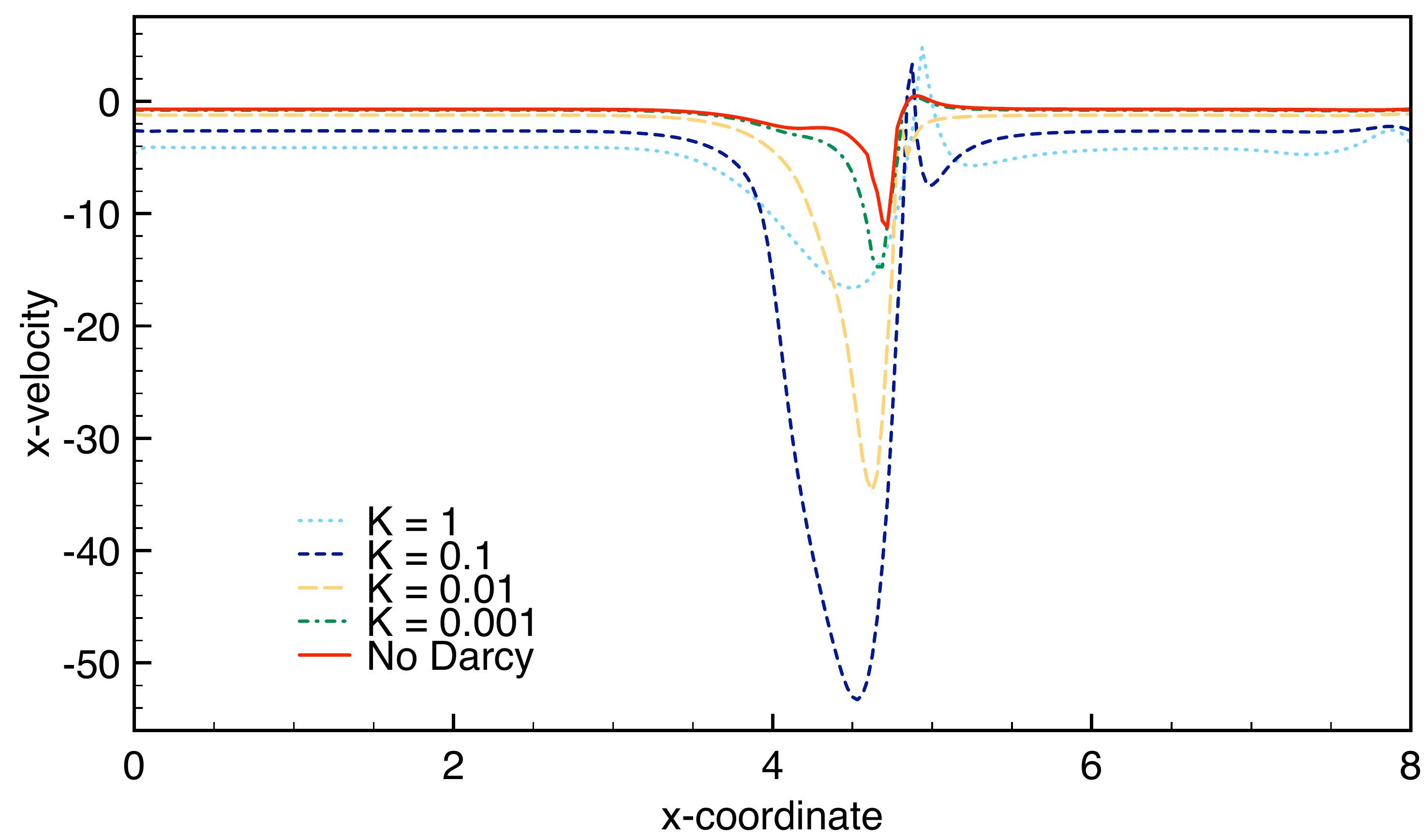}}
		\subfigure[x-velocity at $ t = 0.45 $.]{\includegraphics[width=0.45\linewidth]{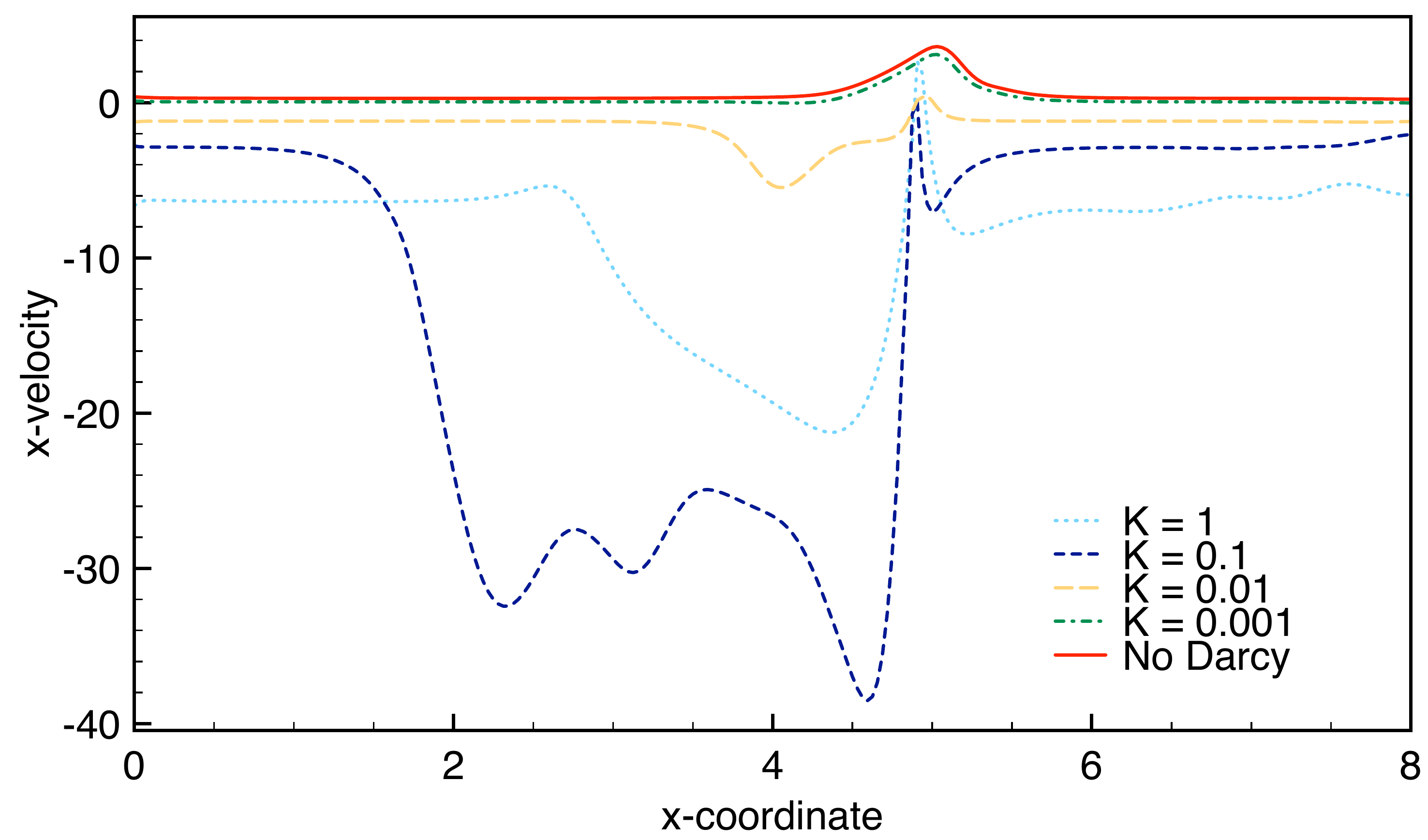}}
		\caption{Fluid velocity on $\Sigma_p.$}
		\label{fig:fluidXVelSigmaP}
	\end{figure}
	Figure~\ref{fig:dispxy} presents the time history of the horizontal, Figure~\ref{fig:dispxy}~(a), and vertical displacement, Figure~\ref{fig:dispxy}~(b), at the upper solid point $B$ for different levels of conductivity. 
	The non-penetration condition with the wall can be seen in Figure~\ref{fig:dispxy}(b), whereas Figure~\ref{fig:dispxy}(a) shows that the structure is sliding over the top wall. The interface is bouncing for all tests except the cases of $K_\tau = K_n = 1$ and $10^{-1}.$ 
	In such cases, the structure reaches contact and the fluid flows abundantly into the porous layer, which prevents the release of contact. When the inlet pressure increases, the valve opens and the flow is restored.
	In all the other tests the interface is bouncing, but with a different reaction time, linked to the conductivity value. There is a slight difference in the first release time, but the more visible differences are on the second bounce. Both, the second contact instant and the final release, are considerably sensitive to the changing in $K_\tau=K_n$. 
	Finally, let notice that taking $K_\tau\rightarrow 0$ and $K_n \rightarrow 0$ we converge to the situation of no porous layer on $\Sigma_p$, as we can see in Figure~\ref{fig:dispxy}. 	
%	These results clearly show that both, Algorithm~\ref{alg:FSIDarcy} and the no-porous model, are able to capture the interface dynamics before and after contact with the upper wall. In addition, the elastic properties of contact are now clearly changing when the porous layer is introduced.
	
	\begin{figure}[ht]
		\centering
		\subfigure[y-velocity at $ t = 0.25 $.]{\includegraphics[width=0.45\linewidth]{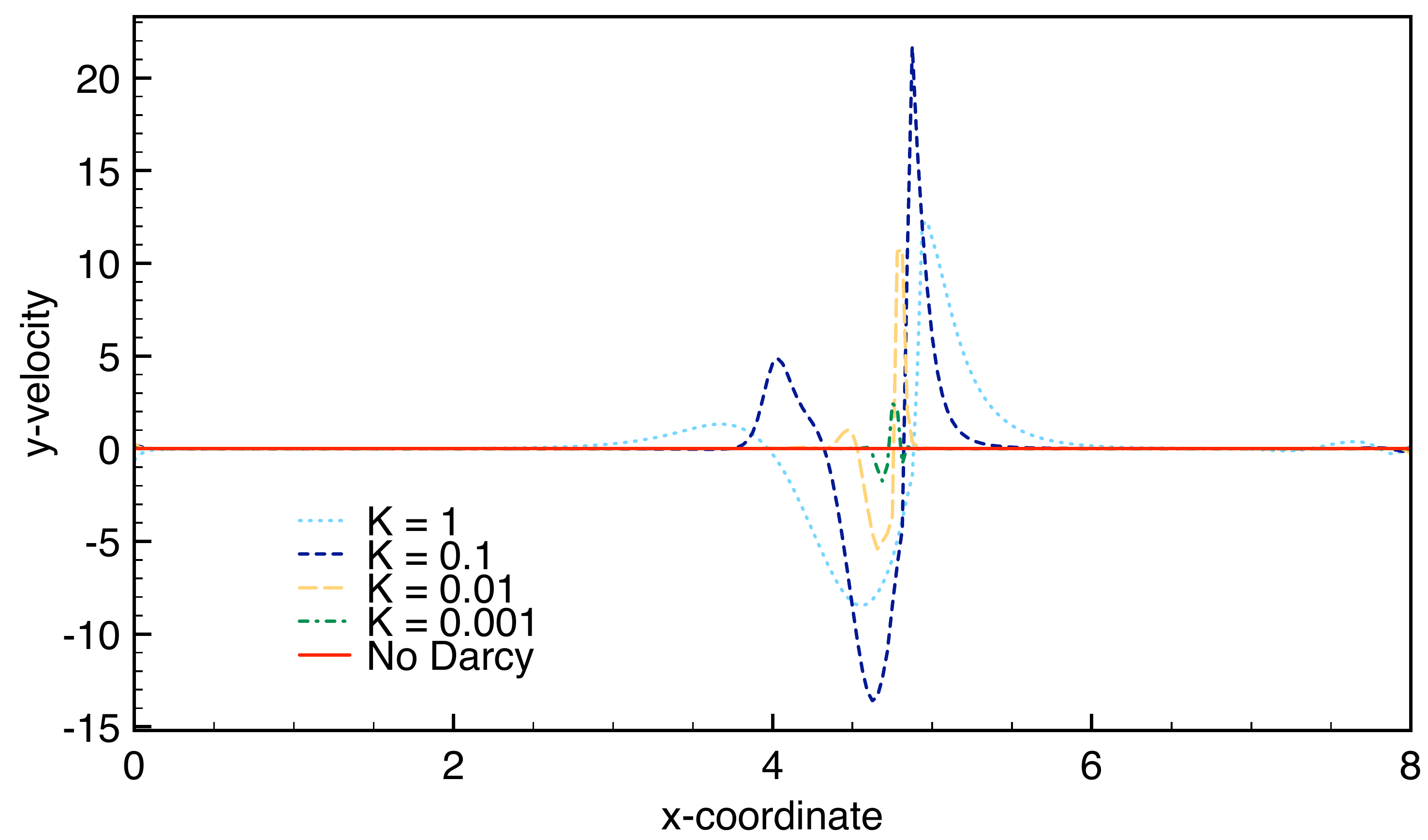}}
		\subfigure[y-velocity at $ t = 0.45 $.]{\includegraphics[width=0.45\linewidth]{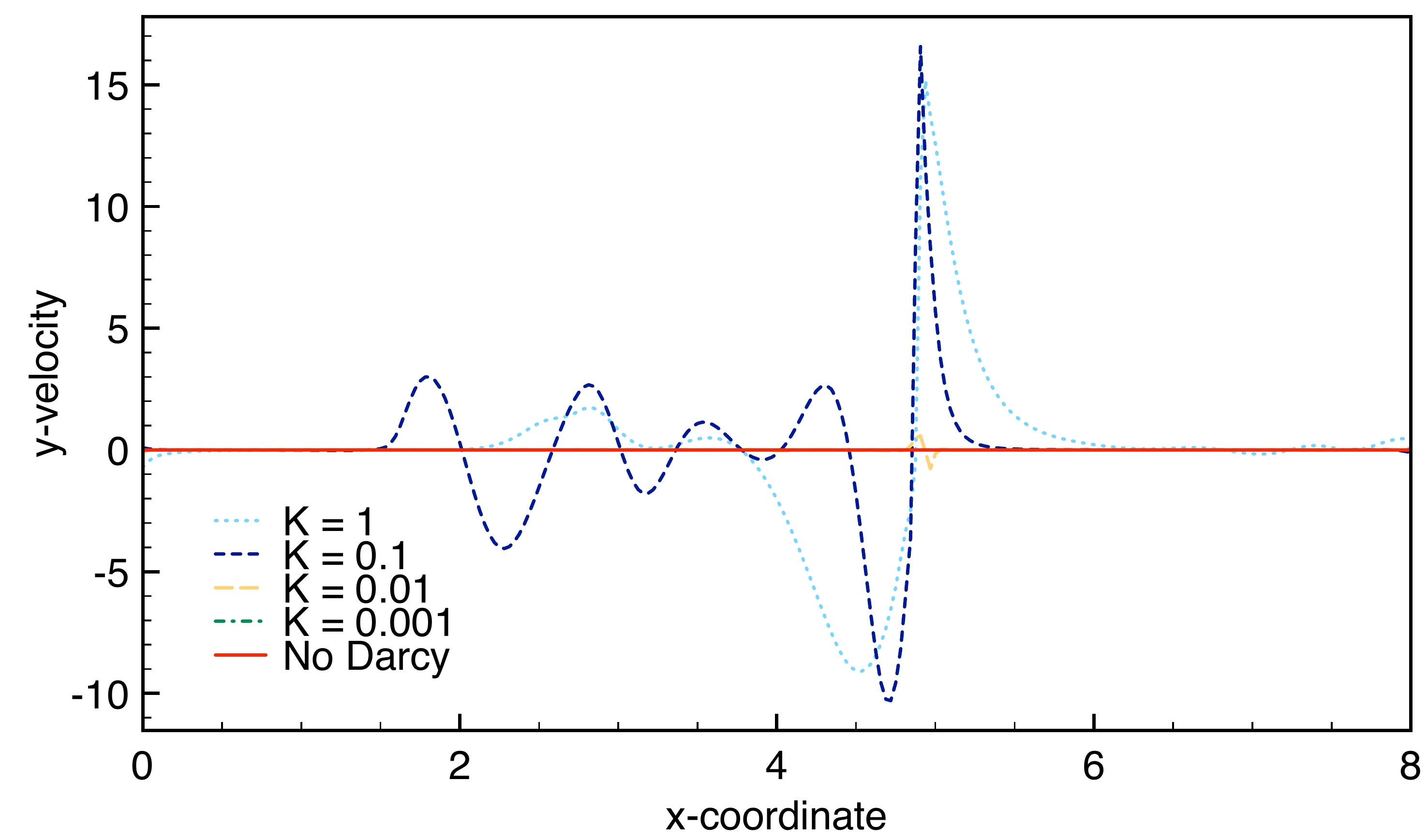}}
		\caption{Fluid velocity on $\Sigma_p.$}
		\label{fig:fluidYVelSigmaP}
	\end{figure}
	\begin{figure}[hbt]
	\centering
	\subfigure[$x$-displacement.]{\includegraphics[width=0.45\linewidth]{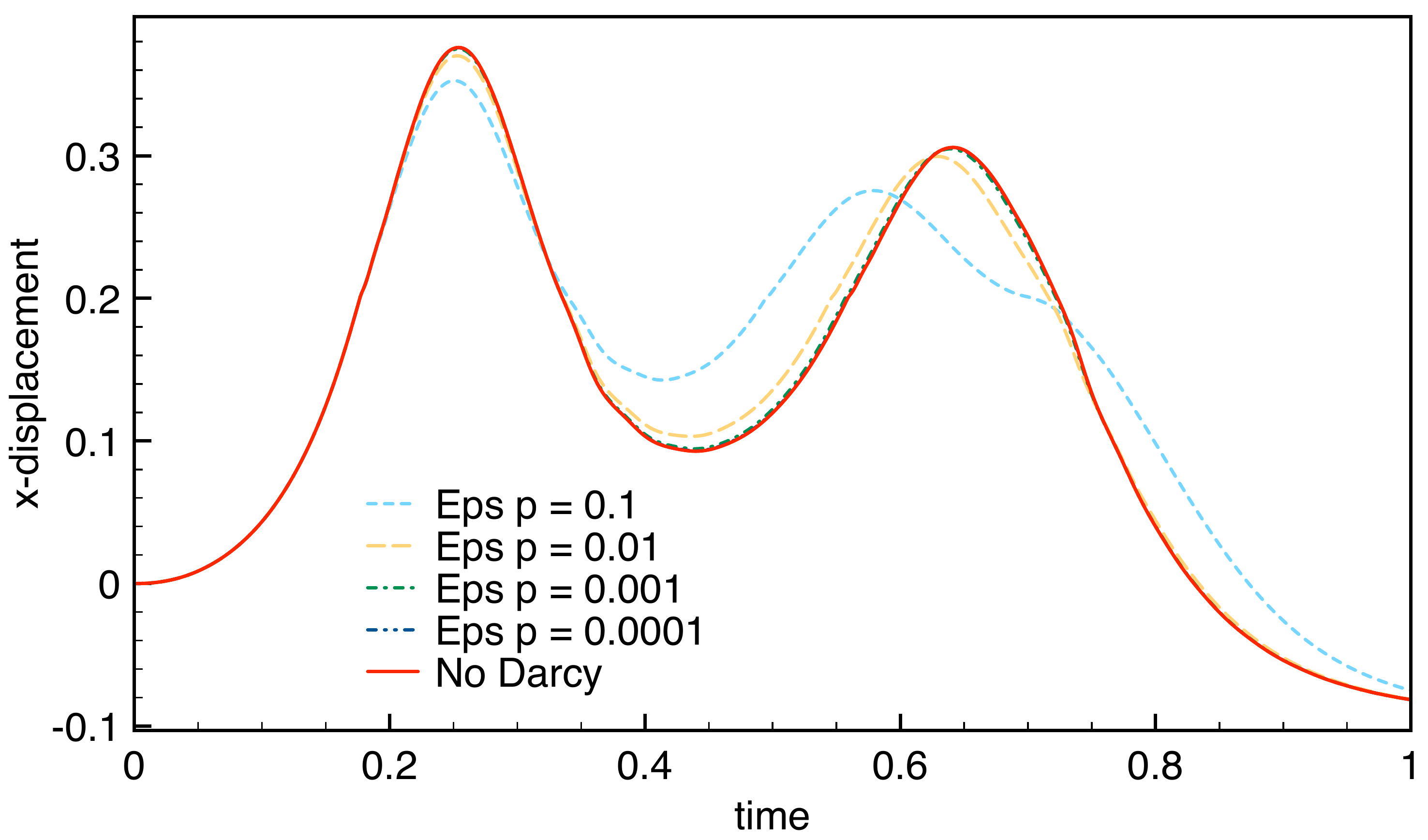}}
	\subfigure[$y$-displacement.]{\includegraphics[width=0.45\linewidth]{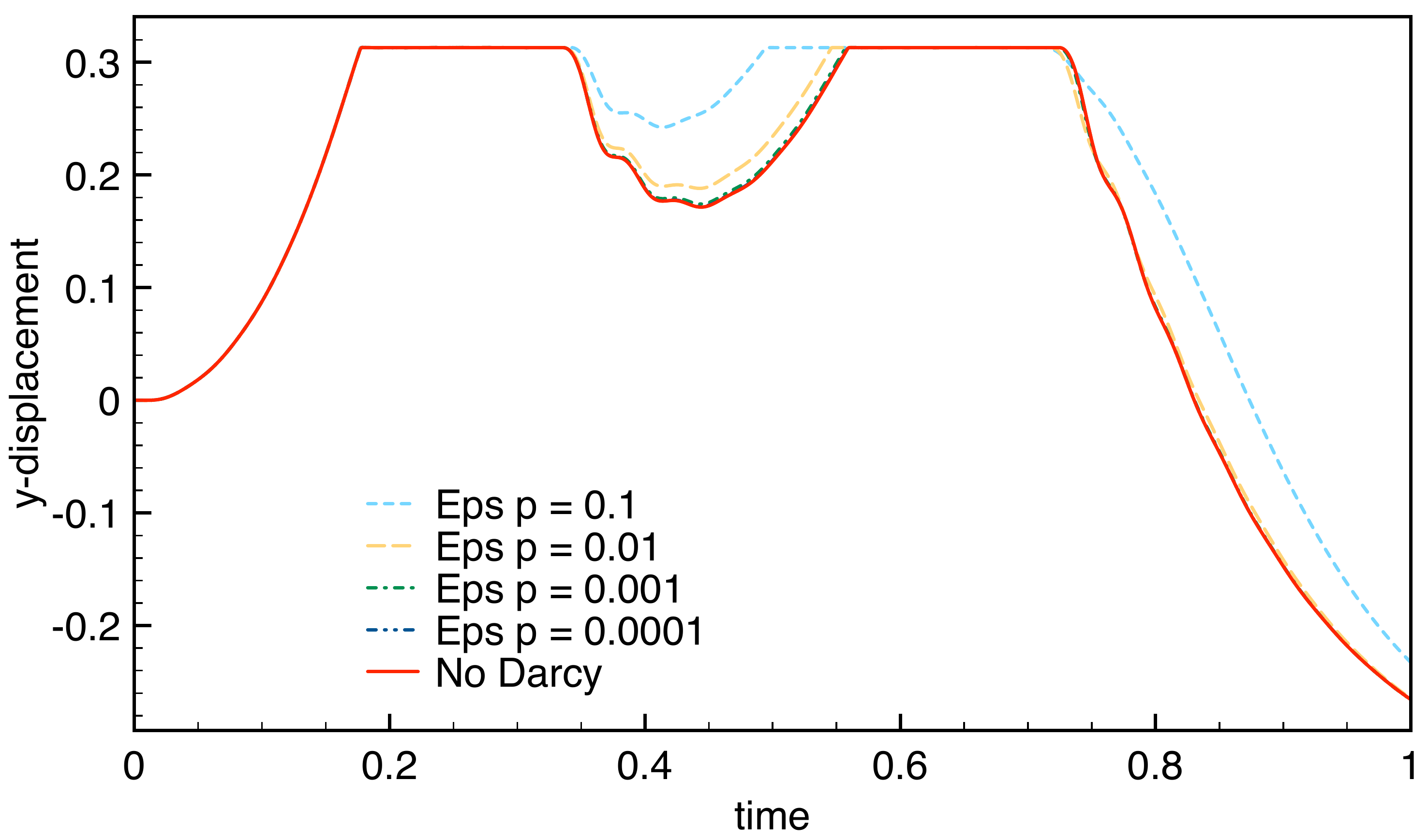}}
	\caption{Time evolution of the $x$-displacement (a) and $y$-displacement (b) for the structure endpoint $B$ for different values of $\varepsilon_p$}
	\label{fig:multiEps}
\end{figure}
\begin{figure}[hbt]
	\centering
	\subfigure[First release.]{\includegraphics[width=0.45\linewidth]{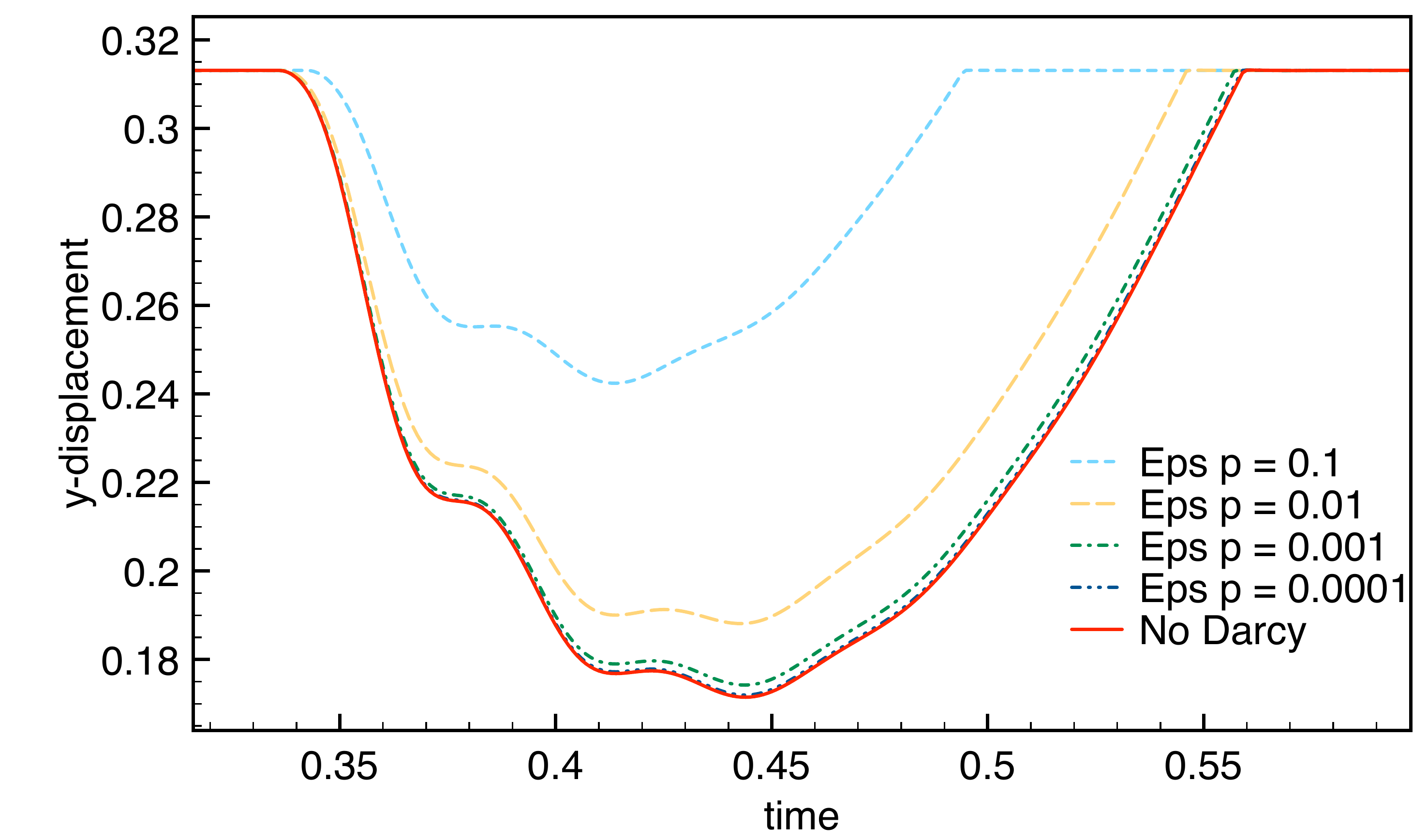}}
	\subfigure[Second release.]{\includegraphics[width=0.45\linewidth]{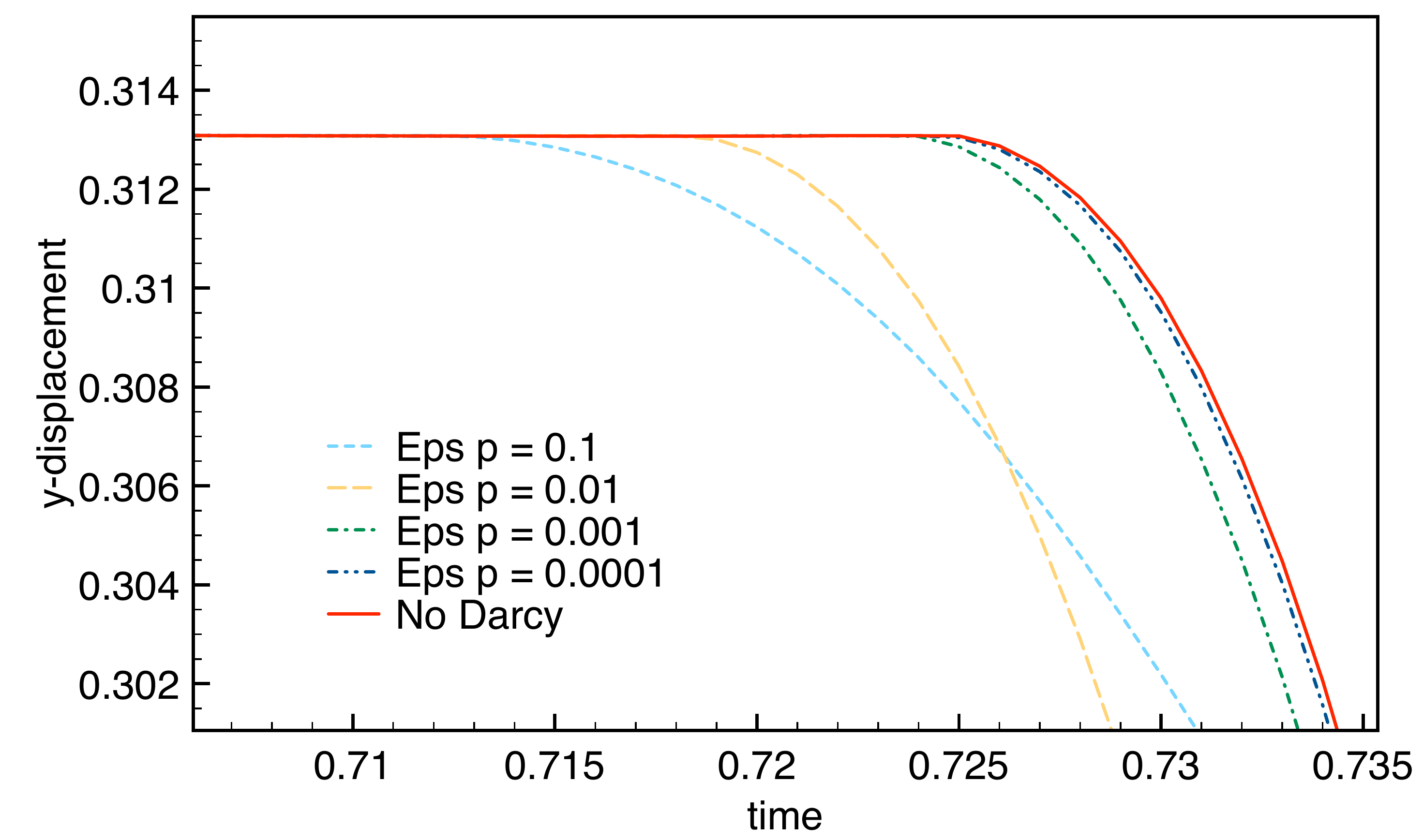}}
	\caption{Time evolution of the $y$-displacement for the structure endpoint $B$, between first release and second contact (a) and after second release (b).}
	\label{fig:multiEpsZoom}
\end{figure}

	Similar observations can be inferred from Figure~\ref{fig:interfacePos}, which shows the interface configuration at time $t=0.25$ (during contact), $t=0.45$ (after the first release) and $t=1$ (when the flow is restored). We can see that for $K_\tau=K_n=1$ and $10^{-1}$ the valve does not bounce, but it only releases once the inlet pressure increases (see Fig.~\ref{fig:interfacePos}(c)).
	Decreasing the conductivity of the porous medium increases the structure sliding at contact (see Fig.~\ref{fig:interfacePos}(a)) and the bouncing force applied on the structure (see Fig.~\ref{fig:interfacePos}(b)). 
	
	Figure~\ref{fig:fluidpressure} displays the fluid pressure (continuous line) and the porous pressures (dashed line) at time $t=0.25$ and $t=0.45$. As expected, both pressures remain close. At time $t=0.25$ the structure is in contact with the upper wall, therefore there is a high pressure gradient that decreases by increasing the conductivity.
%	, while at $t=0.45$ the cases with $K=10^{-2}$ and $10^{-3}$ approach the solution without Darcy getting a smooth pressure. 	
	
	Figure~\ref{fig:fluidXVelSigmaP} shows the fluid $x-$velocity along the porous layer $\Sigma_p$ at two different instants. As we can see, the horizontal velocity is not zero also during contact as effect of the porous layer. 
	As expected, the higher the conductivity the greater is the velocity magnitude and a larger area of the porous layer is leaking or pushing fluid inside the domain. 
	In Figure~\ref{fig:fluidYVelSigmaP} we report the fluid $y-$velocity on $\Sigma_p$. 
%	When contact is occurring the topological changes is happening (up to the contact relaxation parameter) and the porous serves as connection between the two fluid regions. It is visible that when we are in contact, the fluid is entering from the right side of the domain into the porous layer and getting inside the left fluid area through the porous layer. 
	The effect is more localised near the contact area except for cases of $K_\tau = K_n = 1, 0.1$, where the porous layer is still leaking and entering also far from the contact area.

%\begin{remark}
%	The contact relaxation parameter $\varepsilon_h$ must be chosen big enough in order to avoid penetration. No overlaps of the interface and $\Sigma_p$ are allowed in this method, contrarily to \cite{BurmanFernandezFrei}. 
%\end{remark}

%%%%%%%%%%%%%%%%%%%%%%%%%%%%%%%%%%%%%%%%%%%
%\subsubsection{Variation on $\varepsilon_p$}

We now explore the results when variations on the porous thickness $\varepsilon_p$ are considered. The porous hydraulic conductivity is  taken $K_\tau = K_n = 10^{-3}$. We explore results for $\varepsilon_p \in \{10^{-i}\}_{i=1}^4$. The outcome is shown in Figure~\ref{fig:multiEps}.
%We observe two bounces with different heights depending on $\varepsilon_p$. 
For $\varepsilon_p \rightarrow 0$ the curves converge towards the results of no porous layer on the top wall. 
%In fact, taking $\varepsilon_p \rightarrow 0$ in Algorithm~\ref{alg:FSIDarcy}, we force the (porous layer equation), implies $\bu \cdot \bn_f \rightarrow 0$. 
%Hence, we retrieve the boundary condition considered for the valve without the porous model. 
%No visible difference between the solution without Darcy and with $10^{-4}$ are visible. This could be expected, considering that smaller porous layer allows less fluid to diffuse through the layer.

% at 0.805-0.483-0.01 = 0.321 the contact starts to be activated
No particular differences are visible at first contact between the structure and the upper wall.
During contact, the horizontal velocity is lower for higher values of $\varepsilon_p$, therefore, the bouncing force is also lower. 
In addition, the higher is $\varepsilon_p$, the later is the first release, the lower is the rebound force and, consequently, the earlier is the second contact and release.
For illustration purposes, we report in Figure~\ref{fig:multiEpsZoom} a zoom of the $y$-displacement between the first release and the second contact instants.

%\FloatBarrier
%%%%%%%%%%%%%%%%%%%%%%%%%%%%%%%%%%%%%%%%%%%
	\begin{figure}[h!]
	\centering
	\subfigure[$x$-displacement.]{\includegraphics[width=0.48\linewidth]{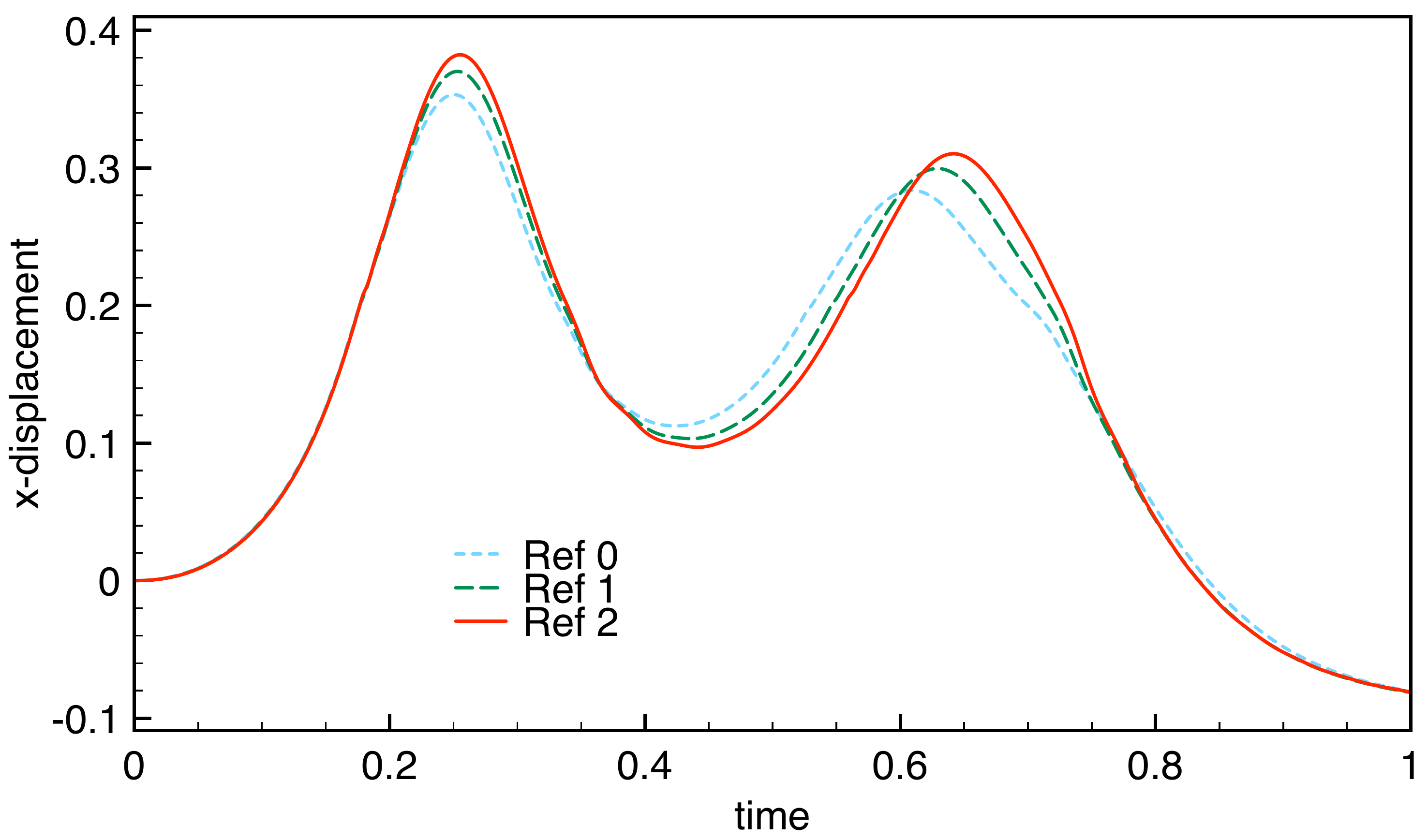}}
	\subfigure[$y$-displacement.]{\includegraphics[width=0.48\linewidth]{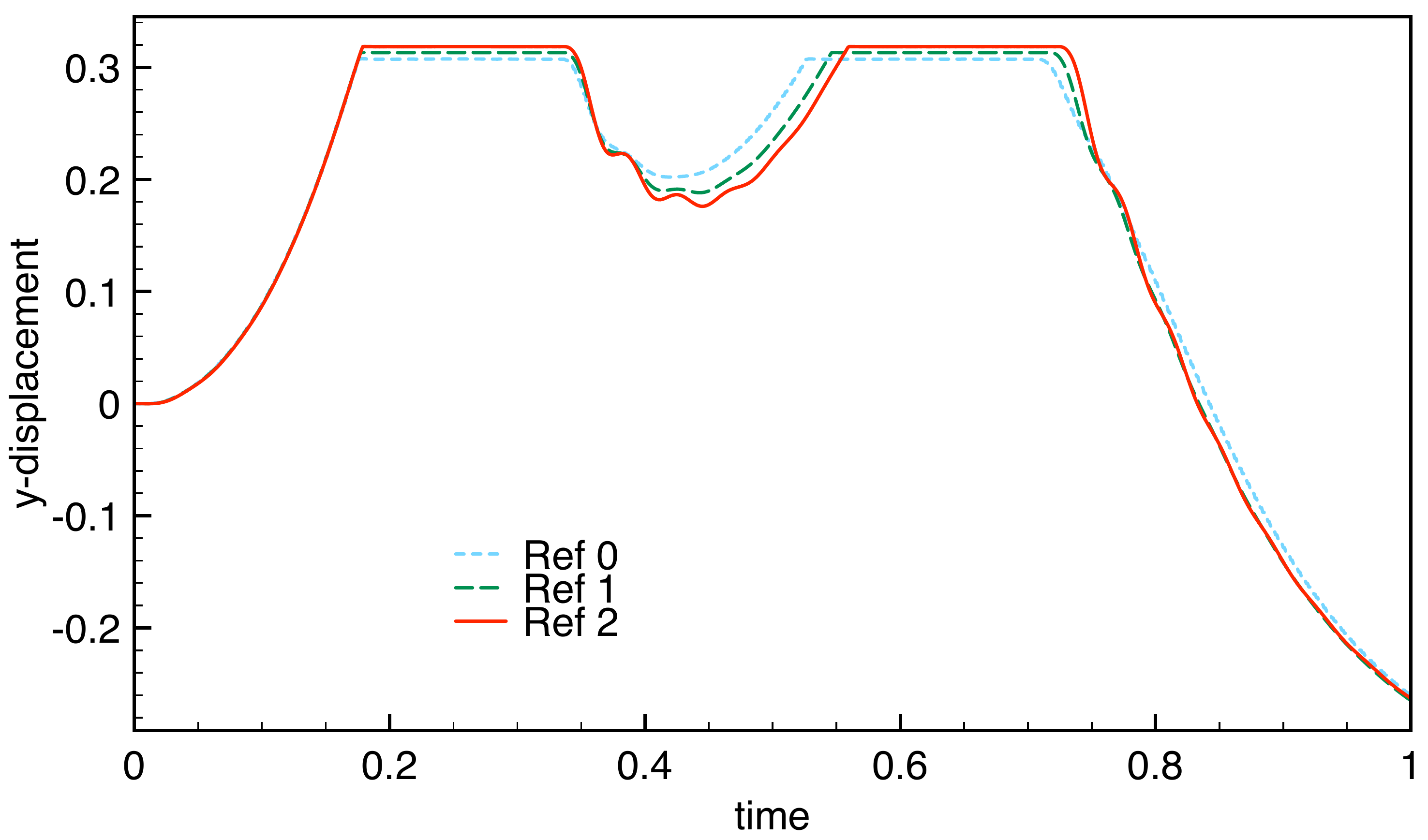}}
	\caption{Time evolution of the $x$-displacement (a) and $y$-displacement (b) for the structure endpoint $B$, with different levels of refinement. }
	\label{fig:Refinement}
\end{figure}

\paragraph{Refinement in space and time}
We explore the convergence behavior taking three levels of space and time refinement, namely $\big(\delta t,h\big) \in \big\{  2 \cdot 10^{-3} \cdot 2^{-i} \, , \,0.07 \cdot 2^{-i} \big\}_{i=0}^2 $. 
The coarser fluid and solid meshes are made of 5\,120 triangles and 26 segments, respectively. 
The second meshes consists of 20\,480 triangles and 50 edges, while the finest one has 81\,920 triangles and 102 segments. 
The porous conductivity is chosen $K_\tau = K_n = 10^{-3}$, and the contact relaxation parameter $\varepsilon_h = \varepsilon_h(h)$, chosen $\varepsilon_h \in \{ 0.02 \cdot 2^{-i} \}_{i=0}^2 $.

We show in Figure~\ref{fig:Refinement} the results with these three refinement levels. 
We observe that the bouncing height is lower for the coarser mesh and that 
%, refining, we are able to catch better the contact dynamics. In particular the small oscillations between the two contact instants.
the intermediate level of refinement provides a reasonable approximation. 
We can also observe that, due to different contact relaxation parameters, contact and release occur at different instants and heights. % $y$-coordinates.

%\clearpage
\subsection{Falling and bouncing elastic ball}
\label{subsec:num.thick}

\begin{figure}
\centering
\resizebox*{0.4\textwidth}{!}{
\begin{picture}(0,0)%
\includegraphics{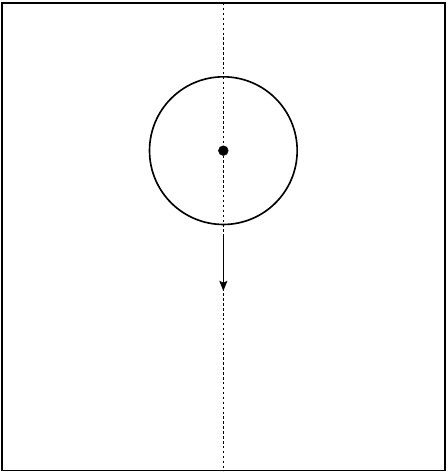}%
\end{picture}%
\setlength{\unitlength}{1036sp}%
\begingroup\makeatletter\ifx\SetFigFont\undefined%
\gdef\SetFigFont#1#2{%
  \fontsize{#1}{#2pt}%
  \selectfont}%
\fi\endgroup%
\begin{picture}(8166,8616)(2218,-9094)
\put(7426,-2221){\makebox(0,0)[lb]{\smash{{\SetFigFont{5}{6.0}{\color[rgb]{0,0,0}$\Sigma$}%
}}}}
\put(6481,-5416){\makebox(0,0)[lb]{\smash{{\SetFigFont{5}{6.0}{\color[rgb]{0,0,0}$f_s$}%
}}}}
\put(6436,-7756){\makebox(0,0)[lb]{\smash{{\SetFigFont{5}{6.0}{\color[rgb]{0,0,0}$\Gamma^{\rm sym}$}%
}}}}
\put(6661,-3796){\makebox(0,0)[lb]{\smash{{\SetFigFont{5}{6.0}{\color[rgb]{0,0,0}$\Omega^s$}%
}}}}
\put(7876,-6226){\makebox(0,0)[lb]{\smash{{\SetFigFont{5}{6.0}{\color[rgb]{0,0,0}$\Omega^f$}%
}}}}
\put(7876,-8926){\makebox(0,0)[lb]{\smash{{\SetFigFont{5}{6.0}{\color[rgb]{0,0,0}$\Sigma_p$}%
}}}}
\end{picture}%
}
\caption{\label{fig.ballgeo} Configuration of the numerical example with the falling ball. The computational domain is reduced to the right half due to symmetry on the axis $\Gamma^{\rm sym}$.}
\end{figure}
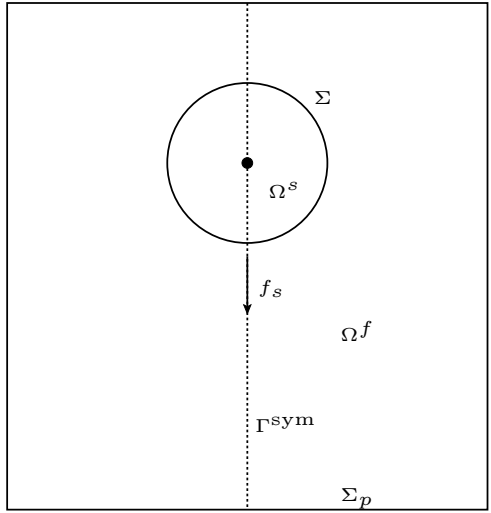

In this section, we consider the example of a falling and bouncing elastic ball in a cylinder, which is filled with a water-glycerin mixture. 
%which is motivated by experiments provided in~%\cite{RichterHagemeier}. 
As we are interested in a detailed numerical study, we restrict ourselves to the two-dimensional case here and consider a box of size $6$cm$\times 7.5$cm. The ball has a radius of 1cm and is kept initially at rest at a distance of 4cm from the bottom. The ball falls down due to gravity $f_s=-9.81\frac{m}{s^2}$ and bounces back after the impact. Due to symmetry, we can reduce the computational domain to the right half by imposing symmetry boundary conditions on the midplane $\Gamma^{\rm sym}$, see Figure~\ref{fig.ballgeo}.

\begin{figure}[t]
\includegraphics[width=0.33\textwidth]{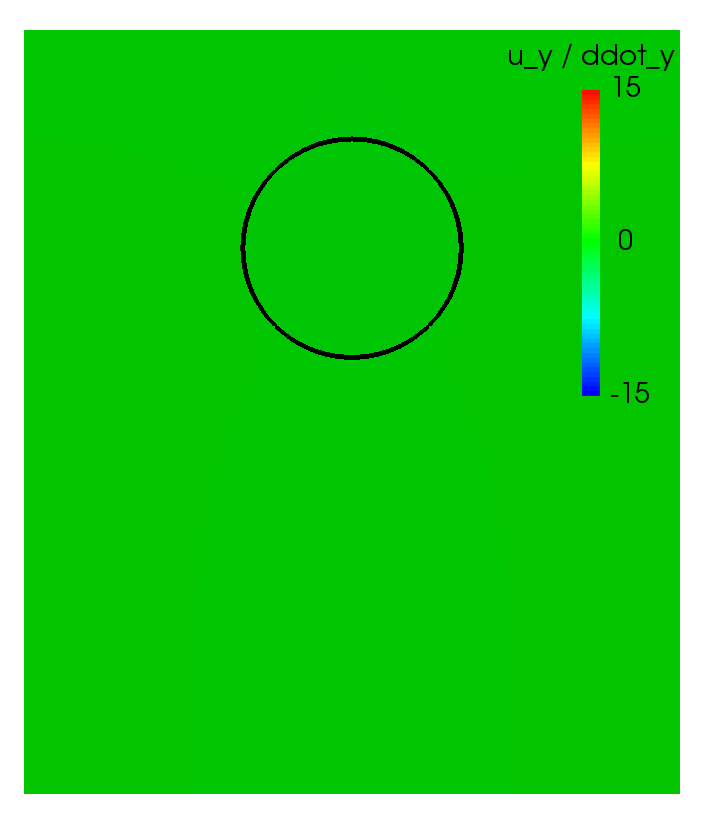}\hfil 
\includegraphics[width=0.33\textwidth]{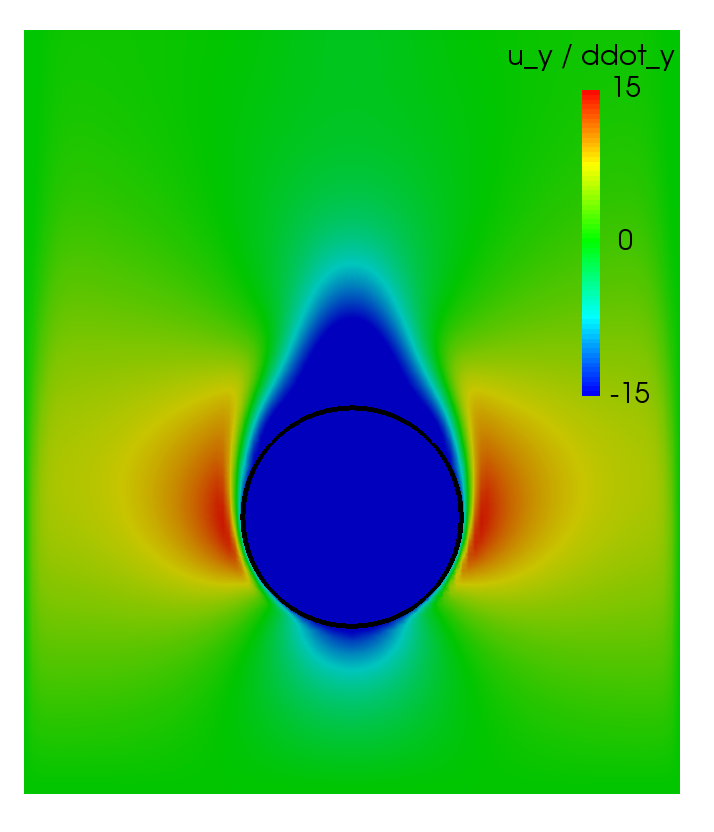}\hfil 
\includegraphics[width=0.33\textwidth]{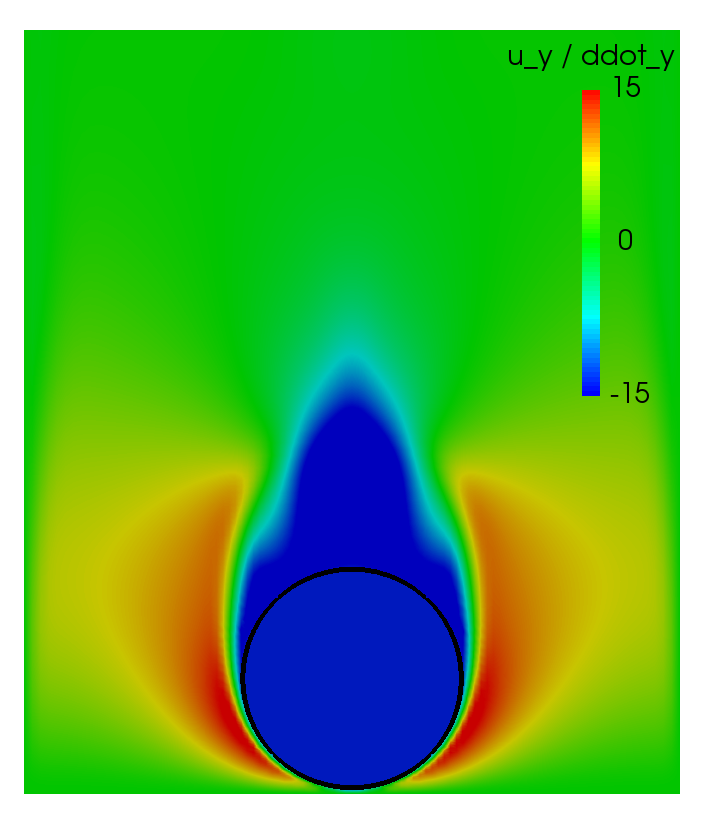}\\ 
\includegraphics[width=0.33\textwidth]{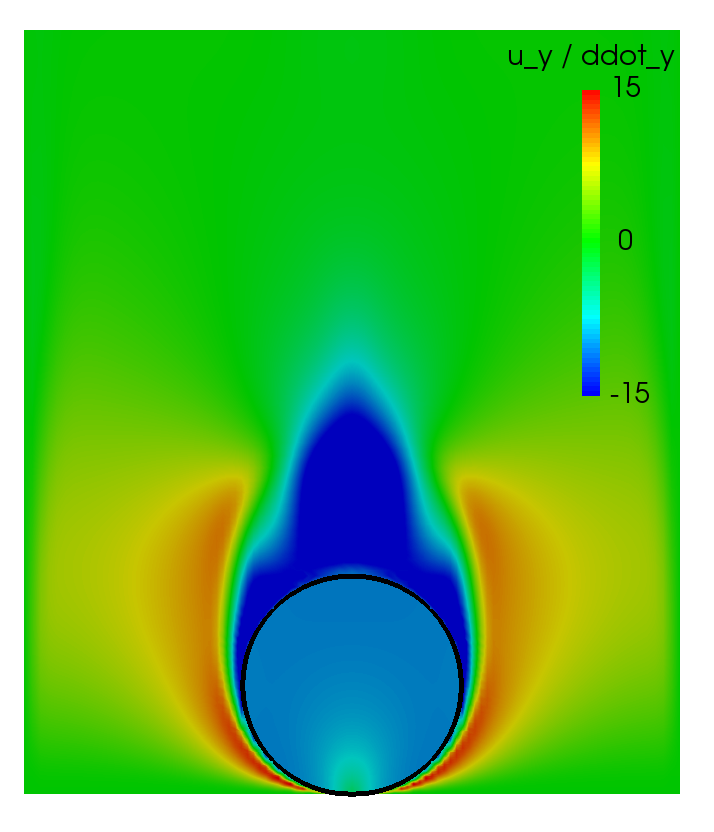}\hfil 
\includegraphics[width=0.33\textwidth]{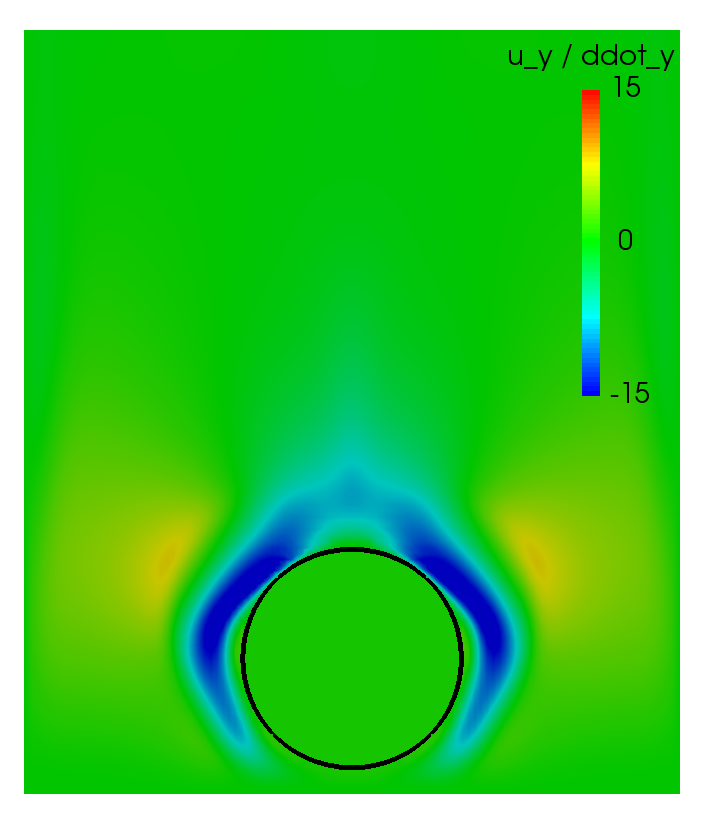}\hfil
\includegraphics[width=0.33\textwidth]{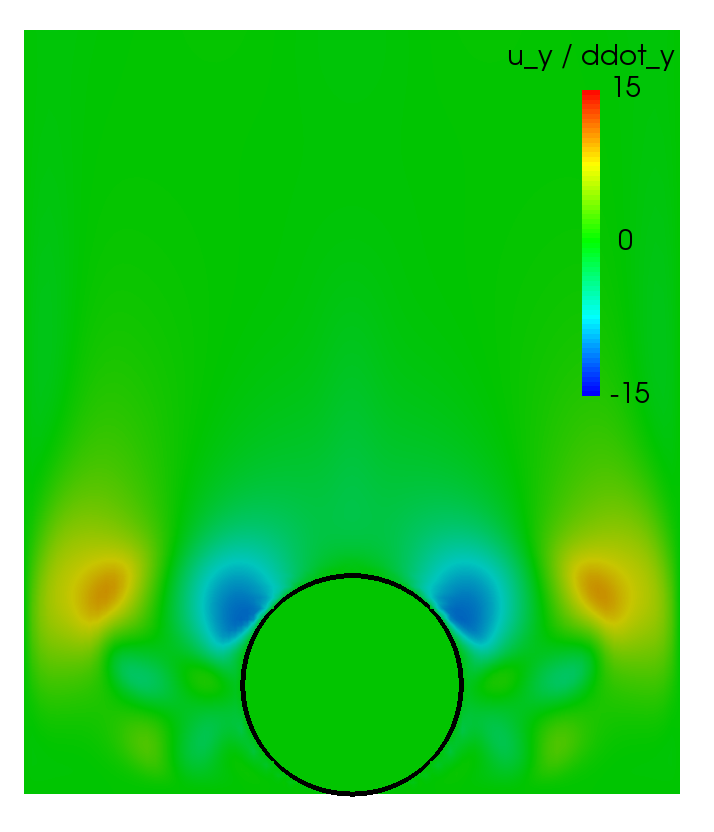}\caption{Illustration of the vertical velocities in the falling ball example at $t=0, 0.28, 0.378$ (top, left to right) and $t=0.382, 0.47$ and $0.63$ (bottom).\label{fig.fall}}
\end{figure}

We use a Fully Eulerian approach to solve the coupled problem, as described in Section~\ref{sec.Eulerian}. For simplicity, we consider here a linear elastic material, %i.e.
%\begin{align*}
%a^s(d_h, w_h) = (\sigma_s(d_h), \nabla w_h)_{\Omega^{\rm s}(t)},
%\end{align*}
where the Cauchy stress tensor $\sstress$ is given by
\begin{align*}
 \sigma_s(d) = 2\mu_s E(d) + \lambda_s \text{tr} (E(d)) I, 
 \quad E(d)=\frac{1}{2} \left( \nabla d + \nabla d^T\right),
 \end{align*}
with Lam\'e parameters $\lambda_s = 7.64\cdot 10^6 \frac{kg}{ms^2}$ and $\nu_s = 1.04\cdot 10^6 \frac{kg}{ms^2}$.

For time discretisation, we use a variant of the backward Euler method, namely a modified dG(0) time-stepping, see~\cite{FreiRichter2017, FreiRichter2017Sammelband}.
We start with a time-step size of $\delta t=2\cdot 10^{-3}$s, which is reduced in a stepwise procedure up to the impact, where a small time-step size of $\delta t=1.25\cdot 10^{-4}$s is reached.

The (kinematic) viscosity of the water-glycerol mixture is $\mu_{\rm f} = 7\cdot 10^{-6} \frac{m}{s^2}$, the density $\rho_f= 1 141\frac{kg}{m^3}$ and the solid density $\rho_s= 1 351\frac{kg}{m^3}$.
Unless stated differently the parameters in the porous medium are chosen as $\epsilon_{\rm p}=10^{-4}$ and $K=K_n = K_\tau= 10^{-2}$ and the
numerical contact parameters are $\gamma_c=30\lambda_s$ and $\epsilon_g=\frac{h}{4}$.
For the Navier-Stokes-Darcy coupling, we use the Beavers-Joseph-Saffman condition in \eqref{eq:darcy} with $\alpha=1$.
  All the results have been obtained with the finite element library Gascoigne3d~\cite{Gascoigne}. We use a structured coarse grid ${\cal T}_{2h}$, which is highly refined close to $\Sigma_p$ with 3 201 vertices in total (unless specified differently).
In Figure~\ref{fig.fall}, we illustrate the vertical velocities $u_y$ and $\dot{d}_y$ of the falling ball at 6 instances of time.

\begin{figure}
\begin{minipage}{0.5\textwidth}
\includegraphics[width=1.1\textwidth]{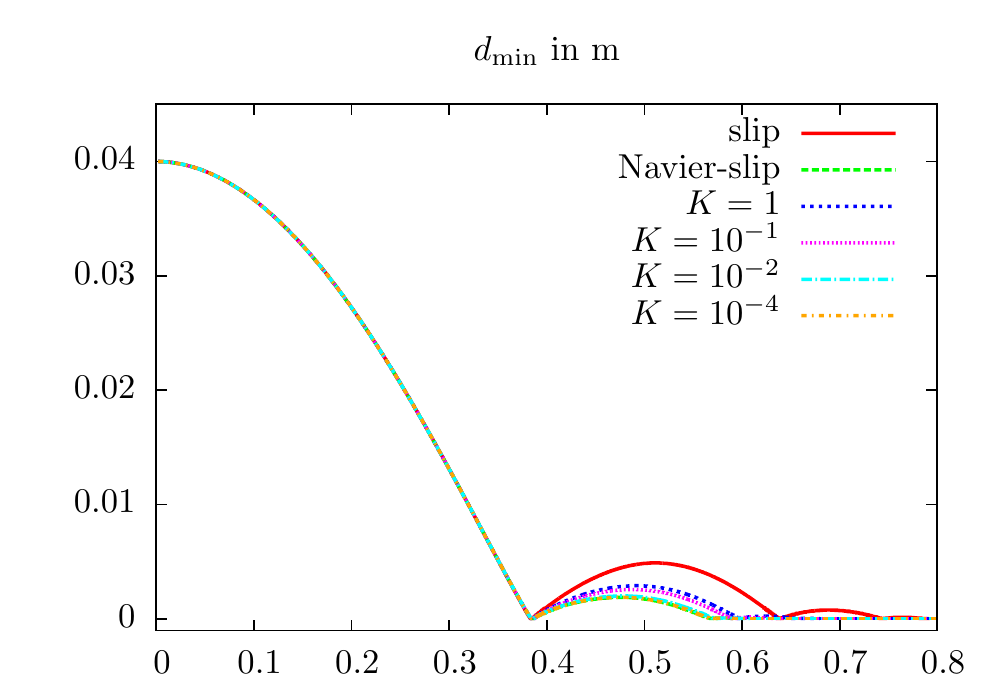}
\end{minipage}
\hfil
\begin{minipage}{0.5\textwidth}
\includegraphics[width=1.1\textwidth]{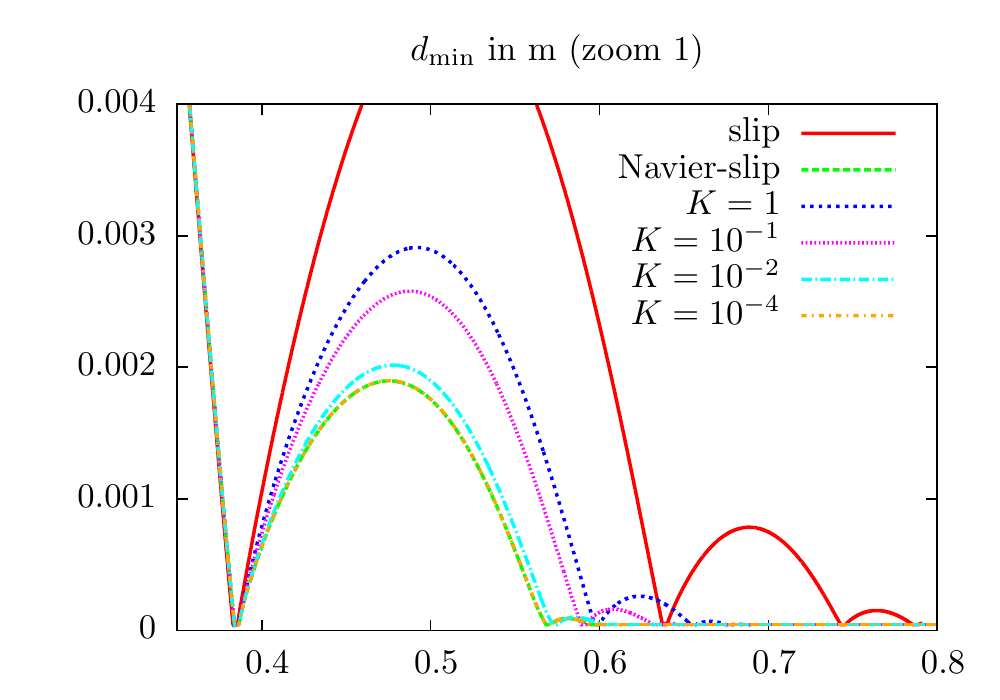}
\end{minipage}
\begin{minipage}{0.5\textwidth}
\includegraphics[width=1.1\textwidth]{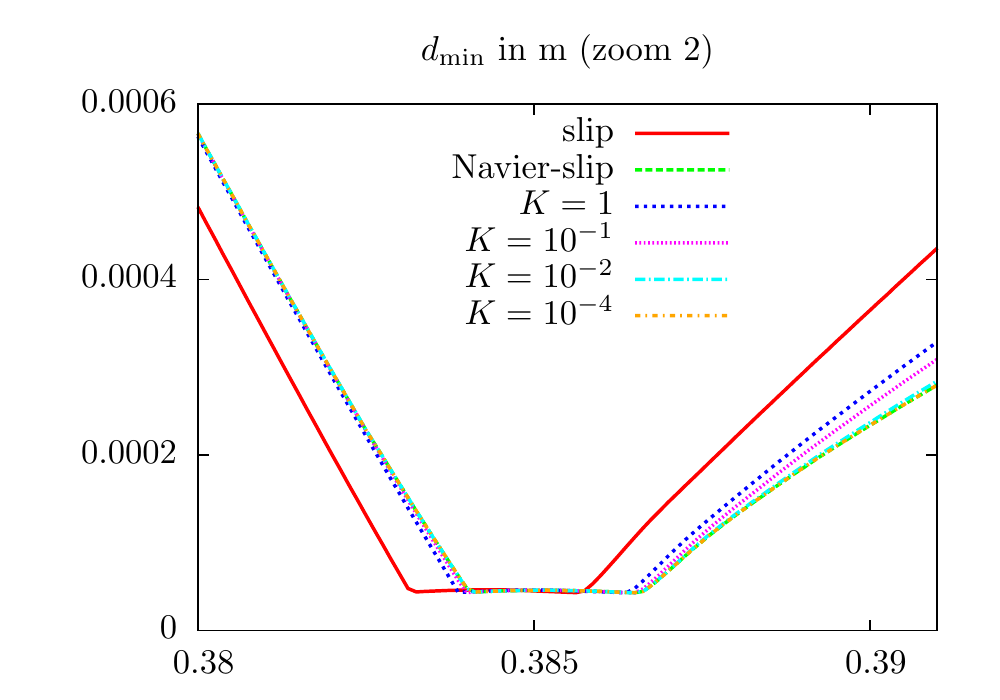}
\end{minipage}
\hfil
\begin{minipage}{0.5\textwidth}
\includegraphics[width=1.1\textwidth]{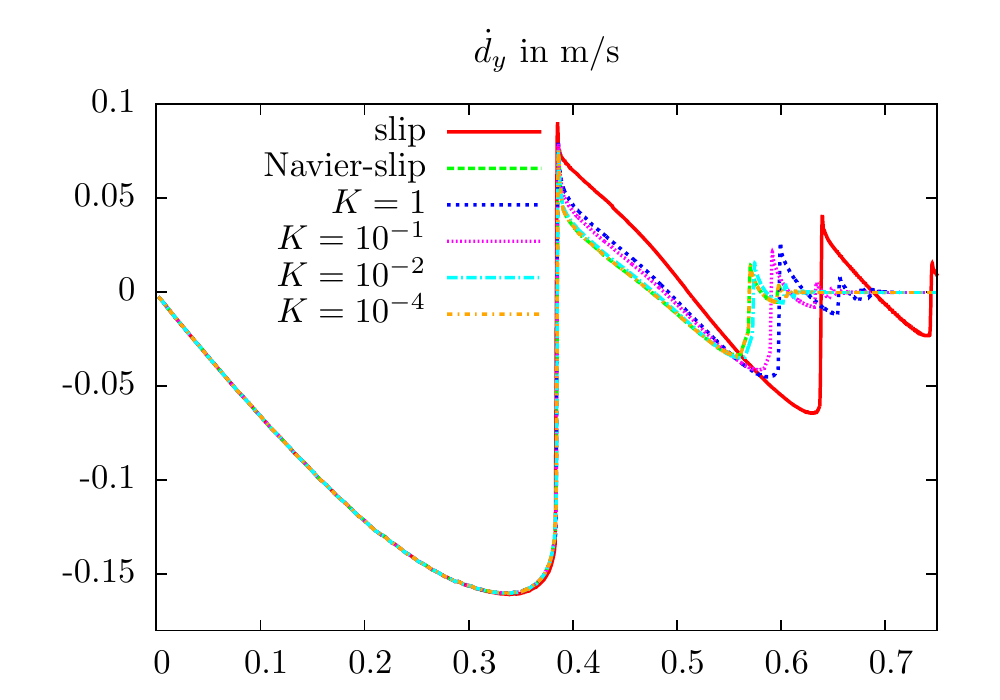}
\end{minipage}
\caption{Minimal distance $d_{\rm min}$ to the ground over time for different conductivities $K$ and compared to pure slip and no-slip conditions on the lower boundary. After the curve over the full time interval $[0,0.8s]$ on the top left, two different zooms are provided to illustrate the bounces (top right) and the time of impact (bottom left). \textit{Bottom right}: Space-averaged vertical velocity of the elastic ball over time.\label{fig.mindist}}
\end{figure}

\paragraph{Variation of the conductivity $K$}
In Figure~\ref{fig.mindist} we compare the minimum distance to the ground during the fall and before and after the impact for different conductivities $K$ with results obtained without a porous model ("No Darcy"), using either a full slip ($\sigma_{f,\tau}=0$) or a Navier-slip condition $\sigma_{f,\tau} = \frac{\alpha}{\sqrt{K_\tau \epsilon_{\rm p}}}$ with $\alpha=1$ on $\Sigma_p$. Note that in the latter case this is exactly the same tangential condition which is imposed by the Beavers-Joseph-Saffman coupling. There, the normal velocity $\bu \cdot \bn$ is however not necessarily zero, as the flow can enter into the porous medium.

 Depending on $K$, the ball bounces 4 or 5 times within the time interval $[0,0.8s]$ with different bouncing heights (top right). The last bounces are barely visible in the graphs shown here due to a very small bouncing height. Moreover, we observe for $K\to 0$ that the results converge towards the results obtained with a Navier-slip-boundary conditions ("No Darcy"), the curve for $K=10^{-4}$ showing no visible differences. For larger $K$ the curves get slightly closer towards the results for a pure slip boundary condition.

Concerning the time of impact (bottom left of Fig.~\ref{fig.mindist}), we observe that the impact happens slightly later, the smaller the conductivity $K$. The latest impact is observed for $K=10^{-4}$ and the pure Navier-slip condition ("No Darcy"). This dependence on $K$ is expected, as the resistive fluid forces that act against the contact are higher for smaller conductivities. The earliest impact is observed for the pure slip condition, followed by $K=1$. Moreover, we observe that a small distance of about $4.3\cdot 10^{-5}m$, which lies slightly below the imposed gap distance of $\epsilon_g=5\cdot 10^{-5}m$, remains in all cases.

From the upper right picture we can infer the bouncing height depending on $K$. It holds that the earlier the impact, the higher the impact velocity and hence, we observe a larger bouncing height. Consequently, we observe the largest bounce for pure slip-conditions, followed by $K=1$, and the smallest one for pure no-slip conditions and $K=10^{-4}$. 
%For $K\geq 10$, we observe again smaller bounces. This is due to the fact that the (porous and fluid) pressure after the impact are smaller, as fluid can enter through the porous part easily and hence, the resistance to the second fall is lower. 

On the bottom right of Figure~\ref{fig.mindist} we show the corresponding vertical velocity $\dot{d}_y$ within the elastic ball, averaged in space. We see that the absolute value of the upwards velocity after the first impact is by more than 30\% smaller than the absolue value of the impact velocity in all cases, which shows the dissipative impact of the fluid. 
%Moreover, we observe that the upwards velocity for the high conductivities ($K\geq 10$) is relatively high right after the impact, but diminishes faster than for smaller $K$ afterwards, which further substantiates the explanation given above.

%\begin{figure}
%\centering
%\begin{minipage}{0.65\textwidth}
%\includegraphics[width=\textwidth]{velocity.pdf}
%\end{minipage}
%\caption{Space-averaged vertical velocity of the elastic ball over time.\label{fig.uy}}
%\end{figure}

%\begin{remark}
%	The contact relaxation parameter $\epsilon_h$ must be chosen big enough in order to avoid penetration. No overlaps of the interface and $\Sigma_p$ are allowed in this method, contrarily to \cite{BurmanFernandezFrei}. 
%\end{remark}

\begin{figure}[bt]
%\begin{minipage}{0.5\textwidth}
%\includegraphics[width=1.1\textwidth]{mindist_epsp.pdf}
%\end{minipage}
%\hfil
\begin{minipage}{0.5\textwidth}
\includegraphics[width=1.1\textwidth]{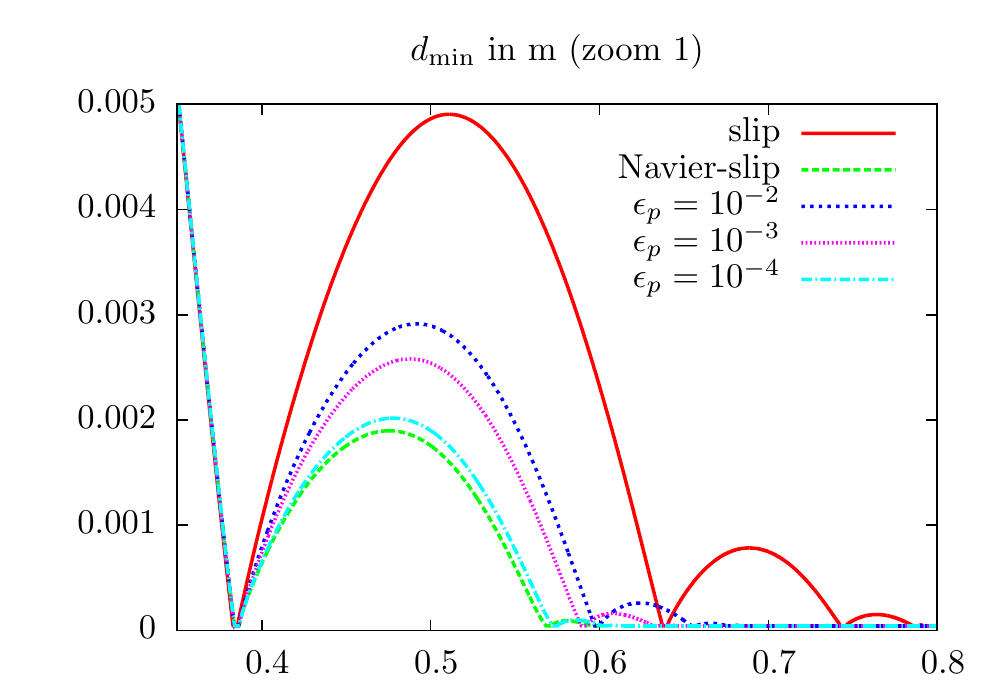}
\end{minipage}
\begin{minipage}{0.5\textwidth}
\includegraphics[width=1.1\textwidth]{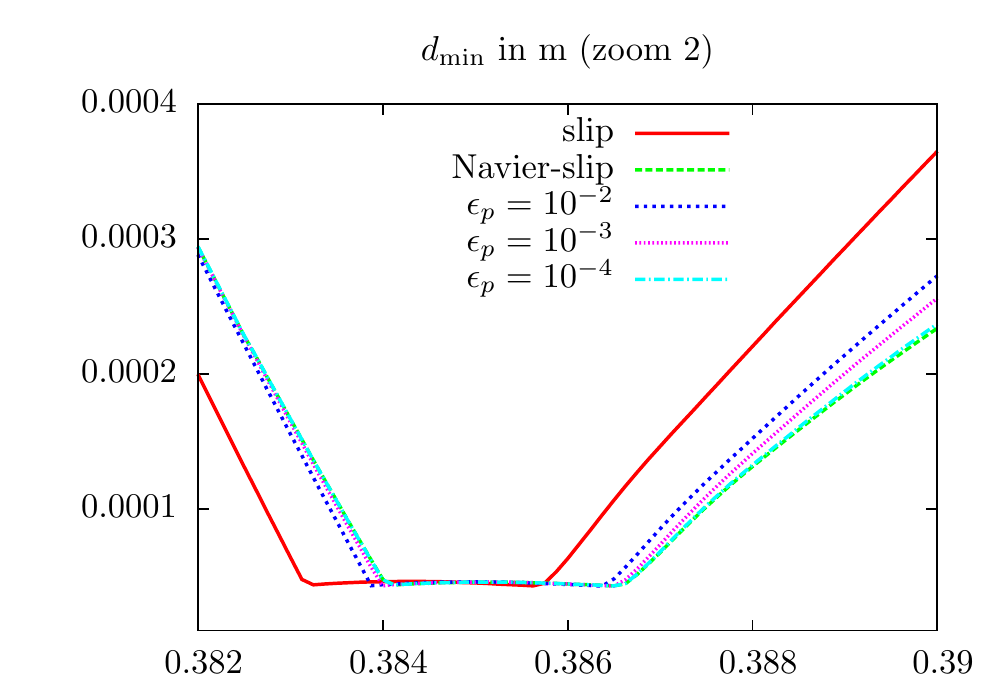}
\end{minipage}
%\hfil
%\begin{minipage}{0.5\textwidth}
%\includegraphics[width=1.1\textwidth]{mindist_zoom_epsp.pdf}
%\end{minipage}
\caption{Minimal distance $d_{\rm min}$ to the ground over time for different thicknesses $\epsilon_{\rm p}$ of the porous layer and compared to pure slip and no-slip conditions on the lower boundary.\label{fig.mindist_epsp}}
\end{figure}

\paragraph{Variation of $\epsilon_{\rm p}$}
In Figure~\ref{fig.mindist_epsp} we vary the thickness $\epsilon_{\rm p}$ of the porous layer. We obtain 3 to 5 bounces with different heights depending on $\epsilon_{\rm p}$. For $\epsilon_{\rm p}\to 0$, the curves converge towards the results for a pure Navier-slip condition ("No Darcy"), as the first equation in \eqref{eq:darcy} implies $\bu \cdot \bn \to 0$. This is also what one expects from the physical model, as a smaller porous layer allows less fluid to diffuse through the layer.
In the curve on the right of Figure~\ref{fig.mindist_epsp}, we see that the contact happens earlier the larger $\epsilon_{\rm p}$ is. This can again be explained by the smaller resistance of the fluid "against" the contact, when this is allowed to escape through the porous layer. The larger impact velocity for larger $\epsilon_{\rm p}$ has again the effect that the bounce is higher for larger $\epsilon_{\rm p}$. 

\paragraph{Variation of the contact parameter $\gamma_c$}

\begin{figure}
%\begin{minipage}{0.5\textwidth}
%\includegraphics[width=1.1\textwidth]{mindist_epsp.pdf}
%\end{minipage}
%\hfil
\begin{minipage}{0.5\textwidth}
\includegraphics[width=1.1\textwidth]{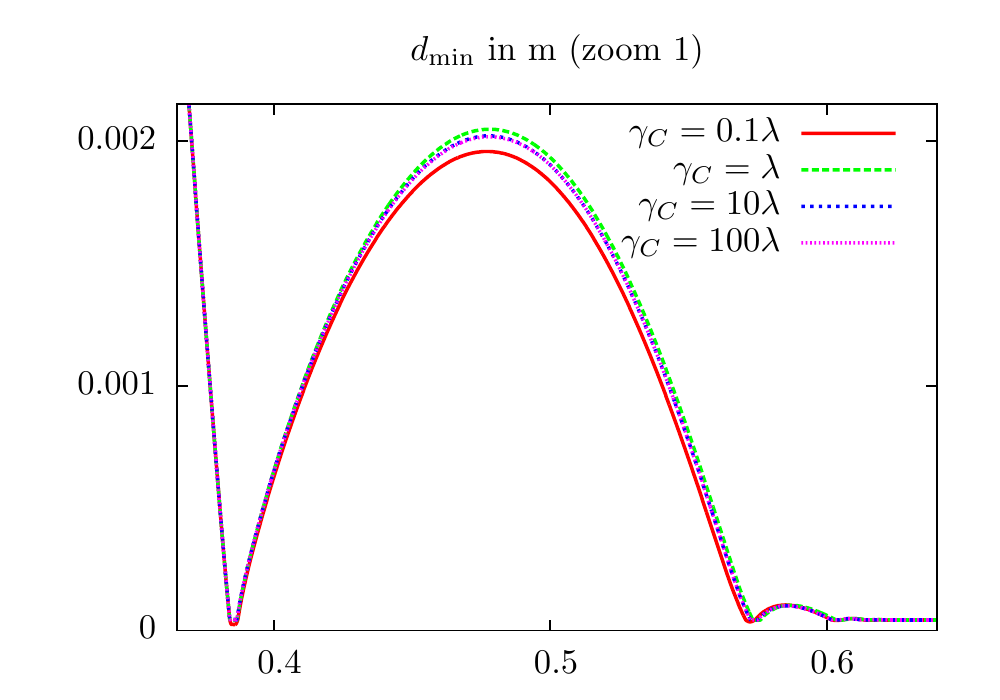}
\end{minipage}\hfil
\begin{minipage}{0.5\textwidth}
\includegraphics[width=1.1\textwidth]{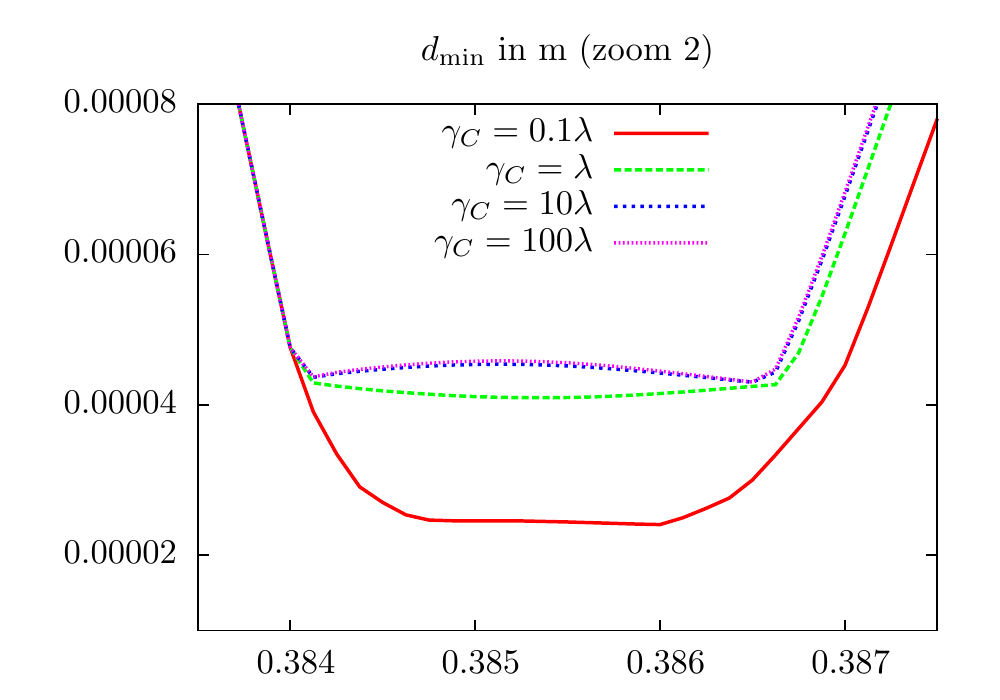}
\end{minipage}\\
\begin{minipage}{0.5\textwidth}
\includegraphics[width=1.1\textwidth]{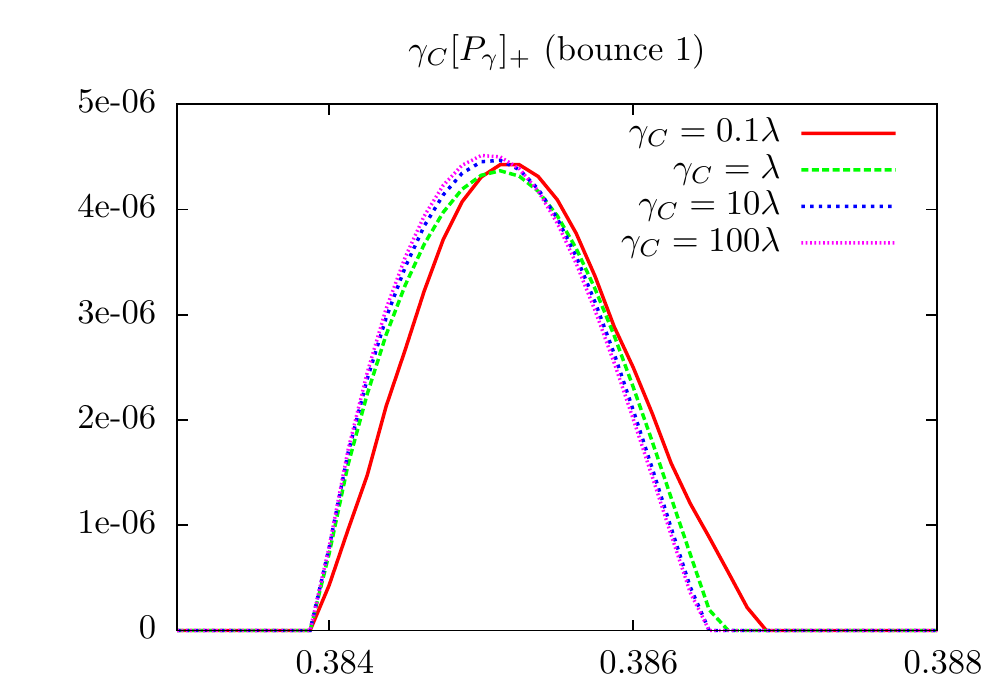}
\end{minipage}\hfil
\begin{minipage}{0.5\textwidth}
\includegraphics[width=1.1\textwidth]{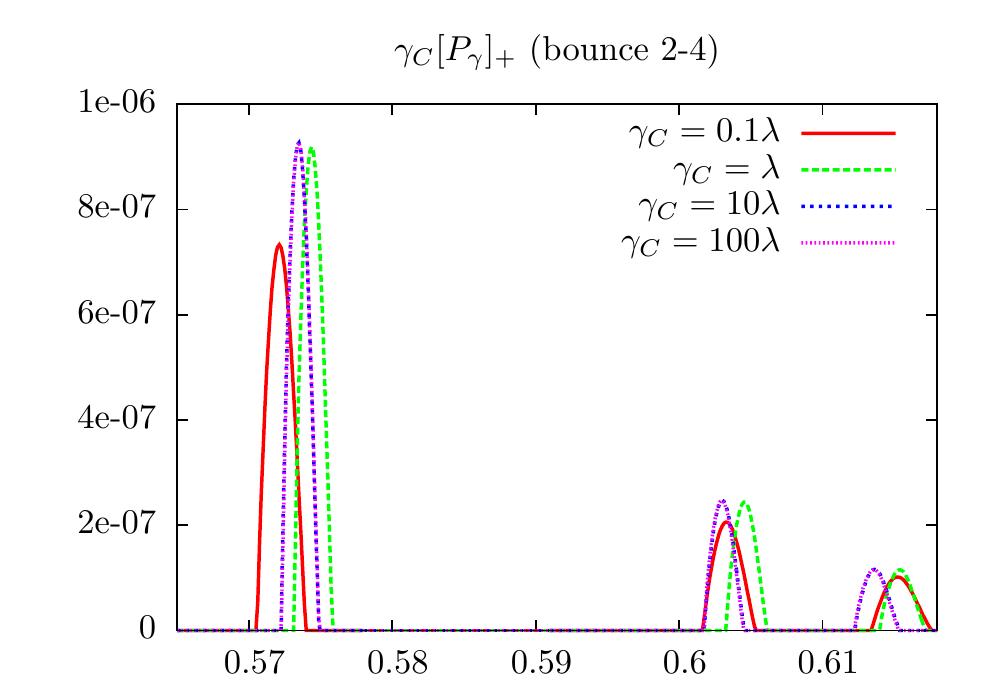}
\end{minipage}
%\hfil
%\begin{minipage}{0.5\textwidth}
%\includegraphics[width=1.1\textwidth]{mindist_zoom_epsp.pdf}
%\end{minipage}
\caption{\textit{Top}: Minimal distance $d_{\rm min}$ to the ground over time for different contact parameters $\gamma_c$. \textit{Bottom}: Contact force $\gamma_c [P_{\gamma_c}]_+$ over time for the first bounce (\textit{left}) and the second to fourth bounce (\textit{right}).\label{fig.mindist_gC}}
\end{figure}

In Figure~\ref{fig.mindist_gC}, we illustrate the influence of the contact parameter $\gamma_c$. As one would expect the violation of the relaxed contact condition is larger for a smaller $\gamma_c$, see the plot on the top right. For $\gamma_c\geq 10\lambda_s$ the curves are almost identical. 

On the bottom, we show the contact force $\gamma_c[P_{\gamma_c}]_+$, which appears on the right-hand side of \eqref{sigmaPgamma}, for the first four bounces. We see that the values of the force are almost independent of the chosen contact parameter. While for the smallest contact parameter $\gamma_c = 0.1\lambda_s$ the contact times are slightly altered, there is (almost) no visible difference between the results for $\gamma_c= 10\lambda_s$ and $100\lambda_s$.

\paragraph{Time discretisation}

\begin{figure}
\begin{minipage}{0.5\textwidth}
\includegraphics[width=1.1\textwidth]{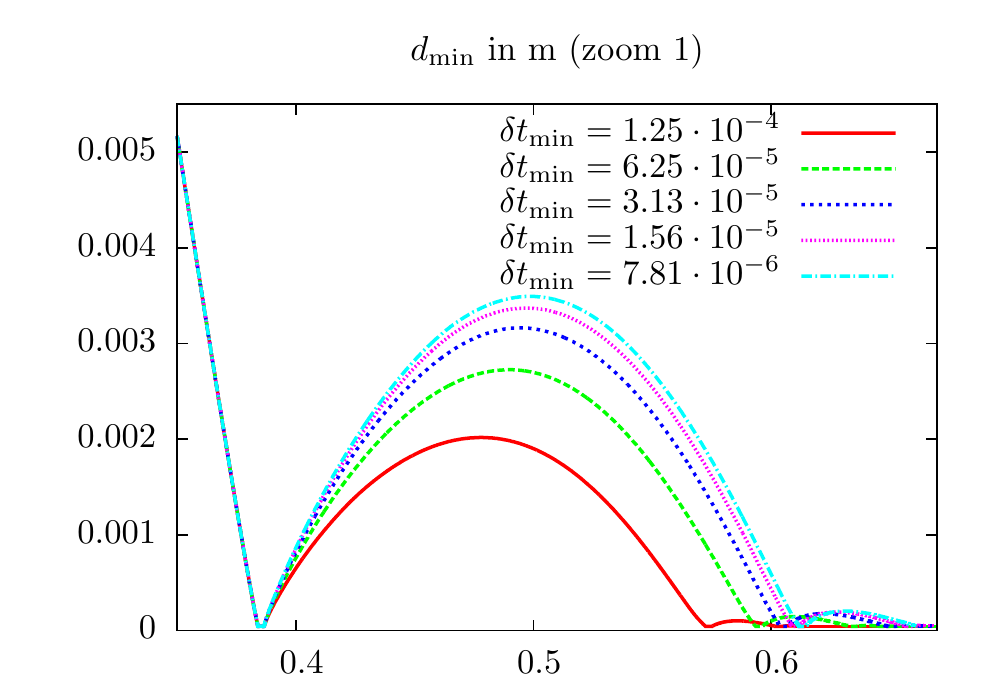}
\end{minipage}\hfil
\begin{minipage}{0.5\textwidth}
\includegraphics[width=1.1\textwidth]{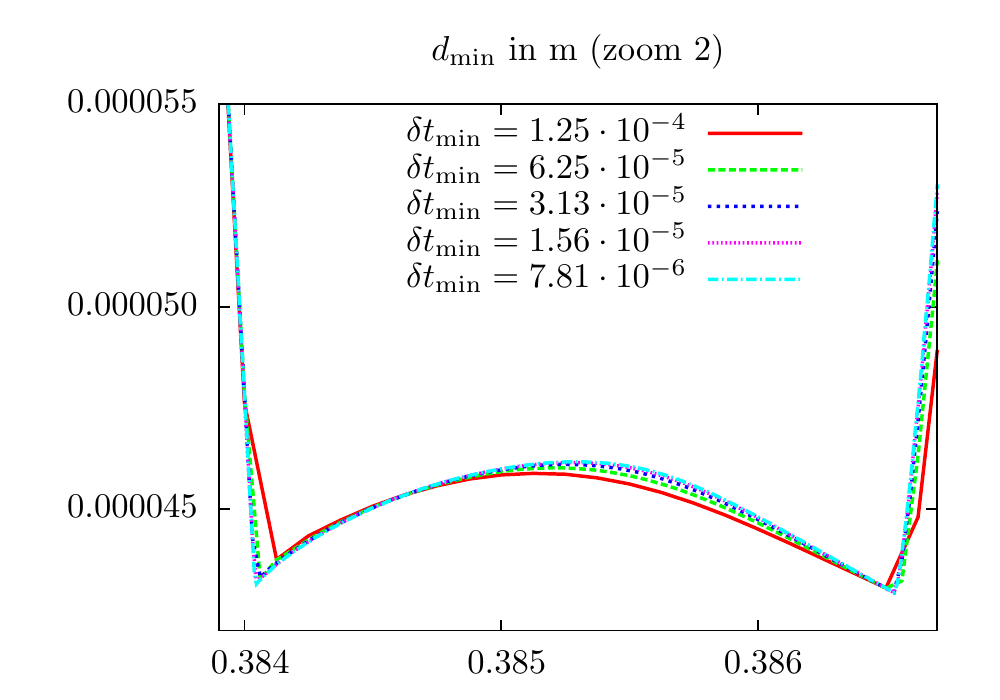}
\end{minipage}
\caption{Minimal distance $d_{\rm min}$ to the ground during the first two bounces (left) and at the first impact period (right)
for different minimal time-step sizes $\delta t_{\min}$.\label{fig.time}}
\end{figure}

In Figure~\ref{fig.time}, we investigate the influence of the time-step size $\delta t$ within and around the contact interval. In each simulation we start with a time-step of $\delta t=2\cdot 10^{-3}$ at $t=0$, which is decreased successively by a factor of 2 depending on the distance to the ground. We see that a very small time-step is necessary to capture the contact dynamics. While for the largest time-step $\delta t_{\min} = 1.25\cdot 10^{-4}$, the bounce is considerably reduced compared to the smaller time-step sizes, the curves seem to converge for $\delta t_{\min}\to 0$. The reason for the deviation can be deduced from the right plot, which shows that for $\delta t_{\min} = 1.25\cdot 10^{-4}$ the time of impact and release, where the curve shows a kink (i.e.\,the solution is non-smooth), is not captured accurately. 

\paragraph{Space discretisation}

\begin{figure}
\begin{minipage}{0.5\textwidth}
\includegraphics[width=1.1\textwidth]{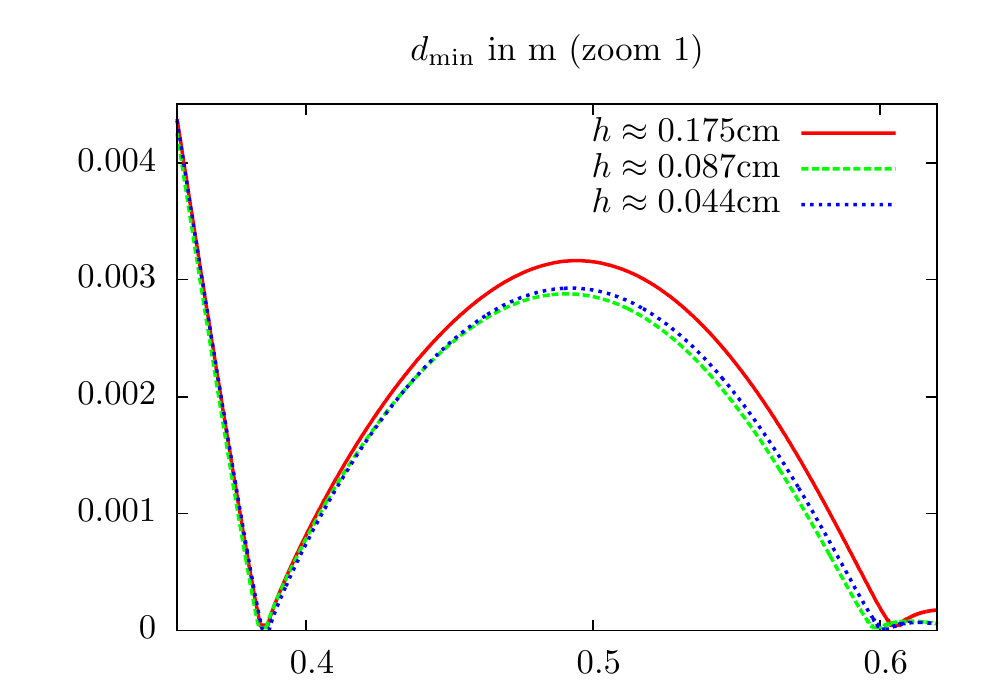}
\end{minipage}\hfil
\begin{minipage}{0.5\textwidth}
\includegraphics[width=1.1\textwidth]{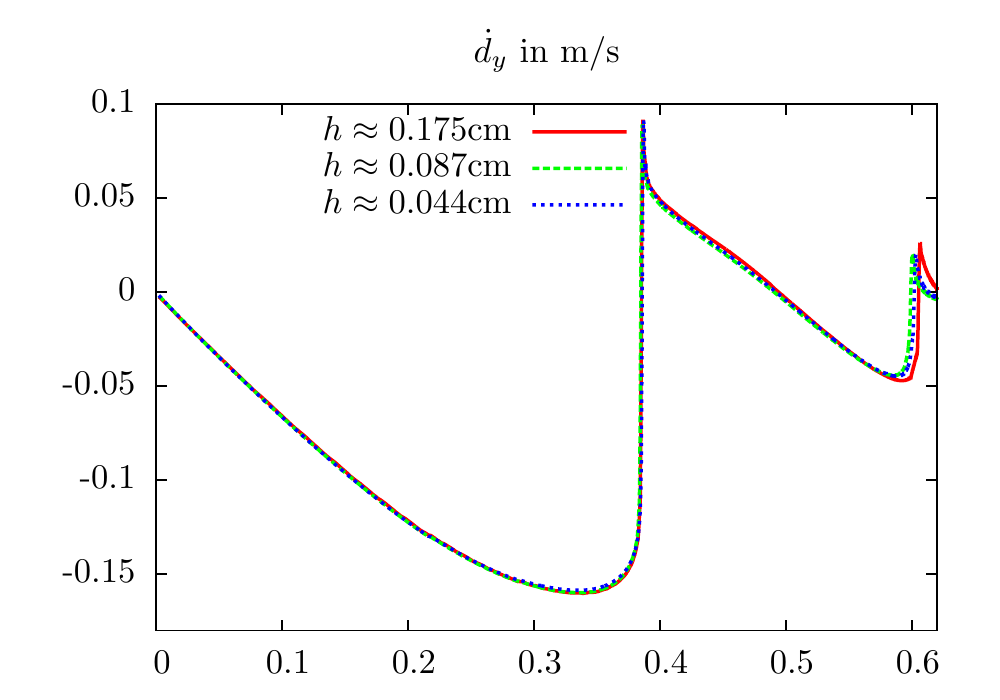}
\end{minipage}
\caption{Minimal distance $d_{\rm min}$ to the ground (left) and space-averaged solid velocity (right) over time on different mesh levels.
\label{fig.space}}
\end{figure}

In Figure~\ref{fig.space} we investigate the convergence behaviour under refinement of the finite element mesh. We fix the minimal time-step to $\delta t_{\min}= 3.125\cdot 10^{-5}$ and consider a coarse mesh with a maximum cell size of $h\approx 0.175$cm and 3 201 vertices and two finer meshes with 12 545 and 49 665 vertices that are constructed from the coarse mesh by global mesh refinement. We observe that the results both concerning minimal distance and vertical velocity are relatively close, even on the coarser mesh level, with an excellent agreement of the results on the finer meshes.

\paragraph{Comparison with a pure no-slip boundary condition}

\begin{figure}
\begin{minipage}{0.5\textwidth}
\includegraphics[width=1.1\textwidth]{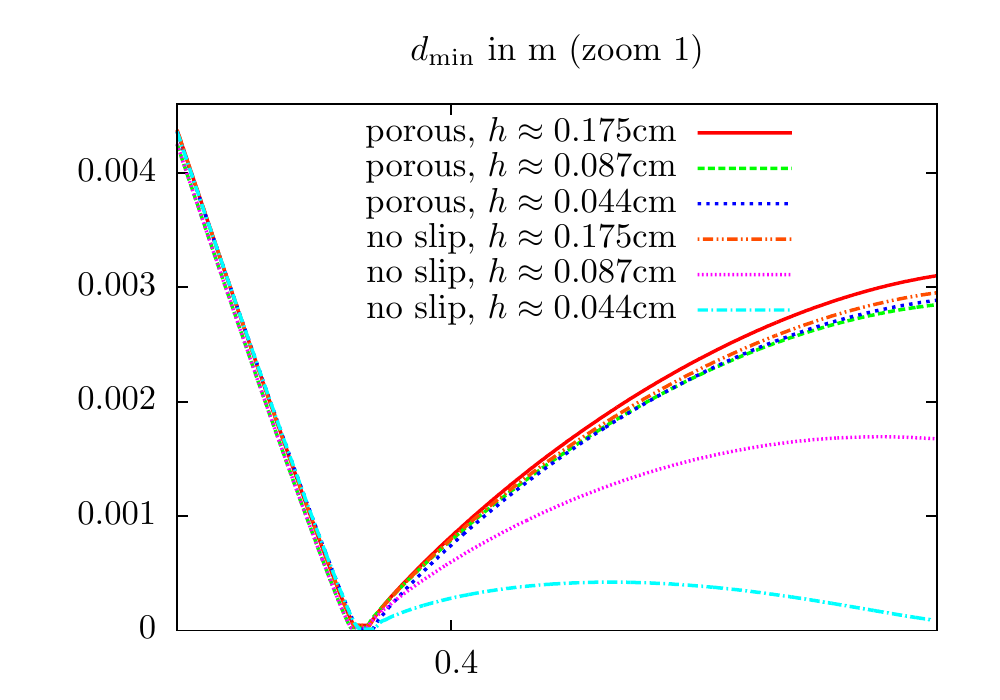}
\end{minipage}\hfil
\begin{minipage}{0.5\textwidth}
\includegraphics[width=1.1\textwidth]{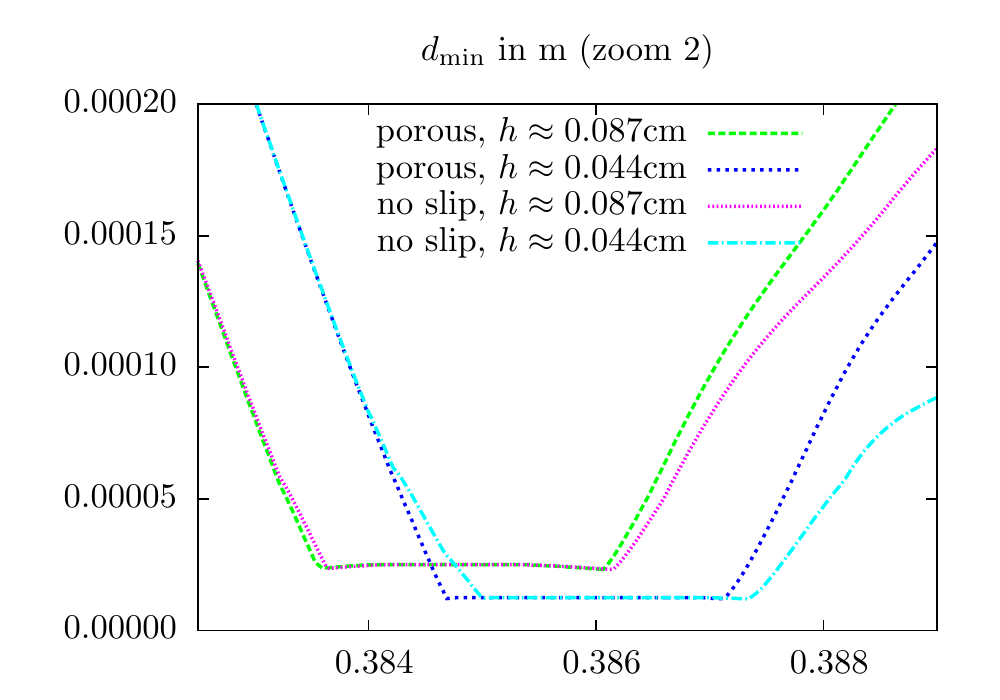}
\end{minipage}
\caption{Comparison of the contact approach with a porous layer with pure no-slip ("No Darcy") conditions.
Minimal distance $d_{\rm min}$ to the ground over time on different mesh levels and two different zooms.
\label{fig.noslipvsporous}}
\end{figure}

In Figure~\ref{fig.noslipvsporous}, we compare the approach presented in this paper with a simple relaxed contact approach without porous medium ("No Darcy"), where a no-slip condition (resp. a Navier-slip condition) is imposed for the fluid on the bottom wall $\Sigma_p$. The no-slip condition is the boundary condition, which is typically used for viscous fluids in absence of contact. First, we note that the curves for the no-slip condition and the Navier-slip condition with (small) slip-length $\left(\frac{\alpha}{\sqrt{K_\tau \epsilon_{\rm p}}}\right)^{-1}$ are almost identical. For this reason the latter curves are omitted in the following graphs.

As observed before, we see in the left picture 
that the curves obtained with the porous medium approach converge towards a certain bouncing height for $h\to 0$. Using a no-slip condition on $\Sigma_p$, the bounce get smaller and smaller and it is to be expected that for $h\to 0$ no bounce takes place at all (which is in agreement with the theoretical works on Navier-Stokes and contact~\cite{HeslaPhD, Hillairet2d, GerardVaretetal2015}).
The reason can be inferred from the zoom given on the right of Fig.~\ref{fig.noslipvsporous}, where we see that the fall is slowed down significantly right before the impact, while the curves for the two variants showed very good agreement until a distance of around $10^{-4}$ is reached. The reason are the strong fluid forces, in particular the pressure, that act against contact, when a pure no-slip condition is used. The finer the mesh, the better these forces are resolved. Interestingly, the results on the coarsest mesh ($h\approx 0.175$) show still a reasonable agreement, which might indicate that (only) on a very coarse mesh no-slip conditions could still yield physical results within a relaxed contact approach.

Finally, we show in Table~\ref{tab.v} the spatially-averaged velocity of the solid at the time of impact $t_i$ and the time of release $t_r$. Here we see quantitatively that the impact velocity is significantly reduced on the finer mesh levels when using a no-slip condition and thus, a much smaller rebound results. 
%The velocities at the time of release are around 20\% smaller compared to the impact velocities in all cases.

\begin{table}\centering
\begin{tabular}{c|cc|cc|cc}
\multicolumn{1}{c}{}&\multicolumn{2}{c}{\textbf{Porous}} &\multicolumn{2}{c}{\textbf{no-slip}} &\multicolumn{2}{c}{\textbf{Navier-slip}}\\
$h$ & $-\overline{\dot{d}_y}(t_i)$ &$\overline{\dot{d}_y}(t_r) $ &  $-\overline{\dot{d}_y}(t_i)$ &$\overline{\dot{d}_y}(t_r)$ &  $-\overline{\dot{d}_y}(t_i)$ &$\overline{\dot{d}_y}(t_r)$\\[0.1cm]
\hline\\[-0.3cm]
$1.75\cdot 10^{-3}$&$1.12\cdot 10^{-1}$ &$8.87\cdot 10^{-2}$ &$1.11\cdot 10^{-1}$ &$8.81\cdot 10^{-2}$ &$1.11\cdot 10^{-1}$ &$8.81\cdot 10^{-2}$
\\
$8.77\cdot 10^{-2}$ &$1.05\cdot 10^{-1}$ &$8.48\cdot 10^{-2}$ & $9.74\cdot 10^{-2}$ &$7.59\cdot 10^{-2}$ & $9.74\cdot 10^{-2}$ &$7.60\cdot 10^{-2}$ \\
$4.39\cdot 10^{-2}$ &$1.03\cdot 10^{-1}$ &$8.91\cdot 10^{-2}$ &$7.25\cdot 10^{-2}$ &$5.98\cdot 10^{-2}$ &$7.26\cdot 10^{-2}$ &$5.98\cdot 10^{-2}$
\end{tabular}
\caption{\label{tab.v} Spatially averaged velocity $\overline{\dot{d}_y}(t):= |\Omega_s(t)|^{-1} \int_{\Omega_s(t)} \dot{d}_y(t) \, dx$ of the ball at the time of impact $t=t_i$ and at the time of release $t=t_r$ for relaxed contact algorithms with a porous medium model and pure no-slip or Navier-slip conditions on $\Sigma_p$ on 3 different mesh levels.
}
\end{table}

\section{Conclusion}
\label{sec:concl}

We have introduced a physically consistent model to describe fluid-structure interactions with contact including seepage. For the latter a Darcy model is used on a thin porous layer of infinitesimal thickness. The approach can be used in a variety of different physical and numerical settings, including thick- and thin-walled solids, Eulerian or immersed (mixed-coordinate) descriptions, unfitted or fitted finite element discretisations, etc.

The numerical results show that the approach is numerically stable and (relatively) insensitive to variations of the numerical parameters, such as $\gamma_{\rm c}$. The model parameters $\epsilon_{\rm p}$ and $K$ of the porous layer need to be chosen depending on the application, e.g., the surface properties of the contacting bodies. Moreover, the results indicate convergence in both space and time. The time-step ${\delta t}_{\min}$ needs to be chosen very small in and around the contact interval to resolve the contact dynamics accurately.

Due to the relaxation of the contact conditions the approach is relatively easy to implement, in particular in comparison to approaches where a full topology change in the discrete fluid domain takes place and small numerical errors can lead to technical issues like unphysical "islands of fluid" appearing within the contact area, see~\cite{Ageretal2020}. 

Future work might focus on the extension to contact between multiple elastic bodies and on further developments within the time discretisation schemes on moving (sub-)domains, including for example adaptive strategies for the time steps $\delta t$.    

\section*{Acknowledgments}
EB was partially supported by the EPSRC grants EP/P01576X/1 and EP/T033126/1.
The second and fourth authors were partially supported by the French National Research Agency (ANR), trough the SIMR project (ANR-19-CE45-0020). 

%%%%%%%%%%%%%%%%%%%%

\bibliographystyle{apalike}
%\bibliography{lit}

\end{document}